\documentstyle[11pt,amssymb]{article}

\begin{document}

\begin{center}
{\large {\bf Arithmetic geometry of algebraic curves and their moduli space}} 
\end{center}
\begin{center}
Takashi Ichikawa (Saga University) 
\end{center}

\begin{center}
\underline{\large {\bf Abstract}} 
\end{center}
We review the following subjects: 
\begin{itemize}

\item 
Basic theory on algebraic curves and their moduli space, 

\item 
Schottky uniformization theory of Riemann surfaces, 
and its extension called arithmetic uniformization theory, 

\item 
Application to these theories to the arithmetic of the moduli space of algebraic curves,  especially to automorphic forms on this space. 

\end{itemize}
\vspace{1ex}
\begin{center}
\underline{\large {\bf Contents}} 
\end{center}

\S 1. Introduction (p.2--3) \\ 
\hspace{5ex} 1.1. Brief history \\ 
\hspace{5ex} 1.2. Plan of this lecture 

\S 2. Algebraic curves and Riemann surfaces (p.4--11) \\ 
\hspace{5ex} 2.1. Riemann's correspondence \\ 
\hspace{5ex} 2.2. Riemann-Roch's theorem \\ 
\hspace{5ex} 2.3. Differential forms, periods and Jacobians 

\S 3. Schottky uniformization (p.12--17) \\ 
\hspace{5ex} 3.1. Degeneration of Riemann surfaces \\ 
\hspace{5ex} 3.2. Schottky uniformization of Riemann surfaces \\ 
\hspace{5ex} 3.3. Explicit formula of periods \\ 
\hspace{5ex} 3.4. Fractal nature of Schottky groups 

\S 4. Arithmetic uniformization (p.18--26) \\ 
\hspace{5ex} 4.1. Periods as power series \\ 
\hspace{5ex} 4.2. Tate curve and Mumford curves \\ 
\hspace{5ex} 4.3. Arithmetic Schottky uniformization 

\S 5. Moduli space of algebraic curves (p.27--35) \\ 
\hspace{5ex} 5.1. Construction of moduli spaces \\ 
\hspace{5ex} 5.2. Stable curves and their moduli space \\ 
\hspace{5ex} 5.3. Intersection theory on the moduli space 

\S 6. Arithmetic theory of modular forms (p.36--48) \\ 
\hspace{5ex} 6.1. Elliptic modular forms \\ 
\hspace{5ex} 6.2. Siegel modular forms \\ 
\hspace{5ex} 6.3. Teichm\"{u}ller modular forms

\newpage

\begin{center}
\underline{\large {\bf \S 1. Introduction}} 
\end{center}
\noindent
\underline{\bf 1.1. Brief history} 
\begin{itemize}

\item 
Around 1800～1830, 
Gauss(1777--1855), Abel(1802--1829) and Jacobi(1804--1851) found 
that the inverse function of the elliptic integral: 
$$
y \ = \ \int dx \left/ \sqrt{f(x)} \right. \ 
\mbox{$(f(x) :$ a polynomial of degree 4 without multiple root)}
$$
is an elliptic function, i.e., a double periodic function of 
the complex variable $y,$ 
and they expressed the function as an infinite product 
and the ratios of theta functions. 
\\
$\Rightarrow$ complex function theory. 

\item 
Riemann(1826--1866) constructed Riemann surfaces 
from algebraic function fields 
$$
{\bf C}(x, y) \ \ 
\mbox{($x$ : a variable, $y$ : finite over ${\bf C}(x)$)}, 
$$
and solved Jacobi's inverse problem 
using Abel-Jacobi's theorem and Riemann's theta functions. 
\\
$\Rightarrow$ complex geometry and algebraic geometry (1857). 

\item Teichm\"{u}ller(1913--1943) constructed analytic theory 
on the moduli of Riemann surfaces. 

\item Mumford constructed the moduli of algebraic curves as 
an algebraic variety (1956), and studied this geometry. 
Further, he and Deligne gave its compactification 
as the moduli of stable curves (1969). 

\item String theory provided a strong relationship between physics 
and the theory of moduli of curves. 

\item Around 1960～1970, Shimura constructed arithmetic theory 
on Shimura models with applications 
to the rationality on Siegel modular forms, 
and further Chai and Faltings extended his result to any base ring (1990). 

\end{itemize}

\newpage
\noindent
\underline{\bf 1.2. Plan of this lecture}
\vspace{2ex}

We will review the following subjects with some proof: 
\begin{itemize}

\item Very classical results on algebraic curves over ${\bf C}$ and 
the associated Riemann surfaces: 
for example, $\wp$-functions and elliptic curves, 
differential $1$-forms and period integrals, 
Riemann-Roch's theorem, Abel-Jacobi's theorem and Jacobian varieties, 
degeneration, Schottky uniformization 
and the description of forms and periods.  

\item 
Rather classical results on moduli and families of algebraic curves: 
for example, moduli of elliptic curves and higher genus curves, 
stable curves and their moduli (Deligne-Mumford's compactification), 
the irreducibility of the moduli, Eisenstein series and Tate curve, 
Mumford curves; 
\vspace{1ex}

\noindent
and recent results on arithmetic version of Schottky uniformization. 

\item 
Recent results on arithmetic geometry of the moduli space of algebraic curves: 
for example, Fourier expansion of (elliptic and Siegel) modular forms 
and their rationality, Teichm\"{u}ller modular forms and the Schottky problem, 
Mumford's isomorphism. 

\end{itemize}

Especially, 
we explain that the classical, but not so familiar Schottky uniformization theory 
which gives an explicit description of differential forms, 
periods and degeneration of Riemann surfaces. 
Furthermore, 
we give its extended version in the category of arithmetic geometry 
(unifying complex geometry and formal geometry over ${\bf Z})$, 
and the application to automorphic forms, 
called Teichm\"{u}ller modular forms,  
on the moduli space of algebraic curves.

\newpage
\begin{center}
\underline{\large {\bf \S 2. Algebraic curves and Riemann surfaces}} 
\end{center}

\noindent
\underline{\bf 2.1. Riemann's correspondence}
\vspace{2ex}

\noindent
\underline{\bf Algebraic curves.}
{\bf Algebraic varieties} are topological spaces obtained by gluing 
zero sets of polynomials of multiple variables, 
and closed subsets of algebraic varieties are defined as 
zero sets of polynomials (Zariski topology). 
These examples are 
\begin{eqnarray*}
\mbox{the {\bf projective {\boldmath $n$}-space}} \ \ 
{\mathbb P}^{n}_{k} 
& \stackrel{\rm def}{=} & 
( k^{n+1} - \{ (0,...,0) \} ) / k^{\times} \\ 
& = & 
\{ (x_{0}: \cdots : x_{n+1}) = ( c x_{0}: \cdots : c x_{n+1}) \ | \ 
c \neq 0 \}, 
\end{eqnarray*}
and its subsets (called {\bf projective varieties} which are {\bf proper} 
over $k$ $(\doteq$ compact)) as the zero sets 
of homogeneous polynomials over an algebraically closed field $k.$ 
\begin{center}
(algebraic) {\bf curves} $\stackrel{\rm def}{=}$ 
1-dimensional algebraic varieties 
\end{center}

\noindent
\underline{\bf Riemann's correspondence.} 
There exists an equivalence (trinity) of the categories: 
$$
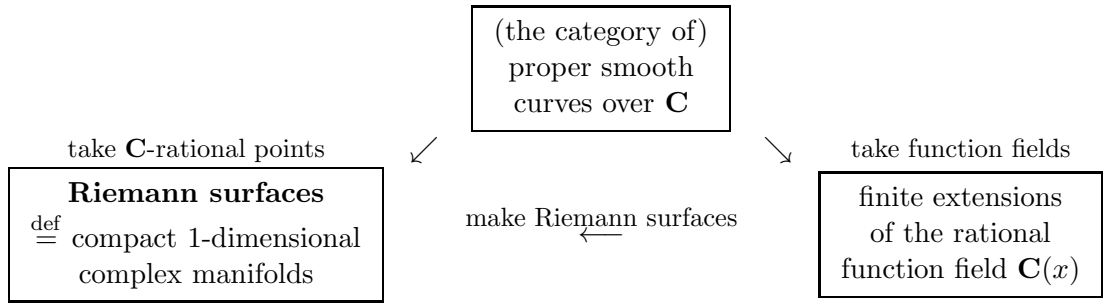
\begin{array}{ccccc}
& & 
\fbox{$
\begin{array}{c} \mbox{(the category of)} \\ 
\mbox{proper smooth} \\ 
\mbox{curves over {\bf C}} 
\end{array}
$}
& \\ 
\mbox{\small take {\bf C}-rational points} & \swarrow & & \searrow & 
\mbox{\small take function fields} \\ 
\fbox{$
\begin{array}{c} 
\mbox{{\bf Riemann surfaces}} \\ 
\mbox{$\stackrel{\rm def}{=}$ compact 1-dimensional} \\ 
\mbox{complex manifolds} 
\end{array}
$}
&  & \stackrel{\mbox{\small make Riemann surfaces}}{\longleftarrow} &  & 
\fbox{$
\begin{array}{c} 
\mbox{finite extensions} \\ 
\mbox{of the rational} \\ 
\mbox{function field ${\bf C}(x)$} 
\end{array}
$}
\end{array} 
$$
\noindent
\underline{\bf Genus.} 
The {\bf genus} of a Riemann surface and the corresponding curve is 
defined as the number of its holes ({\bf Figure}). 

\noindent
\underline{\bf Genus 0 case.} 
$$
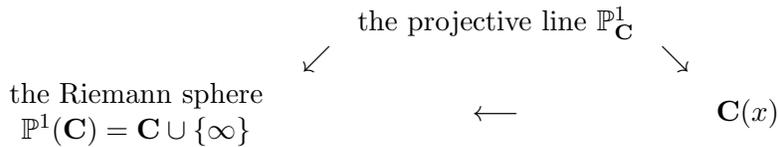
\begin{array}{ccccc}
& & \mbox{the projective line ${\mathbb P}^{1}_{\bf C}$} & \\ 
& \swarrow & & \searrow \\ 
\begin{array}{c} 
\mbox{the Riemann sphere} \\ 
{\mathbb P}^{1}({\bf C}) = {\bf C} \cup \{ \infty \} 
\end{array}
&  & \longleftarrow &  & 
{\bf C}(x) 
\end{array} 
$$
\noindent
\underline{\bf Genus 1 case.} 
For cubic polynomials $f(x) \in {\bf C}[x]$ without multiple root, 
$$
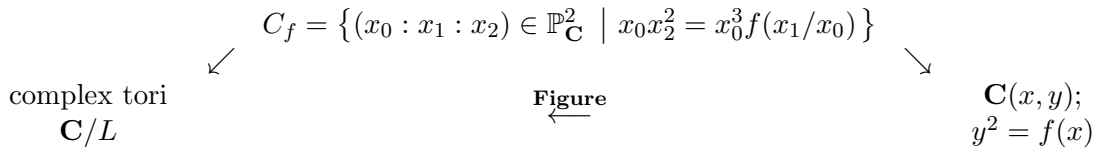
\begin{array}{ccccc}
& & 
C_{f} = 
\left\{ (x_{0}: x_{1}: x_{2}) \in {\mathbb P}^{2}_{\bf C} \ \left| \ 
x_{0} x_{2}^{2} = x_{0}^{3} f(x_{1}/x_{0}) \right. 
\right\} 
& \\ 
  & \swarrow & & \searrow \\ 
\begin{array}{c} 
\mbox{complex tori} \\ 
{\bf C} / L 
\end{array}
&  & \stackrel{\bf Figure}{\longleftarrow} &  & 
\begin{array}{c}
{\bf C}(x, y); \\ 
y^{2} = f(x) 
\end{array}
\end{array} 
$$
Here 
\begin{eqnarray*}
& & 
\mbox{$L$ is a {\bf lattice} in {\bf C}, i.e., a sub {\bf Z}-module of rank 2 
such that $L \otimes_{\bf Z} {\bf R} = {\bf C},$} 
\\ 
& & 
E_{2k}(L) \stackrel{\rm def}{=} \sum_{u \in L-\{0 \}} \frac{1}{u^{2k}} : 
\ \mbox{absolutely convergent series for $k > 1,$} 
\\ 
& & 
f(x) \stackrel{\rm def}{=} 4 x^{3} - 60 E_{4}(L) x - 140 E_{6}(L), 
\\
& & \wp(z) = \wp_{L}(z) \stackrel{\rm def}{=} \frac{1}{z^{2}} + 
\sum_{u \in L-\{0\}} \left( \frac{1}{(z-u)^{2}} - \frac{1}{u^{2}} \right): \ 
\mbox{Weierstrass' {\bf {\boldmath $\wp$}-function}} \\ 
& \Rightarrow & 
\left\{ \begin{array}{l} 
\mbox{$\wp(z)$ is absolutely and uniformly convergent on any compact subset of 
${\bf C} - L,$} 
\\
\mbox{$z \mapsto (1 : \wp(z) : \wp'(z))$ gives a biholomorphic map 
${\bf C}/L \stackrel{\sim}{\rightarrow} C_{f}({\bf C}),$} 
\end{array} \right. 
\end{eqnarray*}
and 
\begin{eqnarray*}
& & \mbox{${\bf C}(x, y)$ is a quadratic extension of ${\bf C}(x)$ 
defined by $y^{2} = f(x)$} \\ 
& \leftrightarrow & 
\mbox{$C_{f}$ is a double cover of ${\mathbb P}^{1}_{\bf C}$ 
ramified at the 3 roots of $f(x)$ and $\infty.$} 
\end{eqnarray*}
An {\bf elliptic curve} is a proper smooth curve $C$ of genus 1 and 
with one marked point $x_{0}$. 
Then $C$ has unique commutative group structure defined algebraically 
with origin $x_{0}.$ 
For example, the above $C_{f}$ with $(0:0:1)$ is an elliptic curve, and 
the map ${\bf C}/L \stackrel{\sim}{\rightarrow} C_{f}({\bf C})$ 
is also a group isomorphism which follows from the addition law of $\wp(z):$ 
$$
\wp(z + w) \ = \ - \wp(z) - \wp(w) 
+ \frac{1}{4} 
\left( \frac{\wp'(z) - \wp'(w)}{\wp(z) - \wp(w)} \right)^{2}.
$$
\underline{\bf Exercise 1.} 
Show the following Laurent expansion of $\wp(z):$ 
$$
\wp(z) = \frac{1}{z^{2}} + \sum_{n=1}^{\infty} (2n+1) E_{2n+2}(L) z^{2n} 
\ \ \mbox{around $z = 0.$} 
$$
Further, using this fact, the periodicity of $\wp(z):$ 
$$
\wp(z + u) = \wp(z) \ \ (u \in L),
$$
and the maximum principle on holomorphic functions, 
prove that 
\begin{eqnarray*}
& & \wp'(z)^{2} = 4 \wp(z)^{3} - 60 E_{4}(L) \wp(z) - 140 E_{6}(L)
\\
& & \left( \mbox{i.e.,} \ \wp(z) = x \ \Rightarrow \ 
z = \int \frac{dx}{\sqrt{4 x^{3} - 60 E_{4} x - 140 E_{6}}} \right) 
\end{eqnarray*}
and that ${\displaystyle E_{8}(L) = \frac{3}{7} E_{4}(L)^{2}.}$ 
\vspace{2ex}

\noindent
\underline{\bf Genus $> 1$ case.} 
For a Riemann surface $R$ of genus $> 1,$ 
by Riemann's mapping theorem, 
its universal cover is biholomorphic to 
$$
H_{1} \ \stackrel{\rm def}{=} \ \left\{ \tau \in {\bf C} \ \left| 
\ {\rm Im}(\tau) > 0 \right. \right\}: 
\mbox{the {\bf Poincar\'{e} upper half plane}.} 
$$
Then we have 
$$
R \cong \ H_{1} / \pi_{1}(R) : \ \mbox{called a {\bf Fuchsian model}}, 
$$
where the fundamental group $\pi_{1}(R)$ of $R$ 
is a cocompact discrete subgroup of 
$$
PSL_{2}({\bf R}) \ \stackrel{\rm def}{=} \ 
\left. \left\{ \left. 
\left( \begin{array}{cc} a & b \\ c & d \end{array} \right) \in M_{2}({\bf R}) 
\ \right| \ \ ad - bc = 1 \right\} \right/ \left\{ \pm E_{2} \right\}
$$ 
which acts on $H_{1}$ by the M\"{o}bius transformation: 
$$
\tau \mapsto \frac{a \tau + b}{c \tau + d} 
$$
(in fact, $PSL_{2}({\bf R})$ is the group ${\rm Aut}(H_{1})$ 
of complex analytic automorphisms of $H_{1}).$ 
\vspace{2ex}

\noindent
\underline{\bf Remark.} 
Let $\Gamma$ be a congruence subgroup of $SL_{2}({\bf Z}),$ 
for example 
\begin{eqnarray*}
SL_{2}({\bf Z}) & \stackrel{\rm def}{=} & 
\left\{ \left. 
\left( \begin{array}{cc} a & b \\ c & d \end{array} \right) \in M_{2}({\bf Z}) 
\ \ \right| \ \ ad - bc = 1 \right\} 
\\
\supset \ 
\Gamma_{0}(N) & \stackrel{\rm def}{=} & \left\{ \left. 
\left( \begin{array}{cc} a & b \\ c & d \end{array} \right) 
\in SL_{2}({\bf Z}) \ \right| \ \ c \equiv 0 \ {\rm mod}(N) \right\} 
\\
\supset \ \ 
\Gamma(N) & \stackrel{\rm def}{=} & \left\{ \left. 
\left( \begin{array}{cc} a & b \\ c & d \end{array} \right) 
\in SL_{2}({\bf Z}) \ \right| \ 
\ a-1 \equiv b \equiv c \equiv d-1 \equiv 0 \ {\rm mod}(N) \right\} 
\\
& : & \mbox{the principal congruence subgroup of $SL_{2}({\bf Z})$ 
of level $N.$} 
\end{eqnarray*}
Then $H_{1} / \Gamma$ is a {\it noncompact} 1-dimensional complex manifold, 
and becomes compact by adding the set ${\mathbb P}^{1}({\bf Q})/\Gamma$ of 
{\bf cusps} of $\Gamma.$ 
$H_{1} / \Gamma$ and 
$\left. \left(H_{1} \cup {\mathbb P}^{1}({\bf Q}) \right) \right/ \Gamma$ are 
called {\bf modular curves}. 
\vspace{2ex}

\noindent
\underline{\bf 2.2. Riemann-Roch's theorem} 
\vspace{2ex}

Let $R$ be a Riemann surface of genus $g$, 
and let $D$ be a divisor on $R$ which is, by definition, 
a finite sum of points on $R$ with coefficients in ${\bf Z}$. 
When $D$ is represented as $\sum_{P \in R} a_{P} \cdot P$, 
the associated invertible sheaf ${\cal O}_{R}(D)$, 
namely line bundle, on $R$ is defined as 
$$
{\cal O}_{R}(D)(U) \stackrel{\rm def}{=}   
\left\{ \mbox{$f$: meromorphic functions on $U$} \ | \ 
{\rm ord}_{P}(f) + a_{P} \geqq 0 \ (P \in U) \right\} 
$$
for open subsets $U$ of $R$. 
Then Riemann-Roch's theorem states 
$$ 
\dim_{\bf C} H^{0}(R, {\cal O}_{R}(D)) - 
\dim_{\bf C} H^{1}(R, {\cal O}_{R}(D)) = \deg(D) + 1 - g, 
$$
where $\deg(D) = \sum_{P \in R} a_{P}$ is the degree of $D$. 
Denote by $\Omega_{R}$ the  invertible sheaf of holomorphic $1$-forms on $R$. 
Then by Serre's duality, 
the residue map gives a nondegenerate pairing 
$$
H^{0}(R, \Omega_{R}(-D)) \times H^{1}(R,  {\cal O}_{R}(D)) \rightarrow {\bf C}, 
$$ 
and hence 
$$ 
\dim_{\bf C} H^{0}(R, {\cal O}_{R}(D)) - 
\dim_{\bf C} H^{0}(R, \Omega_{R}(-D)) = \deg(D) + 1 - g. 
$$

\noindent
\underline{\bf 2.3. Differential forms, periods and Jacobians} 
\vspace{2ex}

Let $R$ be a Riemann surface of genus $g \geq 1.$ 
Then its fundamental group $\pi_{1}(R; x_{0})$ with base point $x_{0} \in R$ 
is represented by 
$$
\left\langle \left. 
\underbrace{\alpha_{1}, \beta_{1},..., \alpha_{g}, \beta_{g}}_{\rm generators} 
\ \right| \ \underbrace{\prod_{i=1}^{g} \left( \alpha_{i} \beta_{i} 
\alpha_{i}^{-1} \beta_{i}^{-1} \right) = 1}_{\rm relation} 
\right\rangle, 
$$
where the generators $\alpha_{i}, \beta_{i}$ are {\bf canonical}, i.e., 
closed oriented paths in $R$ with base point $x_{0}$ 
such that $\alpha_{i}, \beta_{i}$ intersect as the $x, y$-axes and that 
$(\alpha_{i} \cup \beta_{i}) \cap (\alpha_{j} \cup \beta_{j}) = \{ x_{0} \}$ 
if $i \neq j$ ({\bf Figure}). 
Then 
\vspace{2ex}

\noindent
\underline{\bf Theorem 2.1.} (Abel, Jacobi, Riemann, see [Mur]) \begin{it}

{\rm (1)} The space $H^{0} \left(R, \Omega_{R} \right)$ of 
{\it \bfseries holomorphic {\boldmath $1$}-forms} on $R$ is $g$-dimensional, 
and is generated by unique holomorphic $1$-forms $\omega_{1},..., \omega_{g}$ 
satisfying that ${\displaystyle \int_{\alpha_{i}} \omega_{j} = \delta_{ij}.}$ 
Furthermore, 
$\deg \left( \Omega_{R} \right) = 2g - 2$. 

{\rm (2) (Riemann's period relation)} 
The {\it \bfseries period matrix} 
$$
Z \stackrel{\rm def}{=} 
\left( \int_{\beta_{i}} \omega_{j} \right)_{1 \leq i,j \leq g}
$$ 
of $(R; \{ \alpha_{i}, \beta_{i} \}_{1 \leq i \leq g})$ is symmetric, 
and its imaginary part ${\rm Im}(Z)$ is positive definite. 

{\rm (3) (Abel-Jacobi's theorem)} 
Let 
$$
{\rm Cl}^{0}(R) \ = \ 
\left\{ \mbox{divisors with degree $0$ on $R$} \right\} \left/ 
\left\{ \sum_{P \in R} {\rm ord}_{P}(f) \cdot P \right\} \right. 
$$
the {\it \bfseries divisor class group} with degree $0$ of $R,$ 
and let ${\bf C}^{g} / L$ be the $g$-dimensional complex torus obtained from 
the lattice $L \stackrel{\rm def}{=} {\bf Z}^{g} + {\bf Z}^{g} \cdot Z$ 
in ${\bf C}^{g}.$ 
Then the map 
$$
\sum_{j} (P_{j} - Q_{j}) \mapsto 
\left( \sum_{j} \int_{Q_{j}}^{P_{j}} \omega_{i} \right)_{1 \leq i \leq g} 
$$
gives rise to a group isomorphism: 
$$
\mu : {\rm Cl}^{0}(R) \ \stackrel{\sim}{\longrightarrow} \ {\bf C}^{g} / L. 
$$ 
\end{it}

\noindent
\underline{\bf Remark.}
It is clear that $(z_{1},..., z_{g}) \mapsto 
\left( \exp(2 \pi \sqrt{-1} z_{1}),..., \exp(2 \pi \sqrt{-1} z_{g}) \right)$ 
gives the isomorphism 
$$
{\bf C}^{g} / L \ \stackrel{\sim}{\rightarrow} \ 
({\bf C}^{\times})^{g} \left/ \left\langle \left. 
\left( \exp \left( 2 \pi \sqrt{-1} \int_{\beta_{i}} \omega_{j} \right) 
\right)_{1 \leq i \leq g} \ \right| \ 1 \leq j \leq g \right\rangle, \right. 
$$
and then 
$$
\exp \left( 2 \pi \sqrt{-1} \int_{\beta_{i}} \omega_{j} \right) 
\ \ (1 \leq i, j \leq g) 
$$
are called the {\bf multiplicative periods}. 
Let ${\rm Pic}^{0}(R)$ denote the {\bf Picard group} with degree $0$ of $R$ 
which is defined as the group of linear equivalence classes of line bundles 
with degree $0$ over $R.$ 
Then it is known that 
\begin{eqnarray*}
{\rm Cl}^{0}(R) & \cong & {\rm Pic}^{0}(R) 
\\
D & \leftrightarrow & {\cal O}_{R}(D) 
\end{eqnarray*}
becomes an {\bf abelian variety}, i.e., a proper (commutative) algebraic group 
over ${\bf C},$ 
and the isomorphism is also a biholomorphic map. 
This abelian variety is called the {\bf Jacobian variety} of $R$ 
(or of the associated curve), and denoted by ${\rm Jac}(R).$ 
\vspace{2ex}

\underline{\it Proof.} 
(1) By Riemann-Roch's theorem, 
\begin{eqnarray*}
\dim_{\bf C} H^{0}(R, {\cal O}_{R}) - 
\dim_{\bf C} H^{0}(R, \Omega_{R}) 
& = & 
1 - g, 
\\
\dim_{\bf C} H^{0}(R, \Omega_{R}) - 
\dim_{\bf C} H^{0}(R, {\cal O}_{R}) 
& = & 
\deg(\Omega_{R}) + 1 - g. 
\end{eqnarray*}
By the maximum principle on holomorphic functions on $R$, 
$H^{0}(R, {\cal O}_{R})$ consists of constant functions on $R$, 
and hence 
$$
\dim_{\bf C} H^{0}(R, \Omega_{R}) = g, 
\ \deg(\Omega_{R}) = 2g - 2. 
$$

To prove the remains of (1), and (2), (3), 
first we show a generalized form of Riemann's period relation. 
Let ${\cal P}$ the $4g$ oriented sided polygon obtained from $R$ 
by cutting the paths $\alpha_{i}, \beta_{i}$ $(1 \leq i \leq g)$ 
({\bf Figure}). 
Fix $P_{0} \in {\cal P},$ 
and for a holomorphic $1$-form $\phi$ on $R,$ define 
$f(P) = \int_{P_{0}}^{P} \phi.$ 
Then for a meromorphic $1$-form $\psi$ on $R$ whose poles belong to 
the interior ${\cal P}^{\circ}$ of ${\cal P}$ 
(this condition is satisfied by moving slightly $\alpha_{i}, \beta_{i}$ 
if necessary), 
using the function $f^{\pm}$ on the boundary $\partial {\cal P}$ of ${\cal P}$ 
defined by 
\begin{eqnarray*}
f^{+}(P) & \stackrel{\rm def}{=} & 
\int_{P_{0}}^{P} \phi \ \ (P \in \alpha_{i} \cup \beta_{i}), 
\\ 
f^{-}(P) & \stackrel{\rm def}{=} & 
\int_{P_{0}}^{P} \phi \ \ (P \in -\alpha_{i} \cup -\beta_{i}), 
\end{eqnarray*}
we have 
\begin{eqnarray*}
2 \pi \sqrt{-1} \sum_{P \in {\cal P}} {\rm Res}_{P}(f \psi) & = & 
\int_{\partial {\cal P}} f \psi \ \ \mbox{(by the residue theorem)} 
\\ 
& = & 
\sum_{i=1}^{g} \left( \int_{\alpha_{i}} f^{+} \psi 
+ \int_{-\alpha_{i}} f^{-} \psi + 
\int_{\beta_{i}} f^{+} \psi + \int_{-\beta_{i}} f^{-} \psi \right) 
\\
& = & 
\sum_{i=1}^{g} \left( \int_{\alpha_{i}} (f^{+} - f^{-}) \psi + 
\int_{\beta_{i}} (f^{+} - f^{-}) \psi \right) 
\\
& = & 
\sum_{i=1}^{g} \left( \left( - \int_{\beta_{i}} \phi \right) 
\left( \int_{\alpha_{i}} \psi \right) + \left( \int_{\alpha_{i}} \phi \right) 
\left( \int_{\beta_{i}} \psi \right) \right). 
\end{eqnarray*}
Therefore, 
$$
2 \pi \sqrt{-1} \sum_{P \in {\cal P}} {\rm Res}_{P}(f \psi) \ = \ 
\sum_{i=1}^{g} \left( \left( \int_{\alpha_{i}} \phi \right) 
\left( \int_{\beta_{i}} \psi \right) - \left( \int_{\beta_{i}} \phi \right) 
\left( \int_{\alpha_{i}} \psi \right) \right) 
$$
which we call the {\bf generalized Riemann's period relation}. 

In particular, for two holomorphic $1$-forms $\varphi, \varphi',$ 
put $f(P) = \int_{P_{0}}^{P} \varphi$ $(P \in {\cal P}^{\circ}),$ 
and put 
$$
A_{i} = \int_{\alpha_{i}} \varphi , \ \ A'_{i} = \int_{\alpha_{i}} \varphi' , 
\ \ B_{i} = \int_{\beta_{i}} \varphi , \ \ B'_{i} = \int_{\beta_{i}} \varphi' . $$
Then by the above, 
$$
\sum_{i=1}^{g} \left( A_{i} B'_{i} - B_{i} A'_{i} \right) \ = 0. 
$$ 
Further, 
\begin{eqnarray*}
{\rm Im} \left( \sum_{i=1}^{g} \overline{A_{i}} B_{i} \right) 
& = & 
\frac{1}{2 \sqrt{-1}} 
\sum_{i=1}^{g} \left( \overline{A_{i}} B_{i} - \overline{B_{i}} A_{i} \right) 
\\ 
& = & 
\frac{1}{2 \sqrt{-1}} \int_{\partial {\cal P}} \overline{f} \varphi 
= 
\frac{1}{2 \sqrt{-1}} \int_{\cal P} d \left( \overline{f} \varphi \right) 
\\ 
& = & 
\int_{\cal P} du dv 
= 
\int_{\cal P} \left( \left( \frac{\partial u}{\partial x} \right)^{2} + 
\left( \frac{\partial u}{\partial y} \right)^{2} \right) dx dy 
\\
& & 
\left( f = u + \sqrt{-1} v, \ 
z = x + \sqrt{-1} y : \mbox{local coordinates} \right) 
\end{eqnarray*}
by Cauchy-Riemann's relation 
$\partial u / \partial x = \partial v / \partial y$, 
$\partial u / \partial y = - \partial v / \partial x$, 
and hence this value is positive if $\varphi$ is not identically $0$. 
Therefore, any holomorphic $1$-form $\varphi$ with 
$\int_{\alpha_{i}} \varphi = 0$ $(1 \leq i \leq g)$ becomes identically $0,$ 
and hence for any base $\omega'_{i}$ $(1 \leq i \leq g)$ of 
$H^{0}(R, \Omega_{R}),$ 
$\left( \int_{\alpha_{i}} \omega'_{j} \right)_{i,j}$ is a regular matrix. 
This implies (1). 

(2) If $\varphi = \omega_{i}, \varphi' = \omega_{j},$ then by the above, 
$\int_{\beta_{i}} \omega_{j} - \int_{\beta_{j}} \omega_{i} = 0,$ 
hence $Z$ is symmetric. 
Further, 
if $\varphi = \sum_{i=1}^{g} c_{i} \omega_{i} \in H^{0}(R, \Omega_{R})$ 
is not $0,$ then 
$\mbox{\boldmath $c$} = (c_{1},..., c_{g}) \neq \mbox{\boldmath $0$}$, 
and hence  
${\rm Im}(\overline{\mbox{\boldmath $c$}} Z 
\mbox{}^{t} \mbox{\boldmath $c$}) > 0$. 
This implies that ${\rm Im}(Z)$ is positive definite. 
This proves (2).  

(3) If $f$ is a meromorphic function on $R,$ 
then by the generalized period relation, 
\begin{eqnarray*}
& & \sum_{P \in R} 
\left( {\rm ord}_{P}(f) \cdot \int_{P_{0}}^{P} \omega_{j} \right) 
\\
& = & 
\sum_{P \in {\cal P}} 
{\rm Res}_{P} \left( \int_{P_{0}}^{P} \omega_{j} \cdot \frac{df}{f} \right) 
\\
& = & 
\frac{1}{2 \pi \sqrt{-1}} \sum_{i=1}^{g} 
\left( \left( \int_{\alpha_{i}} \omega_{j} \right) 
\left( \int_{\beta_{i}} \frac{df}{f} \right) - 
\left( \int_{\beta_{i}} \omega_{j} \right) 
\left( \int_{\alpha_{i}} \frac{df}{f} \right) \right) 
\\
& \in & L
\end{eqnarray*}
because 
$\int_{\alpha_{i}} df/f, \int_{\beta_{i}} df/f \in 2 \pi \sqrt{-1} {\bf Z}.$ 
Hence the map $\mu$ in (3) is well-defined. 

Next, we show the injectivity of $\mu.$ 
By Riemann-Roch's theorem, 
for $P_{1}, P_{2} \in R$, 
$$
\dim_{\bf C} H^{0}(R, \Omega_{R}(P_{1} + P_{2})) = g + 1 = 
\dim_{\bf C} H^{0}(R, \Omega_{R}) + 1, 
$$ 
and hence there is a meromorphic $1$-form on $R$ which is holomorphic 
except for simple poles at $P_{1}$, $P_{2}$ with residues $1$, $-1$ respectively. 
Let $D$ be a divisor with degree $0$ on $R$ such that $\mu(D) \in L.$ 
Then $D$ is represented as a finite sum 
$\sum_{i} a_{i} \left( P_{1}^{(i)} - P_{2}^{(i)} \right)$ 
for $a_{i} \in {\bf Z}$, $P_{1}^{(i)}, P_{2}^{(i)} \in R$, 
and hence there is a meromorphic $1$-form $\psi$ on $R$ 
such that $\sum_{P \in R} {\rm Res}_{P}(\psi) \cdot P = D.$ 
Hence by the period relation, 
\begin{eqnarray*}
\mu(D) & = & 
\left( \sum_{P \in {\cal P}} {\rm Res}_{P} \left( 
\int_{P_{0}}^{P} \omega_{j} \cdot \psi \right) \right)_{1 \leq j \leq g}
\\
& = & 
\frac{1}{2 \pi \sqrt{-1}} \sum_{i=1}^{g} \left( 
\left( \int_{\alpha_{i}} \omega_{j} \right) 
\left( \int_{\beta_{i}} \psi \right) - 
\left( \int_{\beta_{i}} \omega_{j} \right) 
\left( \int_{\alpha_{i}} \psi \right) \right)_{1 \leq j \leq g}
\\ 
& \in & 
L = \left\{ \sum_{i=1}^{g} \left( \left. 
m_{i} \int_{\alpha_{i}} \omega_{j} - l_{i} \int_{\beta_{i}} \omega_{j} 
\right)_{1 \leq j \leq g} \ \right| \ m_{i}, l_{i} \in {\bf Z} \right\}. 
\end{eqnarray*}
Therefore, there are integers $m_{i}, l_{i}$ $(1 \leq i \leq g)$ such that 
$$
\sum_{i=1}^{g} \left( 
\left( \int_{\beta_{i}} \psi - (2 \pi \sqrt{-1}) m_{i} \right) 
\left( \int_{\alpha_{i}} \omega_{j} \right) - 
\left( \int_{\alpha_{i}} \psi - (2 \pi \sqrt{-1}) l_{i} \right) 
\left( \int_{\beta_{i}} \omega_{j} \right) \right) = 0 
$$ 
for any $1 \leq j \leq g.$ 
By (1), the orthogonal subspace of ${\bf C}^{2g}$ to 
$$
\left( \int_{\alpha_{1}} \omega_{j}, \cdots , \int_{\alpha_{g}} \omega_{j}, 
\int_{\beta_{1}} \omega_{j}, \cdots , \int_{\beta_{g}} \omega_{j} \right) 
\ \ (1 \leq j \leq g)
$$ 
has dimension $g,$ 
and by the period relation, this is generated by 
$$ 
\left( \int_{\beta_{1}} \omega_{j}, \cdots , \int_{\beta_{g}} \omega_{j}, 
- \int_{\alpha_{1}} \omega_{j}, \cdots , 
- \int_{\alpha_{g}} \omega_{j} \right) 
\ \ (1 \leq j \leq g).
$$ 
Hence there are $b_{1},..., b_{g} \in {\bf C}$ such that 
$$
\int_{\alpha_{i}} \psi - (2 \pi \sqrt{-1}) l_{i} = 
\sum_{j=1}^{g} b_{j} \int_{\alpha_{i}} \omega_{j}, 
\ \ 
\int_{\beta_{i}} \psi - (2 \pi \sqrt{-1}) m_{i} = 
\sum_{j=1}^{g} b_{j} \int_{\beta_{i}} \omega_{j}, 
$$
and then 
$$
f \ = \ \exp \left( \int_{P_{0}}^{P} 
\left( \psi - \sum_{j=1}^{g} b_{j} \omega_{j} \right) \right) 
$$
is a meromorphic function on $R$ such that 
$\sum_{P \in R} {\rm ord}_{P}(f) \cdot P = D.$ 
This implies the injectivity of $\mu.$ 

Finally, we show that the surjectivity of $\mu.$ 
Let $\varphi_{1}$ be a nonzero holomorphic $1$-form, 
and $Q_{1}$ be a point on $R$ at which $\varphi_{1}$ does not vanish. 
Then by Riemann-Roch's theorem, 
$$
\dim_{\bf C} H^{0}(R, \Omega_{R}(-Q_{1})) 
= \dim_{\bf C} H^{0}(R, {\cal O}_{R}(Q_{1})) + g - 2 
= g - 1, 
$$
and hence there are nonzero $\varphi_{2} \in H^{0}(R, \Omega_{R}(-Q_{1}))$ 
and $Q_{2} \in R$ at which $\varphi_{2}$ does not vanish. 
By repeating this process, one take a base $\varphi_{1},..., \varphi_{g}$ 
of $H^{0}(R, \Omega_{R})$ and $Q_{1},..., Q_{g} \in R$ such that 
$\varphi_{i+1},..., \varphi_{g},$ but not $\varphi_{i}$ vanish at $Q_{i}$ 
$(1 \leq i \leq g).$ 
Therefore, for $P_{1},..., P_{g}$ in neighborhoods of $Q_{1},..., Q_{g}$ 
respectively, the jacobian of 
$$
\left( P_{1},..., P_{g} \right) \longmapsto 
\left( \sum_{i=1}^{g} \int_{Q_{i}}^{P_{i}} \varphi_{1}, ..., 
\sum_{i=1}^{g} \int_{Q_{i}}^{P_{i}} \varphi_{g} \right) 
$$
is nonzero at $(Q_{1},..., Q_{g}),$ 
and hence by the implicit function theorem, 
the linear map $\mu$ is locally biholomorphic. 
This implies that this image ${\rm Im}(\mu)$ is an open subset of ${\bf C}^{g}/L$. 
By Riemann-Roch's theorem, 
for each divisor $D$ with degree $g$ on $R$, 
there is a nonzero element $f$ of $H^{0} \left( R, {\cal O}_{R}(D) \right)$, 
and hence 
$$
D + \sum_{P \in R} {\rm ord}_{P}(f) \cdot P 
$$
is a sum of $g$ points on $R$. 
Therefore, 
for a fixed point $P_{0} \in R$, 
the map 
$$
\left( P_{1},..., P_{g} \right) \mapsto \sum_{i=1}^{g} P_{i} - g \cdot P_{0}
$$ 
gives a holomorphic surjection from $R^{g}$ onto ${\rm Cl}^{0}(R)$. 
This implies that ${\rm Cl}^{0}(R)$ is compact, 
and hence ${\rm Im}(\mu)$ is closed in the connected set ${\bf C}^{g}/L$. 
Therefore, ${\rm Im}(\mu) = {\rm Cl}^{0}(R)$.  
QED. 
\vspace{2ex}

\noindent
\underline{\bf Exercise 2.} 
Fix $P_{0} \in R.$ Then for each $P \in R,$ 
prove that there is a unique meromorphic $1$-form $w_{P} = w_{P}(z)$ 
on $R$ such that 
\begin{itemize}

\item 
$w_{P}$ is holomorphic except $z = P, P_{0};$ 

\item 
$w_{P}$ has simple poles at $z = P, P_{0}$ with residues $1, -1$ respectively; 

\item 
${\displaystyle \int_{\alpha_{i}} w_{P} = 0 \ (1 \leq i \leq g)}$.  

\end{itemize}
Further, using the generalized Riemann's period relation, prove that 
$$
d \left( \int_{\beta_{i}} w_{z} \right) = 2 \pi \sqrt{-1} \omega_{i}(z). 
$$
\vspace{1ex}

\noindent
\underline{\bf Example.} 
If $R = {\bf C} / L$ : genus 1; $L = {\bf Z} + {\bf Z} \tau$ 
$({\rm Im}(\tau) > 0),$ then 
$$
H^{0}(R, \Omega_{R}) = {\bf C} dz, \ \ 
Z = \int_{0}^{\tau} dz = \tau. 
$$

\newpage
\begin{center}
\underline{\large {\bf \S 3. Schottky uniformization}} 
\end{center}

\noindent
\underline{\bf 3.1. Degeneration of Riemann surfaces} 
\vspace{2ex}

\noindent
\underline{\bf Genus 1 case.} 
If $f(x)$(: degree 3, without multiple root) tends to 
$a (x - \alpha)^{2}(x - \beta)$ $(a \neq 0, \alpha \neq \beta),$ 
then the complex torus $C_{f}({\bf C})$ degenerates to 
a singular space obtained by identifying 2-points on the Riemann sphere 
({\bf Figure}). 

For example, for $f(x) = (x^{2} - \varepsilon^{2})(x + 1),$ 
\begin{eqnarray*}
& & y^{2} = f(x) \Leftrightarrow 
(\sqrt{x^{2} + x^{3}} + y)(\sqrt{x^{2} + x^{3}} - y) = \varepsilon^{2}(1+x) \\ 
& \stackrel{\varepsilon \rightarrow 0}{\longrightarrow} & 
(\sqrt{x^{2} + x^{3}} + y)(\sqrt{x^{2} + x^{3}} - y) = 0 \ \ 
\mbox{around $(x, y) = (0, 0)$}, 
\end{eqnarray*}
where $\sqrt{x^{2} + x^{3}} = \sum_{k=0}^{\infty} 
\left( \begin{array}{c} 1/2 \\ k \end{array} \right) x^{k+1}.$  
\vspace{2ex}

\noindent
\underline{\bf Local degeneration.} 
For a complex number $\varepsilon$ such that 
$0 < |\varepsilon| < 1,$ 
let $D$ be the union of the two annular domains: 
$$
U = \{ x \in {\bf C} \ | \ |\varepsilon| < |x| < 1 \}, \ \ 
V = \{ y \in {\bf C} \ | \ |\varepsilon| < |y| < 1 \} 
$$
by the relation $x y = \varepsilon.$ 
Then under $\varepsilon \rightarrow 0,$ 
$D$ becomes the union of the 2 disks 
$$
\{ x \in {\bf C} \ | \ |x| < 1 \}, \ \ 
\{ y \in {\bf C} \ | \ |y| < 1 \} 
$$
identifying $x = 0$ and $y = 0.$ 
\vspace{2ex}

\noindent
\underline{\bf Ordinary double points.} 
For a point $P$ on a curve $C,$ 
\begin{eqnarray*}
& & \mbox{$P$ is an {\bf ordinary double point} (or {\bf node})} 
\\
& \stackrel{\rm def}{\Longleftrightarrow} & 
\left\{ \begin{array}{l} 
\mbox{the local equation around $P \in C$ is given by $xy = 0$} \\ 
\mbox{for some formal coordinates $x,y$} 
\end{array} \right. 
\\ 
& \Longleftrightarrow & 
\mbox{$P$ is a point of multiplicity $2$ 
with distinct tangent directions} 
\end{eqnarray*}
\vspace{-2ex}

\noindent
\underline{\bf 3.2. Schottky uniformization of Riemann surfaces}
\vspace{2ex}

{\bf Schottky uniformization} is to construct Riemann surfaces of genus $g$ 
from a $2g$ holed Riemann sphere by identifying these holes in pairs 
({\bf Figure}). 
More precisely, let 
$$
PGL_{2}({\bf C}) \ \stackrel{\rm def}{=} \ 
\left. GL_{2}({\bf C}) \right/ {\bf C}^{\times} (\cdot E_{2}) 
$$
which acts on ${\mathbb P}^{1}({\bf C})$ by the M\"{o}bius transformation, 
and let 
\begin{eqnarray*}
& & D_{\pm 1},... , D_{\pm g} \subset {\mathbb P}^{1}({\bf C}) : \ 
\mbox{disjoint closed domains bounded by Jordan curves $\partial D_{i},$} 
\\ 
& & \gamma_{1},..., \gamma_{g} \in PGL_{2}({\bf C}) \ 
\mbox{such that 
$\gamma_{i}({\mathbb P}^{1}({\bf C}) - D_{-i}) =$ 
the interior $D_{i}^{\circ}$ of $D_{i},$}
\\
& & 
\Gamma \ \stackrel{\rm def}{=} \ \langle \gamma_{1},..., \gamma_{g} \rangle 
: \ \mbox{the subgroup of $PGL_{2}({\bf C})$ generated by 
$\gamma_{1},..., \gamma_{g},$} 
\\
& & 
\Omega_{\Gamma} \ \stackrel{\rm def}{=} \ 
\bigcup_{\gamma \in \Gamma} \gamma 
\left( {\mathbb P}^{1}({\bf C}) - \bigcup_{i=1}^{g} 
(D_{i}^{\circ} \cup D_{-i}^{\circ}) \right). 
\end{eqnarray*}
Then the Riemann surface 
\begin{eqnarray*}
R_{\Gamma} & \stackrel{\rm def}{=} & 
\left. \left( {\mathbb P}^{1}({\bf C}) - \bigcup_{i=1}^{g} 
(D_{i}^{\circ} \cup D_{-i}^{\circ}) \right) \right/ 
\partial D_{i} \stackrel{\gamma_{i}}{\sim} \partial D_{-i} \ 
\mbox{(: gluing by $\gamma_{i}$)}
\\ 
& = & 
\Omega_{\Gamma} / \Gamma 
\end{eqnarray*}
is called (Schottky) uniformized by the {\bf Schottky group} $\Gamma.$ 
It is known as Koebe's theorem that any Riemann surface can be Schottky uniformized. 
Counterclockwise oriented loops $\partial D_{i}$ and 
oriented paths from $w_{i} \in \partial D_{-i}$ to 
$\gamma_{i}(w_{i}) \in \partial D_{i}$ $(1 \leq i \leq g)$ 
become canonical generators, 
and we denote them by $\alpha_{i},$ $\beta_{i}$ respectively ({\bf Figure}). 
\vspace{2ex}

\noindent
\underline{\bf Remark.} 
Denote by ${\mathbb H}^{3}$ the $3$-dimensional hyperbolic space. 
Then the quotient hyperbolic $3$-manifold ${\mathbb H}^{3}/\Gamma$ becomes 
a handlebody whose boundary is $R_{\Gamma}$. 
\vspace{2ex}

\noindent
\underline{\bf Exercise 3.} 
Prove that $\Gamma$ is a free group with generators 
$\gamma_{1},..., \gamma_{g},$ 
and that the action of $\Gamma$ on $\Omega_{\Gamma}$ 
is free and properly discontinuous. 
Further, prove that each $\gamma_{i}$ $(1 \leq i \leq g)$ 
is uniquely represented by 
$$
\gamma_{i} = 
\left( \begin{array}{cc} t_{i} & t_{-i} \\ 1 & 1 \end{array} \right) 
\left( \begin{array}{cc} 1 & 0 \\ 0 & s_{i} \end{array} \right) 
\left( \begin{array}{cc} t_{i} & t_{-i} \\ 1 & 1 \end{array} \right)^{-1}
\ {\rm mod}({\bf C}^{\times}), 
$$
where $t_{i} \in D_{i}^{\circ},$ $t_{-i} \in D_{-i}^{\circ}$ 
and $|s_{i} | < 1$ (hence $\gamma_{i}$ is hyperbolic (or loxodromic)), 
and that 
$$
t_{\pm i} \ = \ \lim_{n \rightarrow \infty} \gamma_{i}^{\pm n}(z) \ 
(z \in \Omega_{\Gamma}). 
$$
$t_{i},$ $t_{-i}$ are called the {\bf attractive}, {\bf repulsive} 
fixed point of $\gamma_{i}$ respectively, and 
$s_{i}$ is called the {\bf multiplier} of $\gamma_{i}.$ 
\vspace{2ex}

\noindent
\underline{\bf 3.3. Explicit formula of periods}
\vspace{2ex}

\noindent
\underline{\bf Theorem 3.1.} (Schottky [S]) \begin{it} 
Assume that $\infty \in \Omega_{\Gamma}$ and that 
$\sum_{\gamma \in \Gamma} |\gamma'(z)|$ converges uniformly 
on any compact subset of 
$$
\Omega_{\Gamma} - \bigcup_{\gamma \in \Gamma} \gamma(\infty). 
$$
Then we have 

{\rm (1)} For $n \geq 1$ and a point 
$p \in \Omega_{\Gamma} - \bigcup_{\gamma \in \Gamma} \gamma(\infty),$ 
$$ 
w_{n,p}(z) \ \stackrel{\rm def}{=} \ 
\sum_{\gamma \in \Gamma} \frac{d \gamma(z)}{(\gamma(z) - p)^{n}} = 
\sum_{\gamma \in \Gamma} \frac{\gamma'(z)}{(\gamma(z) - p)^{n}} dz 
$$ 
becomes a meromorphic 1-form on $R_{\Gamma}.$ 
If $n > 1,$ then $w_{n,p}$ is of the 2-nd kind, 
and it has only poles (of order $n$) at the point $\overline{p}$ on $R_{\Gamma}$ 
induced from $p$. 
If $n = 1,$ then $w_{n,p}$ is of the 3-rd kind, 
and it has only simple poles at $\overline{p}, \overline{\infty}$ 
with residues $1, -1$ respectively. 
Furthermore, for $n \geq 0,$ 
$$
\sum_{\gamma \in \Gamma} \gamma(z)^{n} d \gamma(z) = 
\sum_{\gamma \in \Gamma} \gamma(z)^{n} \cdot \gamma'(z) dz 
$$
becomes a meromorphic 1-form on $R_{\Gamma}$ 
which has only pole (of order $n+2$) at $\overline{\infty}.$ 

{\rm (2)} For $i = 1,...,g,$ 
$$
\omega_{i}(z) \ = \ \frac{1}{2 \pi \sqrt{-1}} 
\sum_{\gamma \in \Gamma/\langle \gamma_{i} \rangle} 
\left( \frac{1}{z - \gamma(t_{i})} - \frac{1}{z - \gamma(t_{-i})} \right) 
dz
$$
give a basis of $H^{0}(R_{\Gamma}, \Omega_{R_{\Gamma}})$ satisfying that 
${\displaystyle \int_{\alpha_{i}} \omega_{j} = \delta_{ij}.}$ 

{\rm (3)} For $1 \leq i,j \leq g$ and $\gamma \in \Gamma,$ put
$$
\psi_{ij}(\gamma) = \left\{ \begin{array}{ll} 
s_{i} & \mbox{(if $i=j$ and $\gamma \in \langle \gamma_{i} \rangle$),} \\ 
{\displaystyle \frac{(t_{i}-\gamma(t_{j}))(t_{-i}-\gamma(t_{-j}))}
{(t_{i}-\gamma(t_{-j}))(t_{-i}-\gamma(t_{j}))}} & \mbox{(otherwise),} 
\end{array} \right. 
$$
where $t_{i}$, $t_{-i}$ are the attractive, repulsive fixed points of $\gamma_{i}$ 
respectively, 
and $s_{i}$ is the multiplier of $\gamma_{i}$. 
Then we have 
$$
\exp \left( 2 \pi \sqrt{-1} z_{ij} \right) \ = \ 
\prod_{\gamma \in \langle \gamma_{i} \rangle \backslash \Gamma / 
\langle \gamma_{j} \rangle} \psi_{ij}(\gamma), 
$$ 
where $Z = (z_{ij})_{i,j}$ is the period matrix of 
$(R_{\Gamma}; (\alpha_{i}, \beta_{i})_{1 \leq i \leq g}).$ 
\end{it}
\vspace{2ex}

\underline{\it Proof.} 
The assertion (1) is evident except the convergence of $w_{n,p}(z)$ 
which follows from the assumption and that 
the action of $\Gamma$ on $\Omega_{\Gamma}$ is properly discontinuous. 
Further, $w_{1,p}(z)$ has simple poles at $\overline{p}, \overline{\infty}$ 
with residues $1, -1$ respectively, 
and satisfies that $\int_{\alpha_{i}} w_{1,p} = 0$ $(1 \leq i \leq g).$ 
Then by Exercise 2, 
\begin{eqnarray*}
2 \pi \sqrt{-1} \omega_{i}(z) & = & 
d \left( 
\int_{\zeta_{i}}^{\gamma_{i}(\zeta_{i})} 
\sum_{\gamma \in \Gamma} 
\frac{d \gamma(\zeta)}{\gamma(\zeta) - z} \right) ; \ 
\mbox{$\zeta_{i}$ is a point on $\partial D_{-i}$} 
\\
& = & 
d \left( \sum_{\gamma \in \Gamma} \log 
\left( 
\frac{(\gamma \gamma_{i})(\zeta_{i}) - z}{\gamma(\zeta_{i}) - z} 
\right) \right)  
\\ 
& = & 
\sum_{\gamma \in \Gamma} \left( 
\frac{1}{z - (\gamma \gamma_{i})(\zeta_{i})} - \frac{1}{z - \gamma(\zeta_{i})} 
\right) dz 
\\ 
& = & 
\sum_{\gamma \in \Gamma/\langle \gamma_{i} \rangle} \sum_{n \in {\bf Z}} 
\left( 
\frac{1}{z - (\gamma \gamma_{i}^{n+1})(\zeta_{i})} - 
\frac{1}{z - (\gamma \gamma_{i}^{n})(\zeta_{i})} 
\right) dz, 
\end{eqnarray*}
and since $t_{\pm i} = \lim_{n \rightarrow \infty} \gamma_{i}^{\pm n} (w_{i}) 
\in D_{\pm i}^{\circ}$ (Exercise 3), we have 
$$
\omega_{i}(z) = \frac{1}{2 \pi \sqrt{-1}} 
\sum_{\gamma \in \Gamma/\langle \gamma_{i} \rangle} 
\left( 
\frac{1}{z - \gamma(t_{i})} - \frac{1}{z - \gamma(t_{-i})} 
\right) dz, 
$$
which proves (2). QED. 
\vspace{2ex}

\noindent
\underline{\bf Exercise 4.} 
Prove that $\int_{\alpha_{i}} w_{1,p} = 0$ $(1 \leq i \leq g),$ 
and check that $\omega_{i}$ is $\Gamma$-invariant and 
$\int_{\alpha_{i}} \omega_{j} = \delta_{ij}.$ 
\vspace{2ex} 

\noindent
\underline{\bf Exercise 5.} 
Prove (3) of Theorem 3.1.
\vspace{2ex}

\noindent
\underline{\bf Proposition 3.2.} \begin{it} 
Assume that $\Omega_{\Gamma} \ni \infty,$ and that $t_{\pm i}$ are fixed 
and $s_{i}$ are sufficiently small, 
then the assumption in Theorem 3.1 is satisfied. 
\end{it}
\vspace{2ex}

\underline{\it Proof.} 
For 2 disks $D_{i}, D_{j} \subset {\bf C}$ with radius $r_{i}, r_{j}$ 
respectively, put 
\begin{eqnarray*} 
\rho_{i,j} & : & 
\mbox{the distance between the centers of $D_{i}$ and $D_{j},$} 
\\ 
K_{i,j} & = & 
{\displaystyle \frac{(r_{i}^{2} + r_{j}^{2} - \rho_{i,j}^{2})^{2}}
{4 r_{i}^{2} r_{j}^{2}} - 1 \geq 0,}  
\\ 
L_{i,j} & = & 
{\displaystyle \frac{1}{\sqrt{1 + K_{i,j}} + \sqrt{K_{i,j}}} \leq 1.}  
\end{eqnarray*} 
Then $K_{i,j}$ and $L_{i,j}$ are invariant 
under any M\"{o}bius transformation, 
and $r_{i} \leq L_{i,j} \cdot r_{j}$ if $D_{i} \subset D_{j}.$ 
Under the assumption, one can take disks $D_{\pm 1},..., D_{\pm g}$ such that 
the sum of $L_{i,j}$ $(i,j \in \{ \pm 1,..., \pm g \}, i \neq j)$ 
is smaller than $1.$ 
Hence by the above, there is a positive constant $C$ such that if 
$\gamma = \prod_{s=1}^{l} \gamma_{k(s)} \in \Gamma$ is expressed as 
$$
\left( \begin{array}{cc} a_{\gamma} & b_{\gamma} \\ c_{\gamma} & d_{\gamma} 
\end{array} \right) \ {\rm mod}({\bf C}^{\times}) ; 
\left( \begin{array}{cc} a_{\gamma} & b_{\gamma} \\ c_{\gamma} & d_{\gamma} 
\end{array} \right) \in SL_{2}({\bf C}), 
$$ 
then 
$$
\frac{1}{|c_{\gamma}|^{2}} \leq C \cdot \prod_{s=1}^{l-1} L_{-k(s), k(s+1)}. 
$$
Therefore, 
$$
\sum_{\gamma \in \Gamma - \{ 1 \}} \frac{1}{|c_{\gamma}|^{2}} \leq 
C \cdot \sum_{m=0}^{\infty} \left( \sum_{i \neq j} L_{i,j} \right)^{m} 
< \infty,
$$
and hence 
$$
\sum_{\gamma \in \Gamma} |\gamma'(z)| \leq 1 + \frac{1}{d(z)^{2}} 
\sum_{\gamma \in \Gamma - \{ 1 \}} \frac{1}{|c_{\gamma}|^{2}} 
$$
satisfies the condition 
since $d(z) \stackrel{\rm def}{=} 
\min \{ |z - \gamma^{-1}(\infty)| ; \gamma \in \Gamma \} > 0$ is bounded 
on any compact subset outside $\bigcup_{\gamma \in \Gamma} \gamma(\infty).$ 
QED 
\vspace{2ex}

\noindent
\underline{\bf Remark.} 
Schottky [S] gives a (more geometric) convergence condition 
on $\sum_{\gamma \in \Gamma} |\gamma'(z)|$ as follows: 
all $\partial D_{\pm i}$ can be taken as circles 
(in this case, $\Gamma$ is called classical) and 
there are $2g - 3$ circles $C_{1},..., C_{2g-3}$ in 
$F = {\mathbb P}^{1}({\bf C}) - 
\bigcup_{i=1}^{g} \left( D_{i}^{\circ} \cup D_{-i}^{\circ} \right)$ 
satisfying that 
\begin{itemize}

\item 
$C_{1},..., C_{2g-3},$ $\partial D_{\pm 1},..., \partial D_{\pm g}$ 
are mutually disjoint; 

\item 
$C_{1},..., C_{2g-3}$ divide $F$ into $2g - 2$ domains $R_{1},..., R_{2g-2};$ 

\item 
each $R_{i}$ has exactly three boundary circles. 

\end{itemize}

\noindent
\underline{\bf Variation of forms and periods.} 
Let $\Gamma = \langle \gamma_{1},..., \gamma_{g} \rangle$ be 
a Schottky group of rank $g$ as above, 
and put $\Gamma' = \langle \gamma_{1},...,\gamma_{g-1} \rangle$ 
which is a Schottky group of rank $g-1.$ 
If the multiplier 
\begin{eqnarray*}
s_{g} 
& = & 
\frac{\gamma_{g}(z) - t_{g}}{z - t_{g}} \cdot 
\frac{z - t_{-g}}{\gamma_{g}(z) - t_{-g}} 
\\
& : & 
\mbox{the product of local coordinates around $t_{g}, t_{-g}$ respectively} 
\end{eqnarray*}
of $\gamma_{g}$ tends to $0,$ then 
\begin{itemize}
\item 
$R_{\Gamma} \ \longrightarrow \ 
\left\{ \begin{array}{l} 
\mbox{the singular curve $\widetilde{R}_{\Gamma'}$ with unique singular 
(ordinary double) point} 
\\
\mbox{obtained from $R_{\Gamma'}$ by identifying $t_{g}$ and $t_{-g};$} 
\end{array} \right.$

\item 
${\displaystyle 2 \pi \sqrt{-1} \ \omega_{i}(z) = 
\sum_{\gamma \in \Gamma/\langle \gamma_{i} \rangle} 
\left( \frac{1}{z - \gamma(t_{i})} - \frac{1}{z - \gamma(t_{-i})} \right) dz 
\in H^{0} \left( R_{\Gamma}, \Omega_{R_{\Gamma}} \right)}$ 
\begin{eqnarray*}
& \longrightarrow & 
\left\{ \begin{array}{ll} 
{\displaystyle \sum_{\gamma \in \Gamma'/\langle \gamma_{i} \rangle} 
\left( \frac{1}{z - \gamma(t_{i})} - \frac{1}{z - \gamma(t_{-i})} \right) dz} 
& (i < g), \\ 
{\displaystyle 
\left( \frac{1}{z - t_{g}} - \frac{1}{z - t_{-g}} \right) dz + \cdots} 
& (i = g) 
\end{array} \right. 
\end{eqnarray*} 
which has a pole at the ordinary double point $t_{g} = t_{-g}$ on 
$\widetilde{R}_{\Gamma'}$ if $i = g;$ 

\item (Fay's formula [Fay]) \ 
${\displaystyle p_{ij} \longrightarrow 
\left\{ \begin{array}{ll} 
\mbox{the multiplicative periods of $R_{\Gamma'}$} & (i,j < g), \\ 
0 & (i = j = g). 
\end{array} \right.}$ 

\end{itemize}
Therefore, on the complex geometry of $\widetilde{R}_{\Gamma'},$ 
it is natural to replace the sheaf of holomorphic $1$-forms on $\widetilde{R}_{\Gamma'}$ 
by that of $1$-forms $\eta$ on $R_{\Gamma'}$ 
holomorphic except for simple poles at $t_{g}, t_{-g}$ 
satisfying that ${\rm Res}_{t_{g}}(\eta) + {\rm Res}_{t_{-g}}(\eta) = 0$ 
(see 5.2 below). 
\vspace{2ex}

\noindent
\underline{\bf Remark.} 
We can obtain variational formula under other degenerations (see [I3]). 
\vspace{2ex}

\noindent
\underline{\bf 3.4. Fractal nature of Schottky groups} (cf. [MumSW]) 
\vspace{2ex}

\noindent
\underline{\bf Limit set.} 
The {\bf limit set} $L_{\Gamma}$ of $\Gamma$ is defined to be 
the complement of $\Omega_{\Gamma}$ in ${\mathbb P}^{1}({\bf C})$. 
When we take limits of Schottky groups such as the above domains 
$D_{\pm i} \subset {\mathbb P}^{1}({\bf C})$ are tangent, 
their limit sets become {\it fractal} pictures. 
These limiting process is also important in the study of 
degeneration of $R_{\Gamma}$ and deformation of hyperbolic $3$-manifolds. 

Let $\Gamma$ be a Kleinian group generated by 
$\gamma_{1}, \gamma_{2} \in PGL_{2}({\bf C})$ 
such that there are tangential disks 
$D_{\pm 1}, D_{\pm 2} \subset {\mathbb P}^{1}({\bf C})$ satisfying 
$$
\gamma_{i} \left( {\mathbb P}^{1}({\bf C}) - D_{-i} \right) = D_{i}^{\circ}, \ 
\gamma_{i} (\mbox{\{tangential points\}}) = \mbox{\{tangential points\}} \ 
(i = 1, 2), 
$$ 
and consider the following cases. 
\vspace{2ex}

\noindent
\underline{\bf Case 1.} 
$D_{\pm 1} \cup D_{\pm 2} \subset {\mathbb P}^{1}({\bf C})$ is {\it homeomorphic to} 
$$
\left\{ z \in {\bf C}; |z - \pm (1 + \sqrt{-1})| \leq 1 \right\} \cup 
\left\{ z \in {\bf C}; |z - \pm (1 - \sqrt{-1})| \leq 1 \right\}.  
$$ 
Then $\Gamma$ is called a once-punctured torus group, 
and $R_{\Gamma}$ becomes a union of two tori at one point (cf. [MumSW, p.189--190]). 
\vspace{2ex}

\noindent
\underline{\bf Case 2.} 
Each of $D_{\pm 1}, D_{\pm 2} \subset {\mathbb P}^{1}({\bf C})$ 
is tangent to other three disks. 
Then $L_{\Gamma}$ becomes an Apollonian gasket, 
and $R_{\Gamma}$ becomes a union of two degenerate tori at one point 
(cf. [MumSW, p.205]). 
\vspace{2ex}

\noindent
\underline{\bf Case 3.} 
$D_{\pm 1} \cup D_{\pm 2} \subset {\mathbb P}^{1}({\bf C})$ is {\it homeomorphic to} 
$$
\left\{ z \in {\bf C}; |z - (\pm 1 + \sqrt{-1})| \leq 1 \right\} \cup 
\left\{ z \in {\bf C}; |z - (\pm 1 - \sqrt{-1})| \leq 1 \right\}.  
$$
Then $R_{\Gamma}$ becomes a union of two spheres at three points 
(cf. [MumSW, p.216]).

\newpage
\begin{center}
\underline{\large {\bf \S 4. Arithmetic uniformization}}
\end{center}

\noindent
\underline{\bf 4.1. Periods as power series} 
\vspace{2ex}

First, we calculate the periods $p_{ij}$ given in Theorem 3.1 (3) (cf. Exercise 5) 
as power series over ${\bf Z}$ by regarding 
the fixed points and multiplier $t_{\pm i}, s_{i}$ of $\gamma_{i}$ 
as variables $x_{\pm i}, y_{i}$ respectively. 
Let $\Gamma_{\Delta}$ be a subgroup of $PGL_{2}$ generated by 
$$
\phi_{i} \stackrel{\rm def}{=}  
\left( \begin{array}{cc} x_{i} & x_{-i} \\ 1 & 1 \end{array} \right) 
\left( \begin{array}{cc} 1     & 0  \\ 0 & y_{i} \end{array} \right) 
\left( \begin{array}{cc} x_{i} & x_{-i} \\ 1 & 1 \end{array} \right)^{-1} 
{\rm mod}({\bf G}_{m}) \ \ (1 \leq i \leq g). 
$$
Put 
\begin{eqnarray*}
A_{0} & = & {\bf Z} 
\left[ 
\frac{(x_{i} - x_{j})(x_{k} - x_{l})}{(x_{i} - x_{l})(x_{k} - x_{j})} \ \ 
\left( \begin{array}{lcl} 
i,j,k,l & \in & \{ \pm 1,..., \pm g \} \\ 
        &  :  & \mbox{mutually different} 
\end{array} \right) \right], 
\\ 
A_{\Delta} & = & A_{0} [[ y_{1},..., y_{g} ]], 
\end{eqnarray*}
and let $I_{\Delta}$ be the ideal of $A_{\Delta}$ generated by $y_{1},...,y_{g}$. 
By definition,  
$$
p_{ij} \ = \ 
\prod_{\phi \in \langle \phi_{i} \rangle \backslash \Gamma_{\Delta} / 
\langle \phi_{j} \rangle} \psi_{ij}(\phi), 
$$ 
where 
$$
\psi_{ij}(\phi) = \left\{ \begin{array}{ll} 
y_{i} & \mbox{(if $i=j$ and $\phi \in \langle \phi_{i} \rangle$),} \\ 
{\displaystyle \frac{(x_{i} - \phi(x_{j}))(x_{-i} - \phi(x_{-j}))}
{(x_{i} - \phi(x_{-j}))(x_{-i} - \phi(x_{j}))}} & \mbox{(otherwise).} 
\end{array} \right. 
$$
Put $\phi_{-i} \stackrel{\rm def}{=} \phi_{i}^{-1}$ $(1 \leq i \leq g).$ 
Then 
$$
\Phi_{ij} = \left\{ \phi = \phi_{\sigma(1)} \cdots \phi_{\sigma(n)} \ 
\left| \begin{array}{l} 
\sigma(1) \neq \pm i, \ \sigma(n) \neq \pm j, 
\\
\sigma(k) \neq - \sigma(k+1) \ (1 \leq k \leq n-1) 
\end{array} \right. \right\} 
$$
gives a set of complete representatives of 
$\langle \phi_{i} \rangle \backslash \Gamma_{\Delta} / 
\langle \phi_{j} \rangle.$ 
If $\alpha \in x_{j} + I_{\Delta}$ with $j \neq -i$, 
then 
$$
\phi_{i}(\alpha) = 
\left( x_{i} - \frac{(\alpha - x_{i}) x_{-i} y_{i}}{\alpha - x_{-i}} \right) 
\left( 1 - \frac{(\alpha - x_{i}) y_{i}}{\alpha - x_{-i}} \right)^{-1} 
\in x_{i} + I_{\Delta}. 
$$ 
Hence if $\phi = \phi_{\sigma(1)} \cdots \phi_{\sigma(n)} \in \Phi_{ij}$, 
then $\phi(x_{\pm j}) \in x_{\sigma(1)} + I_{\Delta}$ and 
\begin{eqnarray*}
& & \phi(x_{j}) - \phi(x_{-j}) 
\\ 
& = & 
\frac{(x_{\sigma(1)} - x_{-\sigma(1)})^{2} 
(\phi'(x_{j}) - \phi'(x_{-j})) y_{\sigma(1)}}
{(\phi'(x_{j}) - x_{-\sigma(1)} - 
y_{\sigma(1)}(\phi'(x_{j}) - x_{\sigma(1)}))
(\phi'(x_{-j}) - x_{-\sigma(1)} - 
y_{\sigma(1)}(\phi'(x_{-j}) - x_{\sigma(1)}))} 
\\ 
& & (\phi' \stackrel{\rm def}{=} \phi_{\sigma(2)} \cdots \phi_{\sigma(n)}) 
\\ 
& = & \cdots \ \in I_{\Delta}^{n}. 
\end{eqnarray*}
by inductive calculus. 
Therefore,  
$$
\frac{(x_{i} - \phi(x_{j}))(x_{-i} - \phi(x_{-j}))}
{(x_{i} - \phi(x_{-j}))(x_{-i} - \phi(x_{j}))} 
= 1 + 
\frac{(x_{i} - x_{-i})(\phi(x_{j}) - \phi(x_{-j}))}
{(x_{i} - \phi(x_{-j}))(x_{-i} - \phi(x_{j}))} 
\ \in \ 1 + I_{\Delta}^{n}, 
$$
and hence $p_{ij}$ are elements of $A_{\Delta}$ calculated as 
$$
p_{ij} \ = \ 
c_{ij} \left( 1 + \sum_{|k| \neq i,j} 
\frac{(x_{i} - x_{-i})(x_{j} - x_{-j})(x_{k} - x_{-k})^{2}}
{(x_{i} - x_{k})(x_{-i} - x_{k})(x_{j} - x_{-k})(x_{-j} - x_{-k})} 
y_{|k|} + \cdots \right), 
$$
where 
$$
c_{ij} \ \stackrel{\rm def}{=} \ 
\left\{ \begin{array}{ll}
y_{i} & \mbox{(if $i = j),$} 
\\
{\displaystyle 
\frac{(x_{i} - x_{j})(x_{-i} - x_{-j})}{(x_{i} - x_{-j})(x_{-i} - x_{j})} } 
& \mbox{(if $i \neq j).$} 
\end{array} \right. 
$$
\vspace{1ex}

\noindent
\underline{\bf 4.2. Tate curve and Mumford curves} 
\vspace{2ex}

In order to study geometric meaning of the above calculation, 
we review the theory of the Tate curve and its higher genus version 
called arithmetic uniformization theory. 
\vspace{2ex}

\noindent
\underline{\bf Tate curve.}
Recall that an elliptic curve ${\bf C}/L$ is defined by the equation 
(see 2.1): 
$$
y^{2} = 4 x^{3} - 60 E_{4}(L) x - 140 E_{6}(L). 
$$
Therefore, if 
$$
\begin{array}{lllllll}
x & = & 
{\displaystyle (2 \pi \sqrt{-1})^{2} \left( X + \frac{1}{12} \right),} 
& & 
y & = & {\displaystyle (2 \pi \sqrt{-1})^{3} \left( 2 Y + X \right),} 
\\ 
a_{4} & = & 
{\displaystyle - \frac{15 E_{4}(L)}{(2 \pi \sqrt{-1})^{4}} + \frac{1}{48},} 
& & 
a_{6} & = & {\displaystyle - \frac{35 E_{6}(L)}{(2 \pi \sqrt{-1})^{6}} 
- \frac{5 E_{4}(L)}{4 (2 \pi \sqrt{-1})^{4}} + \frac{1}{1728},} 
\end{array}
$$
then the above equation is equivalent to 
$$
Y^{2} + X Y = X^{3} + a_{4} X + a_{6}. 
$$
Furthermore, if $L = {\bf Z} + {\bf Z} \tau$ and 
$q = e^{2 \pi \sqrt{-1} \tau},$  
then by the calculation of the {\bf Eisenstein series} 
(see Exercise 6 below): 
$$
\sum_{u \in L-\{0\}} \frac{1}{u^{2k}} = 
2 \zeta(2k) + \frac{2 (2 \pi \sqrt{-1})^{2k}}{(2k-1)!} 
\sum_{n=1}^{\infty} \sigma_{2k-1}(n) \ q^{n} \ \ (k > 1), 
$$
where
$$
\zeta(2k) \stackrel{\rm def}{=} 
\sum_{n=1}^{\infty} \frac{1}{n^{2k}}: \ \mbox{the {\bf zeta values},} 
\ \ \mbox{and} \ \ 
\sigma_{2k-1}(n) \stackrel{\rm def}{=} \sum_{d|n} d^{2k-1}, 
$$
we have 
\begin{eqnarray*}
a_{4}(q) & = & - 5 \sum_{n=1}^{\infty} \sigma_{3}(n) \ q^{n} 
= -5 q - 45 q^{2} + \cdots, 
\\ 
a_{6}(q) & = & - \frac{1}{12} 
\sum_{n=1}^{\infty} (5 \sigma_{3}(n) + 7 \sigma_{5}(n)) \ q^{n} 
= - q - 23 q^{2} + \cdots. 
\end{eqnarray*}
\vspace{1ex}

\noindent
\underline{\bf Exercise 6.} 
\noindent
Prove that 
\begin{eqnarray*}
\zeta(2k) & = & 
- \frac{(2 \pi \sqrt{-1})^{2k}}{2(2k)!} B_{2k} 
\\
& & \left( \mbox{$B_{n}$ is the $n$-th {\bf Bernoulli numbers} given by 
${\displaystyle 
\frac{x}{e^{x}-1} = \sum_{n=0}^{\infty} B_{n} \frac{x^{n}}{n!}}$} \right) 
\\
& \Rightarrow & \zeta(2) = \frac{\pi^{2}}{6}, \ 
\zeta(4) = \frac{\pi^{4}}{90}, \ \zeta(6) = \frac{\pi^{6}}{945}, 
\end{eqnarray*}
and 
$$
\sum_{(m,n) \in {\bf Z}^{2}-\{(0,0)\}} \frac{1}{(m + n \tau)^{2k}} = 
2 \zeta(2k) + \frac{2 (2 \pi \sqrt{-1})^{2k}}{(2k-1)!} 
\sum_{n=1}^{\infty} \sigma_{2k-1}(n) \ q^{n} \ \ (k > 1), 
$$
from the well-known formula: 
$$
\pi \cot(\pi a) = 
\frac{1}{a} + \sum_{m=1}^{\infty} 
\left( \frac{1}{a + m} + \frac{1}{a - m} \right) \ 
\left( \Leftrightarrow \ 
\sin z = z \prod_{n=1}^{\infty} \left( 1 - \frac{z^{2}}{n^{2} \pi^{2}} \right) 
\right)
$$
by substituting $x$ to $2 \pi \sqrt{-1} a,$ 
and differentiating the formula successively and substituting $n \tau$ to $a$ 
respectively. 
\vspace{2ex}

\noindent
\underline{\bf Exercise 7.} 
Show that $a_{4}(q)$ and $a_{6}(q)$ belong to the ring 
$$
{\bf Z}[[q]] \stackrel{\rm def}{=} 
\left\{ \left. \sum_{n=0}^{\infty} c_{n} q^{n} \ \right| \ 
c_{n} \in {\bf Z} \right\} 
$$
of formal power series of $q$ with coefficients in ${\bf Z}.$ 
\vspace{2ex}

The {\bf Tate curve} is the curve over ${\bf Z}[[q]]$ 
defined by 
$$
y^{2} + x y = x^{3} + a_{4}(q) x + a_{6}(q). 
$$
Then Tate proved the following: 
\vspace{2ex}

\noindent
\underline{\bf Theorem 4.1.} ([Si, T]) \begin{it} 

{\rm (1)} The Tate curve becomes an elliptic curve over the ring 
$$
{\bf Z}((q)) \stackrel{\rm def}{=} {\bf Z}[[q]] \left[ 1/q \right] 
= \left\{ \left. \sum_{n > m }^{\infty} c_{n} q^{n} \ \right| \ 
m \in {\bf Z}, \ c_{n} \in {\bf Z} \right\} 
$$
of Laurent power series of $q$ with coefficients in ${\bf Z}.$ 

{\rm (2)} Recall $\sigma_{1}(n) = \sum_{d|n} d$, 
and put 
\begin{eqnarray*}
X(u, q) & = & \sum_{n \in {\bf Z}} \frac{q^{n} u}{(1-q^{n}u)^{2}} 
- 2 \sum_{n=1}^{\infty} \sigma_{1}(n) q^{n}, 
\\ 
Y(u, q) & = & \sum_{n \in {\bf Z}} \frac{(q^{n} u)^{2}}{(1-q^{n}u)^{3}} 
+ \sum_{n=1}^{\infty} \sigma_{1}(n) q^{n}. 
\end{eqnarray*}
Then $z \mapsto 
\left( X(e^{2 \pi \sqrt{-1} z}, e^{2 \pi \sqrt{-1} \tau}),  
Y(e^{2 \pi \sqrt{-1} z}, e^{2 \pi \sqrt{-1} \tau}) \right)$ 
gives rise to an isomorphism between ${\bf C}/L$ 
and the elliptic curve $E_{\tau}$ over ${\bf C}$ obtained from the Tate curve 
by substituting $q = e^{2 \pi \sqrt{-1} \tau}.$ 

{\rm (3)} Let $K$ be a complete valuation field 
with multiplicative valuation $| \cdot |,$ 
and let $q \in K^{\times}$ such that $|q| < 1.$ 
Then by the substitution the variable $q \mapsto q \in K^{\times},$ 
the series $a_{4}(q)$ and $a_{6}(q)$ converge in $K,$ 
and the Tate curve gives an elliptic curve $E_{q}$ over $K.$ 
Further, we have an isomorphism: 
\begin{eqnarray*}
K^{\times} / \langle q \rangle 
& \stackrel{\sim}{\longrightarrow} & 
E_{q}(K) 
\\ 
u \ {\rm mod} \langle q \rangle 
& \longmapsto & 
\left\{ \begin{array}{ll} 
(X(u,q), Y(u,q)) & (u \not\in \langle q \rangle), \\ 
0                & (u \in \langle q \rangle). 
\end{array} \right.
\end{eqnarray*}
\end{it} 

\underline{\it Proof.} (1) The discriminant $\Delta$ of the Tate curve 
is given by 
\begin{eqnarray*}
& & 
- a_{6}(q) + a_{4}(q)^{2} + 72 a_{4}(q) a_{6}(q) - 64 a_{4}(q)^{3} 
- 432 a_{6}(q)^{2} 
\\
& = & q - 24 q^{2} + \cdots : \ 
\mbox{a formal power series with integral coefficients} 
\\
& \stackrel{\rm in \ fact}{=} & q \prod_{n=1}^{\infty} (1-q^{n})^{24} : 
\ \mbox{a cusp form of weight $12$ for $SL_{2}({\bf Z}).$}
\end{eqnarray*}
Therefore, the Tate curve is smooth over 
${\bf Z}[[q]] \left[ 1 / \Delta \right] = {\bf Z}((q)).$ 

(2) First, note that the following hold: 
\begin{eqnarray*}
\frac{\wp_{L}(z)}{(2 \pi \sqrt{-1})^{2}}  
& = & 
\sum_{n \in {\bf Z}} \frac{q^{n} u}{(1 - q^{n} u)^{2}} 
+ \frac{1}{12} - 2 s_{1}(q) \ \ 
\left( s_{1}(q) \stackrel{\rm def}{=} \sum_{n=1}^{\infty} \sigma_{1}(n) q^{n} 
\right), 
\\ 
\frac{\wp'_{L}(z)}{(2 \pi \sqrt{-1})^{3}} 
& = & 
\sum_{n \in {\bf Z}} \frac{q^{n} u (1 + q^{n} u)}{(1 - q^{n} u)^{3}}. 
\end{eqnarray*}
Because the right hand sides are $q$-series which are invariant under $u \mapsto q u,$ 
hence invariant under $z \mapsto z + 1, z + \tau,$ 
and they have the expansions of $z$ as 
\begin{eqnarray*}
& & 
\frac{u}{(1-u)^{2}} + \frac{1}{12} + 
\sum_{n \neq 0} \frac{q^{n}}{(1 - q^{n})^{2}} - 2 s_{1}(q) + O(z) 
= \frac{1}{\left( 2 \pi \sqrt{-1} z \right)^{2}} + O(z), 
\\
& & 
\frac{u(1+u)}{(1-u)^{3}} + 
\sum_{n \neq 0} \frac{q^{n} \left(1 + q^{n} \right)}{(1 - q^{n})^{3}} + O(z) 
= \frac{-2}{\left( 2 \pi \sqrt{-1} z \right)^{3}} + O(z)  
\end{eqnarray*}
which are equal to those of the left hand sides respectively. 
Therefore, 
\begin{eqnarray*}
X & = & 
\frac{x}{(2 \pi \sqrt{-1})^{2}} - \frac{1}{12} 
= \frac{\wp_{L}(z)}{(2 \pi \sqrt{-1})^{2}} - \frac{1}{12} 
= X(u, q), 
\\
Y & = & 
\frac{y}{2 (2 \pi \sqrt{-1})^{3}} - \frac{x}{2 (2 \pi \sqrt{-1})^{2}} 
+ \frac{1}{24} 
\\ 
& = & 
\frac{\wp'_{L}(z)}{2 (2 \pi \sqrt{-1})^{3}} 
- \frac{\wp_{L}(z)}{2 (2 \pi \sqrt{-1})^{2}} + \frac{1}{24} 
\\
& = & 
\frac{1}{2} \sum_{n \in {\bf Z}} 
\frac{q^{n} u (1 + q^{n} u)}{(1 - q^{n} u)^{3}} - 
\frac{1}{2} \sum_{n \in {\bf Z}} 
\frac{q^{n} u}{(1 - q^{n} u)^{2}} + s_{1}(q) 
\\
& = & 
Y(u, q). 
\end{eqnarray*}
As seen in 2.1, $z + L \mapsto (x = \wp_{L}(z), y = \wp'_{L}(z))$ 
is an isomorphism from ${\bf C}/L$ onto the elliptic curve 
$y^{2} = 4 x^{3} - 60 E_{4}(L) x - 140 E_{6}(L),$ 
and hence 
$$
z + L \mapsto (X = X(u, q), Y = Y(u, q))
$$ 
gives an isomorphism ${\bf C}/L \stackrel{\sim}{\rightarrow} E_{\tau}.$ 

(3) By substituting the variable $q \mapsto q \in K^{\times}$ with $|q| < 1,$ 
$\Delta = q - 24 q^{2} + \cdots$ satisfies that $|\Delta| = |q| \neq 0,$ 
and hence $E_{q}$ is an elliptic curve over $K.$ 
By (2), $X(u, q)$ and $Y(u, q)$ satisfies the equation of the Tate curve: 
$$
Y(u, q)^{2} + X(u, q) Y(u, q) = X(u, q)^{3} + a_{4}(q) X(u, q) + a_{6}(q)
$$
for all complex numbers $u, q$ in a certain convergence domain, 
and hence this equation holds as formal power series in $q$ 
with coefficients in ${\bf Q}(u).$ 
Therefore, by substituting 
the variable $q \mapsto q \in K^{\times}$ with $|q| < 1,$  
one can see that the map in (3) is well-defined, and is evidently injective. 
The addition law on the Tate curve is given by 
\begin{eqnarray*}
& & P_{i} = (x_{i}, y_{i}) \ (i = 1,2,3), \ \ P_{1} + P_{2} = P_{3} 
\\
& \longrightarrow & 
\left\{ \begin{array}{lll} 
(x_{2} - x_{1})^{2} x_{3} & = & 
(y_{2} - y_{1})^{2} + (y_{2} - y_{1})(x_{2} - x_{1}) - 
(x_{2} - x_{1})^{2} (x_{1} + x_{2}), 
\\
(x_{2} - x_{1}) y_{3} & = & 
\left( -(y_{2} - y_{1}) + (x_{2} - x_{1}) \right) x_{3} - 
(y_{1} x_{2} - y_{2} x_{1}), 
\end{array} \right. 
\end{eqnarray*}
if $x_{1} \neq x_{2}.$ 
Hence by (2), this holds if 
$x_{i} = X(u_{i}, q),$ $y_{i} = X(u_{i}, q),$ $(i = 1, 2, 3)$ with 
$u_{1} u_{2} = u_{3}$ for all complex numbers $u_{1}, u_{2}, q$ 
in a certain convergence domain, 
and hence holds as formal power series in $q$ 
with coefficients in ${\bf Q}(u_{1}, u_{2}).$ 
Therefore, by substituting 
the variable $q \mapsto q \in K^{\times}$ with $|q| < 1,$  
one can see that the map in (3) is a homomorphism. 
We omit the surjectivity of the map which is most hardest part of the proof. 
QED. 
\vspace{2ex}

\noindent
\underline{\bf Remark.} 
Similar argument to the proof of Theorem 4.1 (3) is used in [I1] 
to show that $p$-adic theta functions of Mumford curves 
give solutions to soliton equations. 
\vspace{2ex}

\noindent
\underline{\bf Mumford curves.} 
Mumford [Mu2] gave a higher genus version of the Tate curve 
over complete local domains as an analogy of Schottky uniformization theory, 
i.e., for a complete integrally closed noetherian local ring $R$ 
with quotient field $K,$ 
and a Schottky group $\Gamma \subset PGL_{2}(K)$ over $K$ 
which is {\it flat} over $R,$ 
he constructed a {\bf Mumford curve} over $(R \subset) K$ 
which is a proper smooth curve $C_{\Gamma}$ over $K$ obtained 
as the general fiber of a stable curve over $R$ uniformized by $\Gamma$ 
such that its special fiber consists of (may be singular) projective lines 
and its singularities are all $k$-rational $(k$ is the residue field of $R).$ 
Furthermore, he showed that $\Gamma \mapsto C_{\Gamma}$ gives rise to 
the following bijection: 
$$
\begin{array}{ccc}
\left\{ \begin{array}{l} 
\mbox{Conjugacy classes of flat} \\ 
\mbox{Schottky groups over $(R \subset) K$} 
\end{array} \right\} 
& \stackrel{\sim}{\longleftrightarrow} & 
\left\{ \begin{array}{l} 
\mbox{Isomorphism classes of} \\ 
\mbox{Mumford curves over $(R \subset) K$} 
\end{array} \right\} 
\end{array} 
$$
If $K$ is a complete valuation field, 
then any Schottky group $\Gamma$ over $K$ is flat over its valuation ring, 
and it is shown in [GP] that $C_{\Gamma}$ is given as the quotient 
by $\Gamma$ of its region of discontinuity in $K \cup \{ \infty \}$ 
(important examples of rigid analytic geometry). 
\vspace{2ex}

\noindent
\underline{\bf 4.3. Arithmetic Schottky uniformization}
\vspace{2ex}

\noindent
\underline{\bf Stable curves.}
A {\bf stable curve} of genus $g > 1$ over a scheme $S$ is 
a proper and flat morphism $C \rightarrow S$ whose geometric fibers 
are reduced and connected 1-dimensional schemes $C_{s}$ such that 
\begin{itemize}
\item 
$C_{s}$ has only ordinary double points; 
\item 
${\rm Aut}(C_{s})$ is a finite group, i.e., 
if $X$ is a smooth rational component of $C_{s},$ 
then $X$ meets the other components of $C_{s}$ at least 3 points; 
\item 
the dimension of $H^{1}(C_{s}, {\cal O}_{C_{s}})$ is equal to $g.$ 
\end{itemize}

\noindent
\underline{\bf Degenerate curves and dual graphs.} 
A {\bf degenerate curve} is a stable curve whose irreducible components 
are (may be singular) projective lines. 
For a degenerate curve, by the correspondence: 
$$
\begin{array}{rcl}
\mbox{its irreducible components} & \longleftrightarrow & \mbox{vertices} 
\\ 
\mbox{its singular points} & \longleftrightarrow & \mbox{edges} 
\end{array} 
$$
(an irreducible component contains a singular point if and only if 
the corresponding vertex is contained in (or adjacent to) 
the corresponding edge), we have its {\bf dual graph} 
which becomes a stable graph, i.e., a connected and finite graph 
whose vertices has at least $3$ branches ({\bf Figure}). 
For a degenerate curve $C$ with dual graph $\Delta,$ 
\begin{eqnarray*}
\mbox{the genus of $C$} 
& = & 
{\rm rank}_{\bf Z} H_{1}(\Delta, {\bf Z}) 
\\ 
& = & 
\mbox{the number of generators of the free group $\pi_{1}(\Delta).$} 
\end{eqnarray*}

Since any triplet of distinct points on ${\mathbb P}^{1}$ is 
uniquely translated to $(0, 1, \infty)$ by the action of $PGL_{2},$ 
for a stable graph $\Delta,$ 
the moduli space of degenerate curves with dual graph $\Delta$ has dimension 
$$
\sum_{v : \ {\rm vertices \ of} \ \Delta} \left( \deg(v) - 3 \right), 
$$
where $\deg(v)$ denotes the number of branches ($\neq$ edges) 
starting from $v.$ 
In particular, a stable graph is {\bf trivalent}, i.e., 
all the vertices have just $3$ branches if and only if 
the corresponding curves are {\bf maximally degenerate} 
which means that this moduli consists of only one point. 
\vspace{2ex}

\noindent
\underline{\bf Exercise 8.} 
For any stable graph $\Delta,$ prove that 
$$
\sum_{v : \ {\rm vertices \ of} \ \Delta} (\deg(v) - 3) + 
\mbox{the number of edges of $\Delta$} 
\ = \ 3 \left( {\rm rank}_{\bf Z} H_{1}(\Delta, {\bf Z}) - 1 \right). 
$$
\vspace{1ex}

\noindent
\underline{\bf General degenerating process.} 
(Ihara and Nakamura [IhN]). 
For a stable graph $\Delta$ with orientation on each edge, 
\begin{eqnarray*}
g & \stackrel{\rm def}{=} & {\rm rank}_{\bf Z} H_{1}(\Delta, {\bf Z}), 
\\
P_{v} & \stackrel{\rm def}{=} & 
{\mathbb P}^{1}({\bf C}) \ \ \mbox{$(v :$ vertices of $\Delta).$}
\end{eqnarray*}
and for each oriented edge $e$ $(v_{-e} \stackrel{e}{\rightarrow} v_{e})$ 
of $\Delta,$ let
\begin{eqnarray*}
v_{e} & \stackrel{\rm def}{=} & \mbox{the end point of $e,$}
\\
v_{-e} & \stackrel{\rm def}{=} & \mbox{the starting point of $e,$}
\\
\gamma_{e} & : & \mbox{a hyperbolic element of $PGL_{2}({\bf C})$ 
which gives $\gamma_{e} : P_{v_{-e}} \stackrel{\sim}{\rightarrow} P_{v_{e}},$} 
\\ 
t_{e} & \in & P_{v_{e}} : \mbox{the attractive fixed point of $\gamma_{e},$} 
\\ 
t_{-e} & \in & P_{v_{-e}} : \mbox{the repulsive fixed point of $\gamma_{e}.$} 
\end{eqnarray*}
Fix a vertex $v_{0}$ of $\Delta,$ and put 
\begin{eqnarray*}
\Gamma & \stackrel{\rm def}{=} & 
\left\{ \left. \gamma_{e_{1}}^{i_{1}} \cdots \gamma_{e_{n}}^{i_{n}} \ \right| 
\ \mbox{$e_{k} :$ edges, $i_{k} \in \{ \pm 1 \}$ such that 
$e_{n}^{i_{n}} \cdots e_{1}^{i_{1}} \in \pi_{1}(\Delta; v_{0})$} \right\}. 
\end{eqnarray*}
Then under the assumption that the multipliers $s_{e}$ of all $\gamma_{e}$ 
are sufficiently small, 
\begin{itemize}

\item 
$\Gamma$ is a Schottky group of rank $g;$ 

\item 
If $\infty \in \Omega_{\Gamma},$ 
then $\sum_{\gamma \in \Gamma} |\gamma'(z)|$ converges uniformly 
on any compact subset of 
$\Omega_{\Gamma} - \bigcup_{\gamma \in \Gamma} \gamma(\infty);$ 

\item 
$R_{\Gamma} = \Omega_{\Gamma} / \Gamma$ is a Riemann surface of genus $g$ 
obtained from holed Riemann spheres $P_{v}$ $(v :$ vertices of $\Delta)$ 
gluing by $\gamma_{e}$ $(e :$ edges of $\Delta);$ 

\end{itemize}
and hence 
\begin{eqnarray*}
& & s_{e} \rightarrow 0 \ \mbox{$(e$ : edges of $\Delta)$} 
\\
& \Rightarrow & 
R_{\Gamma} \rightarrow 
\mbox{the degenerate curve 
${\displaystyle \left. C_{0} = \left( \bigcup_{v} P_{v} \right) \right/ 
\begin{array}{l} t_{e} = t_{-e} \\ 
\mbox{\footnotesize ($e$ : edges of $\Delta$)} 
\end{array}}$ with dual graph $\Delta.$} 
\end{eqnarray*}
Since ${\mathbb P}^{1}$ has only trivial deformation, 
$R_{\Gamma}$ gives a universal deformation of $C_{0},$ 
and hence varying $t_{\pm e}$ as the {\bf moduli parameters}, 
$s_{e}$ as the {\bf deformation parameters}, 
$R_{\Gamma}$ make an open subset (of dimension $3g -3$ by Exercise 8) 
of the moduli space of curves of genus $g.$ 
\vspace{2ex}

\noindent
\underline{\bf Arithmetic Schottky uniformization.} 
An extension of this process in terms of arithmetic geometry 
(unifying complex geometry and formal geometry over ${\bf Z},$ 
hence rigid geometry) is the following 
{\bf arithmetic Schottky uniformization theory} 
which also gives a higher genus version of the Tate curve: 
\vspace{1ex}

\noindent
\underline{\bf Theorem 4.2.} 
([I3], (1)--(3) were already proved in [IhN] 
for maximally degenerate case without singular components). 
\begin{it} 
Let 
\begin{eqnarray*}
A_{0} & \stackrel{\rm def}{=} & 
\mbox{the coordinate ring of the moduli space 
(i.e., the ring of moduli parameters)} 
\\ 
& & \mbox{over ${\bf Z}$ of degenerate curves with dual graph $\Delta$,} 
\\
A_{\Delta} & \stackrel{\rm def}{=} & 
A_{0} [[ y_{e} \ (e : \mbox{edges of $\Delta$}) ]]. 
\end{eqnarray*}
Then there exists a stable curve $C_{\Delta}$ 
(called the {\it \bfseries generalized Tate curve}) over $A_{\Delta}$ 
of genus $g \stackrel{\rm def}{=} {\rm rank}_{\bf Z} H_{1}(\Delta, {\bf Z})$ 
satisfying:  

{\rm (1)} $C_{\Delta}$ is a universal deformation 
of the universal degenerate curve with dual graph $\Delta$. 

{\rm (2)} By substituting complex numbers $t_{\pm e}$ to the moduli parameters 
and $s_{e} \in {\bf C}^{\times}$ to $y_{e}$ ($e$ are edges of $\Delta$), 
$C_{\Delta}$ becomes a Schottky uniformized Riemann surface 
if $s_{e}$ are sufficiently small. 

{\rm (3)} $C_{\Delta}$ is smooth over $B_{\Delta} = 
A_{\Delta} \left[ 1/y_{e} \ (e : \mbox{edges of} \ \Delta) \right],$ 
and is Mumford uniformized by a Schottky group over $B_{\Delta}$. 
Furthermore, for a complete integrally closed noetherian local ring $R$ 
with quotient field $K$ and a Mumford curve $C$ over $(R \subset) K$ 
such that $\Delta$ is the dual graph of its degenerate reduction, 
there is a ring homomorphism $A_{\Delta} \rightarrow R$ 
gives rise to $C_{\Delta} \otimes_{A_{\Delta}} K \cong C.$ 

{\rm (4)} Using Mumford's theory [Mu3] on degenerating abelian varieties, 
the generalized Jacobian of $C_{\Delta}$ can be expressed as 
$$
{\bf G}_{m}^{g} \left/ 
\langle (p_{ij})_{1 \leq i \leq g} \ \right| \ 1 \leq j \leq g \rangle; 
\ \ 
{\bf G}_{m} \stackrel{\rm def}{=} 
\mbox{the multiplicative algebraic group,} 
$$
where the multiplicative periods $p_{ij}$ of $C_{\Delta}$ 
(called {\it \bfseries universal periods}) 
are given as computable elements of $B_{\Delta}$. 
\end{it}
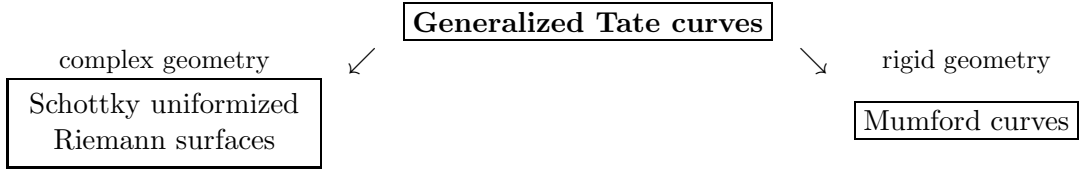
$$
\begin{array}{ccccc}
&  & \fbox{\bf Generalized Tate curves} &  &        \\ 
\mbox{\small complex geometry} & \swarrow &  
& \searrow & \mbox{\small rigid geometry} \\ 
\fbox{$
\begin{array}{c} 
\mbox{Schottky uniformized} \\ \mbox{Riemann surfaces} 
\end{array}
$}
&  &  &  & \fbox{Mumford curves} 
\end{array} 
$$

\underline{\it Sketch of proof.} 
\begin{itemize}

\item 
Step 1 of constructing $C_{\Delta}$ is to give a Schottky group 
$\Gamma_{\Delta}$ over $B_{\Delta}$ 
as in the above {\it general degenerating process,} 
and show that $\Gamma_{\Delta}$ is flat over $A_{\Delta}$ 
(note that this fact together with the result of [Mu2] 
cannot imply the existence of $C_{\Delta}$ 
since $A_{\Delta}$ is not local). 

\item 
Step 2 is, following argument in [Mu2], 
to show that the collection of sets consisting of $3$ fixed points 
in ${\mathbb P}^{1}$ of $\Gamma - \{ 1 \}$ gives rise to a tree 
which is the universal cover of $\Delta$ with covering group $\Delta,$ 
and to construct $C_{\Delta}$ as the quotient by $\Gamma$ 
of the glued scheme of ${\mathbb P}^{1}_{A_{\Delta}}$ 
associated with this tree using Grothendieck's formal existence theorem. 

\item 
In order to give a power series expansion of $p_{ij},$ 
use the infinite product presentation 
by Schottky [S], Manin and Drinfeld [ManD] of the multiplicative periods 
given in Theorem 2.2 (3). 

\end{itemize}

\noindent
\underline{\bf Example.} 
When $\Delta$ consists of one vertex and $g$ loops, 
the universal periods $p_{ij}$ are given in 4.1. 
\vspace{2ex} 

\noindent
\underline{\bf Remark.} 
Denote by 
\begin{eqnarray*}
T_{g} & : & \mbox{the Teichm\"{u}ller space of degree $g,$} 
\\
S_{g} & : & \mbox{the Schottky space of degree $g$} \\
& & \mbox{(the moduli space of Schottky groups with free $g$ generators),} 
\\
H_{g} & : & \mbox{the Siegel upper half space of degree $g.$} 
\end{eqnarray*}
Then 
$$
\begin{array}{ccll}
T_{g} & \stackrel{p}{\longrightarrow} & H_{g} & 
: \mbox{the period map (transcendental)} \\ 
\downarrow & & 
\downarrow \mbox{\small $\exp(2 \pi \sqrt{-1} \cdot)$} &  \\ 
S_{g} & \longrightarrow & H_{g} / {\bf Z}^{g(g+1)/2} & 
: \mbox{computable as power series} \\ 
\downarrow & & \downarrow &  \\ 
{\cal M}_{g}({\bf C}) & \stackrel{\tau}{\longrightarrow} 
& H_{g} / Sp_{2g}({\bf Z}) & 
: \mbox{the Torelli map (algebraic).} 
\end{array}
$$
\vspace{1ex}

\noindent
\underline{\bf Problem.} When any vertex of $\Delta$ has just $3$ branches 
(i.e., the corresponding degenerate curve is maximally degenerate), 
the moduli space of degenerate curves with dual graph $\Delta$ 
consists of one point, and hence $A_{0} = {\bf Z}.$ 
Then express integral coefficients of  
$$
p_{ij} \ \in \ 
A_{\Delta} = {\bf Z} [[y_{e} \ (e : \mbox{edges of $\Delta$})]] 
$$ 
by using some arithmetic functions (cf. [MaT] for the genus $2$ case).

\newpage
\begin{center}
\underline{\large {\bf \S 5. Moduli space of algebraic curves}} 
\end{center}

\noindent
\underline{\bf 5.1. Construction of moduli spaces} 
\vspace{2ex}

A {\bf scheme} is a locally ringed space which is locally given by 
the affine scheme: 
$$
{\rm Spec}(A) \stackrel{\rm def}{=} 
\left\{ \mbox{prime ideals of $A$} \right\} 
\ni {\mathfrak{p}} \mapsto A_{\mathfrak{p}} 
\stackrel{\rm def}{=} 
\left\{ a/s \ | \ a \in A, \ s \in A - {\mathfrak{p}} \right\} 
$$
associated with a commutative ring $A$ with unit $1$ 
(the category of affine schemes is contravariantly equivalent to that 
of commutative rings with unit $1$). 
A scheme $X$ over a scheme $S$ is a scheme with morphism $X \rightarrow S.$ 

The {\bf moduli space of curves} is a space representing 
\begin{center}
the isomorphism classes of curves. 
\end{center}
More precisely, if ${\cal M}_{g}$ is a {\it fine} moduli of curves 
of genus $g,$ 
then 
\begin{eqnarray*}
{\cal M}_{g}(S) & \stackrel{\rm def}{=} & 
\left\{ \mbox{morphisms from $S$ to ${\cal M}_{g}$} \right\} 
\ \ (S: \mbox{schemes}) \\ 
& \stackrel{\rm functorial}{\cong} & 
\left\{ \mbox{isomorphism classes of curves over $S$ of genus $g$} \right\} 
\end{eqnarray*}

\noindent
\underline{\bf Caution!} 
There is no fine moduli as an scheme 
since there are curves with nontrivial automorphism 
(for example, {\bf hyperelliptic curves} defined by $y^{2} = f(x)$ 
has a nontrivial automorphism $x \mapsto x,$ $y \mapsto -y$). 
Because if ${\cal M}_{g}$ is a fine moduli scheme, 
then the identity map on ${\cal M}_{g}$ corresponds to 
a curve ${\cal C}$ over ${\cal M}_{g}$ which is universal, 
i.e., for each scheme $S,$
$$
\begin{array}{ccc}
\left\{ \mbox{morphisms from $S$ to ${\cal M}_{g}$} \right\} 
& \longleftrightarrow & 
\left\{ \mbox{curves over $S$ with morphisms to ${\cal C}$} \right\} 
\\ 
S \rightarrow {\cal M}_{g} & \longmapsto & 
{\cal C} \times_{{\cal M}_{g}} S \rightarrow {\cal C} 
\end{array}
$$
Therefore, any automorphism on a curve over $S$ must be the identity map.  
\vspace{2ex}

\noindent
\underline{\bf Solutions.}
\begin{itemize}

\item[{\bf (S1)}] Construct the fine moduli as an scheme by considering 
additional structures on curves. 

\item[{\bf (S2)}] Taking the categorical quotient of the above fine moduli, 
construct the fine moduli as an {\bf algebraic stack}, 
the scheme-theoretic analog of {\bf orbifolds}, which is represented as 
$$
\left[ U/R \right] : \ \mbox{the quotient of $U$ by $R,$} 
$$
where $U, R$ are schemes with \'{e}tale morphisms $s, t : R \rightarrow U$ 
and a morphism $\mu : R \times_{U,t,s} R \rightarrow R$ such that 
$(s, t) : R \rightarrow U \times U$ is finite 
and $s, t, \mu$ form a groupoid. 
For a scheme $S,$ $\left[ U/R \right](S) = {\rm Hom}(S, U/R)$ 
is the category given by 
\begin{eqnarray*}
{\rm Ob} \left( \left[ U/R \right](S) \right) 
& \stackrel{\rm def}{=} & 
{\rm Hom}(S, U), 
\\
{\rm Mor} \left( \left[ U/R \right](S) \right) 
& \stackrel{\rm def}{=} & 
\left\{ \mbox{$\alpha \in {\rm Hom}(S, R)$ giving 
$s \circ \alpha \stackrel{\sim}{\rightarrow} t \circ \alpha$} \right\}, 
\end{eqnarray*}
and $R$ gives the equivalence relation by $\mu.$  

\item[{\bf (S3)}] Taking the geometric quotient of the above fine moduli, 
construct the {\it coarse} moduli as an scheme. 

\end{itemize}
\vspace{2ex}

\noindent
\underline{\bf Moduli of elliptic curves.} 
\begin{itemize}

\item {\bf Analytic construction:} 

\noindent 
{\bf (S1)} If $E$ is an elliptic curve over ${\bf C},$ 
and $\iota$ is an isomorphism 
${\bf Z}^{\oplus 2} \stackrel{\sim}{\rightarrow} H_{1}(E, {\bf Z})$ 
such that $\iota$ is {\it canonical}, i.e., 
$\iota(\mbox{\boldmath $e$}_{1}), \iota(\mbox{\boldmath $e$}_{2})$ 
intersects as the $x, y$-axes, 
then the ratio 
$$
\left. \left( \int_{\iota(\mbox{\boldmath $e$}_{2})} \omega \right) \right/ 
\left( \int_{\iota(\mbox{\boldmath $e$}_{1})} \omega \right) 
$$
is independent of $\omega \in H^{1}(E, \Omega_{E}) - \{ 0 \}$ 
and belongs to the Poincar\'{e} upper half plane $H_{1}.$ 
Therefore, by the correspondence: 
$$
H_{1} \ni \tau \ \leftrightarrow \ 
\left( {\bf C}/({\bf Z} + {\bf Z} \tau); \ 
\iota(\mbox{\boldmath $e$}_{1}) = 1, \iota(\mbox{\boldmath $e$}_{2}) = \tau 
\right), 
$$ 
$H_{1}$ becomes the fine moduli space of elliptic curves $E$ over ${\bf C}$ 
with canonical isomorphism 
${\bf Z}^{\oplus 2} \stackrel{\sim}{\rightarrow} H_{1}(E, {\bf Z}).$ 
\vspace{1ex}

\noindent
{\bf (S2)} By (S1), the fine moduli stack of elliptic curves over ${\bf C}$ 
is given by the complex analytic stack, i.e. {\bf orbifold} 
$$
\left[ H_{1} / SL_{2}({\bf Z}) \right]. 
$$
{\bf (S3)} Since 
\begin{eqnarray*}
& & 
\mbox{the elliptic curves 
$y_{i}^{2} = 4 x_{i}^{3} - \alpha_{i} x_{i} - \beta_{i}$ $(i = 1,2)$ 
over ${\bf C}$ are isomorphic} 
\\
& \Leftrightarrow & 
\mbox{there are $a, b, c, d, e \in {\bf C}$ with $a, c \neq 0$ such that}
\\
& & 
\left\{ \begin{array}{ll} 
\mbox{$x_{2} = a x_{1} + b $ : order $\geq -2$ at the origin,} 
\\
\mbox{$y_{2} = c y_{1} + d x_{1} + e$ : order $\geq -3$ at the origin} 
\end{array} \right. 
\\
& \Leftrightarrow & 
\mbox{there are $a, c \in {\bf C}^{\times}$ such that 
$a^{3} = c^{2}, x_{2} = a x_{1}, y_{2} = c y_{1}$} 
\\
& \Leftrightarrow & 
\mbox{the {\bf {\boldmath $j$}-invariants} 
${\displaystyle \frac{\alpha_{i}^{3}}{\alpha_{i}^{3} - 27 \beta_{i}^{2}}}$ of 
$y_{i}^{2} = 4 x_{i}^{3} - \alpha_{i} x_{i} - \beta_{i}$ $(i = 1,2)$ 
are equal} 
\\
& & 
\left(\mbox{note that $\alpha_{i}^{3} - 27 \beta_{i}^{2} \neq 0$} \right), 
\end{eqnarray*}
the coarse moduli scheme of elliptic curves over ${\bf C}$ 
becomes the affine line over ${\bf C},$ 
and the $j$-function 
\begin{eqnarray*}
j(\tau) 
& \stackrel{\rm def}{=} & 
\frac{(60 E_{4}({\bf Z} + {\bf Z} \tau))^{3}}
{(60 E_{4}({\bf Z} + {\bf Z} \tau))^{3} - 
27 (140 E_{6} ({\bf Z} + {\bf Z} \tau))^{2}} 
\\
& = & 
\frac{1}{1728} \left( \frac{1}{q} + 744 + 196884 q + 21493760 q^{2} + \cdots 
\right) \ \left( q \stackrel{\rm def}{=} e^{2 \pi \sqrt{-1} \tau} \right) 
\end{eqnarray*} 
gives a biholomorphic map from the geometric quotient 
$H_{1} / SL_{2}({\bf Z})$ onto ${\bf C}.$ 

\item {\bf Algebraic construction:} 

\noindent
{\bf (S1)} For complex numbers 
$\mu \neq 1, \ \zeta_{3} = e^{2 \pi \sqrt{-1}/3}, \ \zeta_{3}^{2},$ 
put
\begin{eqnarray*}
E(\mu) & \stackrel{\rm def}{=} & 
\left\{ (x_{0}:x_{1}:x_{2}) \in {\mathbb P}^{1}({\bf C}) \ | \ 
x_{0}^{3} + x_{1}^{3} + x_{2}^{3} = 3 \mu x_{0} x_{1} x_{2} \right\} : 
\ \mbox{\bf Hesse's cubic}
\\
& : & \mbox{an elliptic curve over ${\bf C}$ with origin $(1:-1:0)$ 
containing}
\\ 
& & \mbox{$3$-division points} \ (1:-\beta:0), (0:1:-\beta), (-\beta:0:1) \ 
\left( \beta = 1, \zeta_{3}, \zeta_{3}^{2} \right) 
\end{eqnarray*}
Then $\mu \mapsto (E(\mu)$ with the $3$-division points) 
gives a bijection: 
\begin{eqnarray*}
{\bf C} - \left\{ 1, \zeta_{3}, \zeta_{3}^{2} \right\} 
& \stackrel{\sim}{\rightarrow} & 
\left\{ \begin{array}{l} 
\mbox{isomorphism classes of elliptic curves over ${\bf C}$} 
\\ 
\mbox{with symplectic level $3$ structure}
\\
\mbox{$\left( {\bf Z}/ 3{\bf Z} \right)^{\oplus 2} 
\stackrel{\sim}{\rightarrow}  E[3] \stackrel{\rm def}{=} 
\{ P \in E \ | \ 3P = 0 \} 
$} 
\end{array} \right\} 
\\
& \cong & 
H_{1} / \Gamma(3), 
\end{eqnarray*}
where $\Gamma(3)$ denotes the principal congruence subgroup 
of $SL_{2}({\bf Z})$ of level $3.$ 
Therefore, ${\bf C} - \left\{ 1, \zeta_{3}, \zeta_{3}^{2} \right\}$ 
has a natural model over ${\bf Z} \left[ 1/3, \zeta_{3} \right]$ 
as the fine moduli scheme of elliptic curves with level $3$ structure, 
and this can be compactified to ${\mathbb P}^{1}$ by adding the $4$ points 
$1, \zeta_{3}, \zeta_{3}^{2}, \infty$ which correspond degenerate curves. 
Nakamura [N] gave this higher dimensional version, i.e., 
constructed a compactification of the moduli of principally polarized 
abelian varieties with level structure as an moduli space. 
\vspace{1ex}

\noindent
{\bf (S2)} The fine moduli stack over 
${\bf Z} \left[ 1/3, \zeta_{3} \right]$ is given by the quotient stack 
of the above model in (S1) by $SL_{2} \left( {\bf Z} / 3 {\bf Z} \right).$ 

\item {\bf Scheme theoretic construction:} 

\noindent 
{\bf (S1)} If $E$ is an elliptic curve over a scheme $S$ with $0$-section 
$e : S \rightarrow E,$ 
then 
$$
H^{0} \left( E, {\cal O}(3 \cdot e(S)) \right) : \ 
\mbox{the space of generalized elliptic functions} 
$$ 
defines an embedding $E \hookrightarrow {\mathbb P}^{2}_{S}.$ 
Therefore, by the theory of Hilbert schemes, 
there is a fine moduli scheme ${\cal H}_{1}$ over ${\bf Z}$ classifying 
elliptic curves with embedding into ${\mathbb P}^{2}_{S}$ as above. 
\vspace{1ex}

\noindent
{\bf (S2)} The fine moduli stack over ${\bf Z}$ of elliptic curves is given by 
the quotient stack 
$$
\left[ {\cal H}_{1} / {\rm Aut}({\mathbb P}^{2}) \right] = 
\left[ {\cal H}_{1} / PGL_{3} \right]; \ \ 
PGL_{n} \stackrel{\rm def}{=} GL_{n}/{\bf G}_{m}. 
$$
\end{itemize}
\vspace{1ex}

\noindent
\underline{\bf Exercise 9.} 
Prove that for $\tau \in H_{1},$ 
\begin{eqnarray*}
& & 
\left\{ \gamma \in SL_{2}({\bf Z}) \ | \ \gamma(\tau) = \tau \right\} 
\\
& = & 
\left\{ \begin{array}{lll}
{\displaystyle \left\langle \rho 
\left( \begin{array}{rr} 1 & 1 \\ -1 & 0 \end{array} \right) 
\rho^{-1} \right\rangle } & 
\mbox{: order $6$} & 
\mbox{(if $\exists \rho \in SL_{2}({\bf Z})$ such that 
$\tau = \rho(\zeta_{3})),$} 
\\
& & 
\\
{\displaystyle \left\langle \rho 
\left( \begin{array}{rr} 0 & 1 \\ -1 & 0 \end{array} \right) 
\rho^{-1} \right\rangle } & 
\mbox{: order $4$} & 
\mbox{(if $\exists \rho \in SL_{2}({\bf Z})$ such that 
$\tau = \rho(\sqrt{-1})),$} 
\\
& & 
\\
{\displaystyle \left\langle 
\left( \begin{array}{rr} -1 & 0 \\ 0 & -1 \end{array} \right) \right\rangle} & 
\mbox{: order $2$} & 
\mbox{(otherwise).} 
\end{array} \right.
\end{eqnarray*}  
\vspace{1ex}

\noindent
\underline{\bf Construction of moduli for genus {\boldmath$> 1.$}}
There are 3 approaches using 
\vspace{2ex}

\noindent
{\bf 1. Teichm\"{u}ller theory:} 
Fix a Riemann surface $R_{0}$ of genus $g > 1.$ 
Then the {\bf Teichm\"{u}ller space} of degree $g$ is defined by 
\begin{eqnarray*}
T_{g} & \stackrel{\rm def}{=} & 
\left\{ (R, h) \ \left| \ 
\begin{array}{l} 
\mbox{$R$ : Riemann surfaces of genus $g$} \\ 
\mbox{$h$ : orientation preserving diffeomorphisms $R_{0} \rightarrow R$} 
\end{array} \right\} \right/ \sim 
\\ 
& ; & (R, h) \sim (R', h') \stackrel{\rm def}{\Longleftrightarrow} 
h' \circ h^{-1} \ \mbox{is homotopic to a biholomorphic map,} 
\end{eqnarray*}
and the {\bf Teichm\"{u}ller modular group} or {\bf mapping class group} of 
degree $g$ is defined by 
\begin{eqnarray*}
\Pi_{g} & \stackrel{\rm def}{=} & 
\left\{ \mbox{homotopy classes of orientation preserving 
diffeomorphisms $R_{0} \rightarrow R_{0}$} \right\} 
\end{eqnarray*}
which acts on $T_{g}$ as $\mu (R, h) = (R, h \circ \mu)$ $(\mu \in \Pi_{g})$ 
properly discontinuously. 
Then the quotient orbifold $\left[ T_{g} / \Pi_{g} \right]$ exists 
and becomes the moduli space of Riemann surfaces of genus $g.$ 
Teichm\"{u}ller proved that $T_{g}$ is homeomorphic to ${\bf R}^{6g-6}$ 
and becomes naturally a complex manifold of dimension $3g - 3$ 
by using the theory of quasiconformal maps (see [IT]). 
Since $T_{g}$ is connected and simply connected, 
$$
\pi_{1} \left( \left[ T_{g} / \Pi_{g} \right] \right) \ \cong \ \Pi_{g}, 
$$
and these are canonically isomorphic to 
$$
\left. {\rm Aut}^{+} \left( \pi_{1}(R_{0}) \right) \right/ 
{\rm Inn} \left( \pi_{1}(R_{0}) \right), 
$$
where ${\rm Aut}^{+} \left( \pi_{1}(R_{0}) \right)$ denotes 
the automorphism group of $\pi_{1}(R_{0})$ preserving 
the (alternating and bilinear) intersection form on 
$H_{1} \left( R_{0}, {\bf Z} \right) = 
\pi_{1}(R_{0}) / \left[ \pi_{1}(R_{0}), \pi_{1}(R_{0}) \right],$ 
and ${\rm Inn} \left( \pi_{1}(R_{0}) \right)$ denotes 
the inner automorphism group. 
\vspace{2ex}

\noindent
\underline{\bf Caution!} 
Royden showed that if $g > 1,$ then ${\rm Aut}(T_{g}) = \Pi_{g},$ 
and hence the $T_{g}$ is not a homogeneous space. 
Therefore, one cannot regard the Teichm\"{u}ller modular group 
as a discrete subgroup of a Lie group. 
\vspace{2ex}

\noindent
{\bf 2. Moduli theory of abelian varieties:} 
A {\bf principally polarized abelian variety} $(A, \varphi)$ 
is a pair of an abelian variety $A,$ i.e., 
a proper (commutative) algebraic group and 
an isomorphism $A \rightarrow \widehat{A}$ 
(: the dual abelian variety of $A)$ 
induced from an ample divisor on $A.$ 
There exists a moduli space ${\cal A}_{g}$ of 
principally polarized $g$-dimensional abelian varieties, 
and 
$$
{\cal A}_{g}({\bf C}) \ \cong \ 
\left[ H_{g} \left/ Sp_{2g}({\bf Z}) \right. \right]. 
$$
Here 
\begin{eqnarray*}
H_{g} & \stackrel{\rm def}{=} & 
\left\{ Z \in M_{g}({\bf C}) \ | \ Z : \mbox{symmetric}, \ 
{\rm Im}(Z) > 0 \right\} 
\\ 
& : & \mbox{the {\bf Siegel upper half space} of degree $g,$} 
\\ 
Sp_{2g}({\bf Z}) & \stackrel{\rm def}{=} & 
\left\{ G \in M_{2g}({\bf Z}) \ \left| \ 
G \left( \begin{array}{cc} 0 & E_{g} \\ -E_{g} & 0 \end{array} \right) 
\mbox{}^{t}G = 
\left( \begin{array}{cc} 0 & E_{g} \\ -E_{g} & 0 \end{array} \right) 
\right. \right\} 
\\ 
& : & \mbox{the {\bf integral symplectic group} of degree $g$ over ${\bf Z}$} 
\\
& & \mbox{acts on $H_{g}$ as $Z \mapsto (AZ+B)(CZ+D)^{-1}$ for 
${\displaystyle G = 
\left( \begin{array}{cc} A & B \\ C & D \end{array} \right),}$} 
\end{eqnarray*}
and $Z / \sim \ \in H_{g} / Sp_{2g}({\bf Z})$ corresponds to 
the pair of an abelian variety ${\bf C}^{g} / L,$ 
where $L = {\bf Z}^{g} + {\bf Z}^{g} \cdot Z$ 
is the lattice in ${\bf C}^{g}$ generated by 
the unit vectors $\mbox{\boldmath$e$}_{i}$ and 
the $i$-th row vectors $\mbox{\boldmath$z$}_{i}$ of $Z),$ 
and the polarization associated with the alternating bilinear form 
$\psi$ on $L \times L$ such that 
$$
\psi(\mbox{\boldmath $e$}_{i}, \mbox{\boldmath $e$}_{j}) = 
\psi(\mbox{\boldmath $z$}_{i}, \mbox{\boldmath $z$}_{j}) = 0, \ \ 
\psi(\mbox{\boldmath $e$}_{i}, \mbox{\boldmath $z$}_{j}) = \delta_{ij}. 
$$
By Torelli's theorem, by the correspondence: 
\begin{eqnarray*} 
& & \mbox{proper smooth curves $C$} 
\\
& \longmapsto & 
\mbox{their Jacobian varieties ${\rm Jac}(C)$ with principally polarization} 
\\ 
& & 
\mbox{induced from the theta divisor 
$\left\{ P_{1} + \cdots + P_{g-1} - (g-1) P_{0} \ | \ P_{i} \in C \right\},$} 
\end{eqnarray*}
the (coarse) moduli of proper smooth curves is realized as 
a subvariety of ${\cal A}_{g}.$ 
This fact gives rise to the {\bf Schottky problem} 
which means to characterize Jacobian varieties 
among general abelian varieties, 
or to describe explicitly the subvariety of ${\cal A}_{g}$ 
consisting of Jacobian varieties. 
\vspace{2ex}

\noindent
{\bf 3. Geometric invariant theory:} 
For a proper smooth curve $C$ over $S$ of genus $g > 1,$ 
the spaces $H^{0} \left( C_{s}, \Omega_{C_{s}}^{\otimes 3} \right)$ $(s \in S)$ 
have dimension $5(g - 1)$ by Riemann-Roch's theorem, 
and give an embedding $C \hookrightarrow {\mathbb P}_{S}^{5g - 6}.$ 
Then by the theory of Hilbert schemes, 
there exists a fine moduli scheme ${\cal H}_{g}$ over ${\bf Z}$ 
classifying tricanonically embedded curves 
$C \hookrightarrow {\mathbb P}_{S}^{5g - 6},$ 
and hence the quotient stack 
$$
{\cal M}_{g} \ \stackrel{\rm def}{=} \ 
\left[ {\cal H}_{g} / PGL_{5g - 5} \right] 
$$
is the fine moduli space of proper smooth curves of genus $g.$ 
Since $PGL_{5g - 6}$ is smooth and the functor 
$S \mapsto {\rm Isom}_{S}(C, C')$ is represented by 
a finite and unramified scheme over $S$ for curves $C, C'$ over $S,$ 
by an \'{e}tale slice argument, 
${\cal M}_{g}$ becomes an algebraic stack. 
Furthermore, by showing that each point on ${\cal H}_{g}$ is {\it stable} 
under the action of $PGL_{5g - 5},$ 
it follows from {\bf geometric invariant theory} [FKM] by Mumford 
that the geometric quotient ${\cal H}_{g} / PGL_{5g - 5}$ exists and 
gives the coarse moduli scheme of proper smooth curves of genus $g.$ 
\vspace{2ex}

\noindent
\underline{\bf Dictionary for the moduli stack.} 
In what follows, 
$$
\mbox{\boldmath ${\cal M}_{g} \ \stackrel{\rm def}{=}$ \ 
{\bf the moduli stack over Z of proper smooth curves of genus} 
$g > 1.$}
$$
Then 
$$
{\cal M}_{g}({\bf C}) \ = \ \mbox{the quotient orbifold 
$\left[ T_{g} / \Pi_{g} \right],$}
$$
and for schemes (more generally algebraic stacks) $S,$ 
\begin{eqnarray*}
& & 
{\cal M}_{g}(S) = \mbox{the category of proper smooth curves over $S$ 
of genus $g$}
\\
& \Rightarrow & 
\mbox{The identity map on ${\cal M}_{g}$ gives 
the {\bf universal curve} ${\cal C}$ over ${\cal M}_{g}.$} 
\end{eqnarray*}
Furthermore, 
\begin{eqnarray*}
& & \mbox{an object $\alpha$ on (over) ${\cal M}_{g}$} 
\\
& \Longleftrightarrow & 
\mbox{a system $\{ \alpha_{S} \}$ of objects on $S$ 
for proper smooth curves over $S$ of genus $g$} 
\\ 
& & \mbox{such that $\{ \alpha_{S} \}$ are functorial for $S.$}
\end{eqnarray*}

\noindent
\underline{\bf Dimension of the moduli.}
\begin{itemize}

\item {\bf Analytic method:} 
Since ${\rm Aut}(H_{1}) = PSL_{2}({\bf R}),$ by the theory of Fuchsian models, 
for Riemann surfaces $R = H_{1} / \pi_{1}(R),$ $R' = H_{1} / \pi_{1}(R')$ 
of genus $g > 1,$ 
\begin{eqnarray*}
& & 
\left( R; \pi_{1}(R) \hookrightarrow PSL_{2}({\bf R}) \right) 
\cong 
\left( R; \pi_{1}(R') \hookrightarrow PSL_{2}({\bf R}) \right) 
\\
& \Longleftrightarrow & 
\mbox{$\pi_{1}(R)$ and $\pi_{1}(R')$ are conjugate in $PSL_{2}({\bf R}).$} 
\end{eqnarray*}
Therefore, under fixing a Riemann surface $R_{0}$ of genus $g,$ 
$$
T_{g} \cong 
\left\{ \begin{array}{l} 
\mbox{conjugacy classes of injective homomorphisms} 
\\
\mbox{$\iota : \pi_{1}(R_{0}) \rightarrow PSL_{2}({\bf R})$ satisfying that}
\\
\mbox{$H_{1} / \iota(\pi_{1}(R_{0}))$ are Riemann surfaces of genus $g$} 
\end{array} \right\}, 
$$
and the real dimension of the right hand side is 
\begin{eqnarray*}
& & \dim_{\bf R} (PSL_{2}({\bf R}) \times 
\left( \sharp \{ \mbox{generators of $\pi_{1}(R_{0})$} \} - 
\sharp \{ \mbox{relations in $\pi_{1}(R_{0})$} \} - 1 \right) 
\\
& = & 6g - 6. 
\end{eqnarray*}

Furthermore, 
under the assumption that for Schottky uniformized Riemann surfaces $R, R'$ 
of genus $g >1,$ 
\begin{eqnarray*}
& & 
\left( R = \Omega_{\Gamma}/\Gamma; 
\Gamma \hookrightarrow PGL_{2}({\bf C}) \right) 
\cong 
\left( R' = \Omega_{\Gamma'}/\Gamma'; 
\Gamma' \hookrightarrow PGL_{2}({\bf C}) \right) 
\\
& \stackrel{\mbox{may be}}{\Longleftrightarrow} & 
\mbox{$\Gamma$ and $\Gamma'$ are conjugate in $PGL_{2}({\bf C}),$} 
\end{eqnarray*}
by letting $F_{g}$ be the free group of rank $g,$ 
we have 
$$
{\cal M}_{g}({\bf C}) \cong 
\left. \left\{ \begin{array}{l} 
\mbox{conjugacy classes of injective homomorphisms} 
\\
\mbox{$\iota : F_{g} \rightarrow PGL_{2}({\bf C})$ satisfying that}
\\
\mbox{$\iota(F_{g})$ are Schottky groups} 
\end{array} \right\} \right/ {\rm Aut}(F_{g}), 
$$
and the complex dimension of the right hand side is 
$$
\dim_{\bf C} (PGL_{2}({\bf C}) \times 
\left( \sharp \{ \mbox{generators of $F_{g}$} \} - 1 \right) 
= 3g - 3. 
$$

\item {\bf Algebraic method (deformation theory [HM]):} 
For a field $k,$ 
\begin{eqnarray*}
A_{0} 
& \stackrel{\rm def}{=} & 
k[\varepsilon]/(\varepsilon^{2}), 
\\
C 
& : & 
\mbox{a proper smooth curve over $k$ of genus $g > 1,$} 
\\
\{ U_{\alpha} \} 
& : & 
\mbox{an affine open cover of $C,$} 
\end{eqnarray*}
and let $\varphi_{\alpha \beta}$ be a first-order infinitesimal deformation 
of $C,$ i.e., $A_{0}$-linear ring homomorphisms 
$$ 
{\cal O}_{U_{\alpha} \times {\rm Spec}(A_{0})} 
|_{(U_{\alpha} \cap U_{\beta})} \rightarrow 
{\cal O}_{U_{\beta} \times {\rm Spec}(A_{0})} |_{(U_{\alpha} \cap U_{\beta})}
$$
satisfying that 
$$
\left\{ \begin{array}{l} 
\mbox{$\varphi_{\alpha \gamma} = 
\varphi_{\beta \gamma} \circ \varphi_{\alpha \beta}$ on 
$U_{\alpha} \cap U_{\beta} \cap U_{\gamma}$ (: the cocycle condition),} 
\\
\mbox{$\varphi_{\alpha \beta}|_{(U_{\alpha} \cap U_{\beta}) 
\times {\rm Spec}(k)}$ is the identity.}
\end{array} \right. 
$$ 
Then the $k$-linear homomorphisms 
$D_{\alpha \beta} : {\cal O}_{(U_{\alpha} \cap U_{\beta})} \rightarrow 
{\cal O}_{(U_{\alpha} \cap U_{\beta})}$ 
given by $\varphi_{\alpha \beta}(f) = f + \varepsilon D_{\alpha \beta}(f)$ 
satisfies that 
$$
D_{\alpha \beta}(f \cdot g) = f \cdot D_{\alpha \beta}(g) + 
g \cdot D_{\alpha \beta}(f), 
\ \ 
D_{\alpha \gamma}(f) = D_{\beta \gamma}(f) \cdot D_{\alpha \beta}(f), 
$$ 
and hence $\{ D_{\alpha \beta} \}$ defines an element of 
the first cohomology group $H^{1}(C, {\cal T}_{C})$ 
of the tangent bundle ${\cal T}_{C}$ on $C.$ 
Since $\dim_{k}(C) = 1,$ the obstruction space is 
$H^{2}(C, {\cal T}_{C}) = \{ 0 \},$ 
and hence the tangent space of ${\cal M}_{g} \otimes_{\bf Z} k$ 
at the point $[C]$ corresponding to $C$ is isomorphic to 
$H^{1}(C, {\cal T}_{C}).$ 
Therefore, 
\begin{eqnarray*}
& & \mbox{the dimension of the tangent space of 
${\cal M}_{g} \otimes_{\bf Z} k$ at $[C]$} 
\\
& = & 
\dim_{k} H^{1} \left( C, {\cal T}_{C} \right) 
\\
& = & 
\dim_{k} H^{0} \left( C, \Omega_{C}^{\otimes 2} \right) \ \ 
\mbox{(by Serre's duality)} 
\\
& = & 3g - 3 \ \ 
\mbox{(by Riemann-Roch's theorem and that $\deg(\Omega_{C}) = 2g - 2 > 0$).} 
\end{eqnarray*}
\end{itemize}

\noindent
\underline{\bf Remark.} 
For proper smooth curves $C,$ 
$$
H^{1}(C, {\cal T}_{C}) \cong {\rm Ext}^{1}({\cal O}_{C}, {\cal T}_{C}) 
\cong {\rm Ext}^{1}(\Omega_{C}, {\cal O}_{C}), 
$$
and the last group also classifies first-order infinitesimal deformations 
of stable curves. 
\vspace{2ex}

\noindent
\underline{\bf 5.2. Stable curves and their moduli space}
\vspace{2ex}

\noindent
\underline{\bf Stable curves.}
Recall that a {\bf stable curve} of genus $g > 1$ over a scheme $S$ is defined to be 
a proper and flat morphism $C \rightarrow S$ whose geometric fibers 
are reduced and connected 1-dimensional schemes $C_{s}$ such that 
\begin{itemize}
\item 
$C_{s}$ has only ordinary double points; 
\item 
${\rm Aut}(C_{s})$ is a finite group, i.e., 
if $X$ is a smooth rational component of $C_{s},$ 
then $X$ meets the other components of $C_{s}$ at least 3 points; 
\item 
the dimension of $H^{1}(C_{s}, {\cal O}_{C_{s}})$ is equal to $g.$ 
\end{itemize}
For a stable curve $C$ over $S$ (may not be smooth), 
it is useful to consider the {\bf dualizing sheaf} 
(or canonical invertible sheaf) $\omega_{C/S}$ on $C$ 
which is defined as the following conditions: 
\begin{itemize}

\item 
$\omega_{C/S}$ is functorial on $S;$ 

\item 
if $S = {\rm Spec}(k)$ $(k$ is an algebraically closed field), 
$f : C' \rightarrow C$ be the normalization (resolution) of $C,$ 
$x_{1},..., x_{n},$ $y_{1},..., y_{n},$ are the points of $C'$ 
such that $z_{i} = f(x_{i}) = f(y_{i})$ $(1 \leq i \leq n)$ 
are the ordinary double points on $C,$ 
then $\omega_{C/S}$ is the sheaf of $1$-forms $\eta$ on $C'$ 
which are regular except for simple poles at $x_{i}, y_{i}$ 
such that 
$$
{\rm Res}_{x_{i}}(\eta) + {\rm Res}_{y_{i}}(\eta) = 0.
$$ 

\end{itemize}
Then it is shown by Rosenlicht and Hartshorne that 
$\omega_{C/S}$ is a line bundle on $C,$ 
Riemann-Roch's theorem holds for the canonical divisor 
corresponding to $\omega_{C},$ and 
$$
\dim H^{1}(C_{s}, {\cal O}_{C_{s}}) \ = \ \dim H^{0}(C_{s}, \omega_{C_{s}}). 
$$

\noindent
\underline{\bf Theorem 5.1.} (Deligne and Mumford [DM]) 
\begin{it} 
There exists the fine moduli space $\overline{\cal M}_{g}$ 
(called {\it \bfseries Deligne-Mumford's compactification} of ${\cal M}_{g}$) 
as an algebraic stack over {\bf Z} 
classifying stable curves of genus $g > 1.$ 
$\overline{\cal M}_{g}$ is proper smooth over ${\bf Z},$ 
and contains ${\cal M}_{g}$ as its open dense substack.
\end{it}
\vspace{2ex}

\underline{\it Sketch of proof.} 
The construction of $\overline{\cal M}_{g}$ is similar to that of 
${\cal M}_{g}$ by replacing $\Omega_{C}$ with dualizing sheaves $\omega_{C}.$ 
The properness of $\overline{\cal M}_{g}$ follows from the valuative criterion 
and the {\bf stable reduction theorem}: 
Let $R$ be a discrete valuation ring 
with quotient field $K,$ 
and let $C$ be a proper and smooth curve over $K$ of genus $g > 1.$ 
Then there exists a finite extension $L$ of $K$ and a stable curve 
${\cal C}$ over the integral closure $R_{L}$ of $R$ in $L$ 
such that ${\cal C} \otimes_{R_{L}} L \cong C \otimes_{K} L.$
\vspace{2ex}

\noindent
\underline{\bf Irreducibility of the moduli.}
\vspace{2ex}

As an application of Theorem 5.1, Deligne and Mumford [DM] 
proved the irreducibility of any geometric fibers of $\overline{\cal M}_{g}$ 
by applying Enriques-Zariski's connectedness theorem 
to the proper and smooth stack $\overline{\cal M}_{g}$ over ${\bf Z}$ 
whose fiber over ${\bf C}$ is connected (by Teichm\"{u}ller's theory). 
Therefore, 
\begin{center}
{\bf Any geometric fiber of {\boldmath ${\cal M}_{g}$} is irreducible.} 
\end{center}
This fact is essentially used in 6.3 to study automorphic forms 
on the moduli of curves. 
\vspace{2ex}

\noindent
\underline{\bf 5.3. Intersection theory on the moduli space} 
\vspace{2ex}

A cycle class in an algebraic variety $X$ is defined to be a rational equivalence class 
of ${\bf Z}$-linear finite sums of subvarieties of $X$, 
and the Chow ring ${\rm CH}^{*}(X)$ denotes the group of cycle classes in $X$ 
whose ring structure is given by intersection products. 
The structure of ${\rm CH}^{*} \left( {\cal M}_{g} \right)$ 
is an important subject in algebraic geometry and mathematical physics, 
and was studied by Mumford, Witten, Kontsevich, Faber, Mirzakhani and others. 
A basic tool to this study is Grothendieck-Riemann-Roch's theorem 
for families of algebraic curves. 
\vspace{2ex}

\noindent
\underline{\bf Grothendieck-Riemann-Roch's theorem (GRR).} 
This theorem states the following: 
If $\pi : X \rightarrow B$ is a proper smooth morphism over a smooth base, 
and $E$ is a coherent sheaf on $X,$ 
then 
$$
{\rm ch} \left( \pi_{!}(E) \right) = 
\pi_{\ast} \left( {\rm ch}(E) \cdot {\rm td} \left( {\cal T}_{C/B} \right) \right) 
$$ 
in the Chow ring ${\rm CH}^{\ast}(B) \otimes_{\bf Z} {\bf Q}$ 
with ${\bf Q}$-coefficients, 
where ${\rm ch}$ denotes the exponential Chern character, 
and ${\cal T}_{C/B} = \Omega_{C/B}^{\otimes (-1)}$ 
denotes the tangent bundle on $C$ over $B$. 
In order to apply this theorem to a proper smooth curve $\pi : C \rightarrow B$ 
of genus $g > 1$, and $E = \Omega_{C/B}^{n}$ $(n \geq 1)$,  
put $\gamma = c_{1}(\Omega_{C/B})$. 
Then 
$$
{\rm ch} \left( \pi_{\ast} \left( \Omega_{C/B}^{n} \right) \right) = 
\pi_{\ast} \left( \left( 1 + \gamma + \frac{\gamma^{2}}{2} + \cdots \right)^{n} 
\cdot \left( 1 - \frac{\gamma}{2} + \frac{\gamma^{2}}{12} + \cdots \right) \right), 
$$ 
and hence 
$$
c_{0} \left( \pi_{!} \left( \Omega_{C/B}^{n} \right) \right) = (2n - 1)(g - 1), \ 
c_{1} \left( \pi_{\ast} \left( \Omega_{C/B}^{n} \right) \right) = 
\frac{6n^{2} - 6n + 1}{12} \pi_{\ast} \left( \gamma^{2} \right). 
$$ 
The first equality means the original Riemann-Roch's theorem given in 2.2. 
Furthermore, by putting $n = 1$ in the second equality 
$$
12 c_{1} \left( \pi_{\ast} \left( \Omega_{C/B}^{n} \right) \right) = 
(6n^{2} - 6n + 1) \cdot \pi_{\ast} \left( \gamma^{2} \right) = 
12 (6n^{2} - 6n + 1) \cdot c_{1} \left( \pi_{\ast} \left( \Omega_{C/B} \right) \right). 
$$ 
Since the Picard group of ${\cal M}_{g}$ is torsion-free, 
$$
c_{1} \left( \pi_{\ast} \left( \Omega_{C/B}^{n} \right) \right) = 
(6n^{2} - 6n + 1) \cdot c_{1} \left( \pi_{\ast} \left( \Omega_{C/B} \right) \right), 
$$ 
and hence we have {\bf Mumford's isomorphism} [Mu4]: 
$$
\det \left( \pi_{\ast} \left( \Omega_{C/B}^{n} \right) \right) \cong  
\det \left( \pi_{\ast} \left( \Omega_{C/B} \right) \right)^{\otimes (6n^{2} - 6n + 1)} 
$$
between line bundles over  $B$, 
where $\det(E)$ denotes the determinant line bundle associated with 
a vector bundle $E$. 
\vspace{2ex}

\noindent
\underline{\bf Remark.}
Morita [Mo] and Mumford [Mu5] conjectured that the stable cohomology groups 
defined for the moduli spaces of curves 
over ${\bf C}:$ 
\begin{eqnarray*}
H^{k}({\cal M}) 
& \stackrel{\rm def}{=} & 
H^{k}({\cal M}_{g}({\bf C}), {\bf Q}) = H^{k}(\Pi_{g}, {\bf Q}) \ \ 
(g \geq 3k -1) 
\\
& : & \mbox{independent of $g \geq 3k -1$ by Harer's result [H2]} 
\end{eqnarray*}
satisfies that 
\begin{eqnarray*}
\bigoplus_{k \geq 0} H^{k}({\cal M}) 
& = & 
{\bf Q} \left[ \kappa_{1}, \kappa_{2}, ...\right] 
: \ \mbox{freely generated over ${\bf Q}$} 
\\
& & \mbox{by the {\bf tautological classes} 
$\kappa_{i} = \pi_{\ast} \left( \left( 
c_{1} \left( \Omega_{{\cal C}/{\cal M}_{g}} \right) \right)^{i+1} \right)$.} 
\end{eqnarray*}
The free generatedness is proved by Miller [Mi] and Morita [Mo], 
and the whole conjecture is proved by Madsen and Weiss [MadW].

\newpage
\begin{center}
\underline{\large {\bf \S 6. Arithmetic theory of modular forms}} 
\end{center}

\noindent
\underline{\bf 6.1. Elliptic modular forms} 
\vspace{2ex}

The Eisenstein series of even degree $2k \geq 4$ 
(appeared in the Laurent coefficients of the $\wp$-function 
$\wp_{{\bf Z} + {\bf Z} \tau}(z)$): 

\begin{eqnarray*}
E_{2k}(\tau) \stackrel{\rm def}{=} 
\sum_{(m, n) \in {\bf Z}^{2} - \{(0,0)\}} \frac{1}{(m + n \tau)^{2k}} 
& \stackrel{\bf Ex. 6}{=} & 
2 \zeta(2k) + \frac{2 (2 \pi \sqrt{-1})^{2k}}{(2k-1)!} 
\sum_{n=1}^{\infty} \sigma_{2k-1}(n) \ q^{n} \\ 
& : & \mbox{the {\bf Fourier expansion}} \ 
\left( q = e^{2 \pi \sqrt{-1} \tau} \right) 
\end{eqnarray*}
is a holomorphic function of $\tau \in H_{1}$ which satisfies 
the following 2 conditions for $SL_{2}({\bf Z})$: 
\begin{itemize}

\item 
Automorphic condition of weight $2k$ : \\ 
${\displaystyle E_{2k} \left( \frac{a \tau + b}{c \tau + d} \right) 
= (c \tau + d)^{2k} E_{2k}(\tau) \ \mbox{for any} \ 
\left( \begin{array}{cc} a & b \\ c & d \end{array} \right) 
\in SL_{2}({\bf Z});}$ 

\item 
Cusp condition : \\ 
$E_{2k}(\tau)$ is holomorphic at $q = 0$ 
$\left( \Leftrightarrow \tau =  \ 
\mbox{the unique cusp $\sqrt{-1} \cdot \infty$ of 
$SL_{2}({\bf Z})$} \right).$

\end{itemize} 

{\bf (Elliptic) modular forms} are holomorphic functions on $H_{1}$ 
satisfying the automorphic and cusp conditions 
for a congruence subgroup of $SL_{2}({\bf Z}).$ 
\vspace{2ex}

\noindent
\underline{\bf Fourier expansion and number theory.} 
The theory of elliptic modular forms and their Fourier expansions 
has the following applications to number theory:  
\begin{enumerate}

\item
${\displaystyle \sigma_{7}(n) = \sigma_{3}(n) + 
120 \sum_{i=1}^{n-1} \sigma_{3}(i) \sigma_{3}(n-i) \ \ 
\left( \Leftarrow E_{8}(\tau) = \frac{3}{7} E_{4}(\tau)^{2} 
\ \mbox{in Exercise 1} \right).}$

\item 
Jacobi's theorem : ${\displaystyle 
\sharp \left\{ (a_{i})_{1 \leq i \leq 4} \in {\bf Z}^{4} \ \left| \ 
\sum_{i=1}^{4} a_{i}^{2} = n \right. \right\} = 
8 \sum_{d|n, 4 \nmid d} d}$ \\ 
$\left( \Leftarrow \ \mbox{the theta series} \ 
\left( \sum_{n \in {\bf Z}} q^{n^{2}} \right)^{4} \ 
\mbox{is expressed by Eisenstein series for} \ 
\Gamma(2) \right).$ 

\item 
Deligne-Serre's theorem [D, DS]:  
For a normalized Hecke eigenform $f = \sum_{n} a(n) q^{n} $ of weight $k$ 
and character $\varepsilon$ for $\Gamma_{0}(N),$ 
there is a 2-dimensional Galois representation $\rho_{f}$ such that 
${\rm tr}(\rho_{f}(F_{\overline{p}})) = a(p)$ and 
$\det(\rho_{f}(F_{\overline{p}})) = \varepsilon(p) p^{k-1}$ 
for any Frobenius automorphism $F_{\overline{p}}$ for unramified primes $p.$ 

\item 
Serre's example [Se]: 
Let $L$ be the decomposition field of $x^{3} - x - 1$ 
which is a Galois extension over ${\bf Q}$ and contains $K = {\bf Q}(\sqrt{-23})$ 
such that the Galois group ${\rm Gal}(L/{\bf Q})$ is isomorphic to 
the symmetric group $S_{3}$ of degree $3$. 
Put  
$$
f(\tau) = 
\frac{1}{2} \left( \sum_{m,n \in {\bf Z}} q^{m^{2} + mn + 6 n^{2}} - 
\sum_{m,n \in {\bf Z}} q^{2 m^{2} + mn + 3 n^{2}} \right) = 
\sum_{n=1}^{\infty} a(n) q^{n}.
$$
Then $f(\tau)$ is a normalized Hecke eigenform of weight $1,$ 
and hence by Deligne-Serre's theorem, for any prime $p \neq 23,$ 
${\rm tr}(\rho_{f}(F_{\overline{p}})) = a(p),$ 
${\displaystyle \det(\rho_{f}(F_{\overline{p}})) 
= \left( \frac{-23}{p} \right) = \left( \frac{p}{23} \right)}$ 
and $\sharp \langle \rho(F_{\overline{p}}) \rangle$ 
is equal to the residue index $f_{L/{\bf Q}}(p)$ of $p$ in $L/{\bf Q}$. 
Therefore, this gives an example of nonabelian class field theory. 

\end{enumerate}

\noindent
\underline{\bf Exercise 10.} 
Show the above 1. 
\vspace{2ex}

\noindent
\underline{\bf Exercise 11.} 
Let the notation be as in the above 4. Serre's example. 
Then prove that for $p \neq 23,$ 
one of the following cases necessarily happens: 
\begin{eqnarray*}
a(p) = 2, \ \left( \frac{p}{23} \right) = 1 & \Longleftrightarrow & 
f_{L/{\bf Q}}(p) = 1, 
\\
a(p) = 0, \ \left( \frac{p}{23} \right) = -1 & \Longleftrightarrow & 
f_{K/{\bf Q}}(p) = 2, \ f_{L/{\bf Q}}(p) = 2, 
\\
a(p) = -1, \ \left( \frac{p}{23} \right) = 1 & \Longleftrightarrow & 
f_{K/{\bf Q}}(p) = 1, \ f_{L/{\bf Q}}(p) = 3, 
\end{eqnarray*}
and describe the decomposition of primes $2, 3, 5, 59$ in $K$ and $L$ 
respectively. 
\vspace{2ex}

\noindent
\underline{\bf Rationality of modular forms.}
For $\tau \in H_{1},$ 
\begin{eqnarray*}
E_{\tau} & \stackrel{\rm def}{=} & 
{\bf C}/({\bf Z} + {\bf Z} \tau) \ 
\mbox{define a family of elliptic curves over $H_{1},$} 
\\ 
z_{\tau} & \stackrel{\rm def}{=} & 
\mbox{the natural coordinate of ${\bf C}$} 
\\ 
& \Rightarrow & 
d z_{\tau} : \mbox{a canonical base of $H^{0}(E_{\tau}, \Omega_{E_{\tau}}),$} 
\end{eqnarray*}
and 
\begin{eqnarray*}
& & 
\left( \begin{array}{cc} a & b \\ c & d \end{array} \right) 
\in SL_{2}({\bf Z}) 
\\ 
& \Rightarrow & 
E_{\frac{a \tau + b}{c \tau + d}} 
\stackrel{\times (c \tau + d)}{\longrightarrow} 
{\bf C} / ({\bf Z}(c \tau + d) + {\bf Z}(a \tau + b)) = 
{\bf C} / ({\bf Z} + {\bf Z} \tau) = E_{\tau} 
\\ 
& \Rightarrow & 
d z_{\frac{a \tau + b}{c \tau + d}} = \frac{1}{c \tau + d} dz_{\tau}. 
\end{eqnarray*}
If $f(\tau)$ is a modular form of weight $k$ for $SL_{2}({\bf Z}),$ 
then 
$$
f \left( \frac{a \tau + b}{c \tau + d} \right) 
\left( d_{\frac{a \tau + b}{c \tau + d}} \right)^{\otimes k} = 
(c \tau + d)^{k} f(\tau) \left( \frac{1}{c \tau + d} \right)^{k} 
(d z_{\tau})^{\otimes k} 
= f(\tau) (dz_{\tau})^{\otimes k}, 
$$
and hence $f(\tau) d z_{\tau}$ $(\tau \in H_{1})$ is 
$SL_{2}({\bf Z})$-invariant, i.e., defines a holomorphic section of 
the line bundle on $[H_{1} / SL_{2}({\bf Z})]$ 
whose fiber over $\tau \in H_{1}$ 
is given by $H^{0}(E_{\tau}, \Omega_{E_{\tau}})^{\otimes k}.$ 

Let ${\cal M}_{1}$ be the moduli stack of elliptic curves, 
$\pi : {\cal E} \rightarrow {\cal M}_{1}$ be the universal elliptic curve, 
and $\pi_{\ast}(\Omega_{{\cal E}/{\cal M}_{1}})$ denote 
a line bundle on ${\cal M}_{1}$ defined by the direct image of 
the sheaf $\Omega_{{\cal E}/{\cal M}_{1}}$ of relative $1$-forms 
on ${\cal E}/{\cal M}_{1},$ i.e., 
$$
\pi_{\ast}(\Omega_{{\cal E}/{\cal M}_{1}})(S) \stackrel{\rm def}{=} 
H^{0}(E, \Omega_{E/S}), 
$$
for elliptic curves $E$ over schemes $S.$ 
Then an {\bf integral modular form} $f$ of weight $k$ is defined as 
an element of 
$$
H^{0} \left( {\cal M}_{1}, 
\pi_{\ast} \left( \Omega_{\cal E}/{\cal M}_{1} \right)^{\otimes k} \right), 
$$
i.e., a global section of 
$\pi_{\ast}(\Omega_{{\cal E}/{\cal M}_{1}})^{\otimes k}$ 
on ${\cal M}_{1}$ which is, by the above dictionary on the moduli stack, 
a system of 
$$
\left\{ \left. 
\mbox{sections $f_{S}$ of $H^{0} \left( E, \Omega_{E/S} \right)^{\otimes k}$} 
\ \right| \ \mbox{$E:$ elliptic curves over $S$} \right\} 
$$
which are functorial for schemes $S.$ 
Hence 
\begin{eqnarray*}
& & E/S \ : \ 
\mbox{the Tate curve $y^{2} + xy = x^{3} + a_{4}(q) x + a_{6}(q)$ 
over ${\bf Z}((q))$} 
\\ 
& \Rightarrow & 
\frac{du}{u} = \frac{d X(u, q)}{X(u, q) + 2 Y(u, q)} = 
\frac{dx}{x + 2y} : \ \mbox{a base of $1$-forms on the Tate curve} 
\\ 
& \Rightarrow & 
f \ \mbox{is represented as} \ 
F(f) \left( \frac{dx}{x + 2y} \right)^{\otimes k}, 
\end{eqnarray*}
where $F(f) \in {\bf Z}((q))$ is called the {\bf evaluation} of $f$ 
on the Tate curve under the trivialization of 
$\pi_{\ast} \left( \Omega_{{\cal E}/{\cal M}_{1}} \right)$ on ${\bf Z}((q)).$ 
By Theorem 4.1 (2), 
\begin{eqnarray*}
q = e^{2 \pi \sqrt{-1} \tau} 
& \Rightarrow & 
{\bf C}/({\bf Z} + {\bf Z} \tau) \cong {\bf C}^{\times} / \langle q \rangle
\\ 
& \Rightarrow & 
\frac{dx}{x + 2y} = 2 \pi \sqrt{-1} \ 
\frac{d \wp_{{\bf Z} + {\bf Z} \tau}(z_{\tau})}
{\wp'_{{\bf Z} + {\bf Z} \tau}(z_{\tau})} = 2 \pi \sqrt{-1} d z_{\tau} 
\\ 
& \Rightarrow & 
f(\tau) = (2 \pi \sqrt{-1})^{k} F(f) (d z_{\tau})^{\otimes k}. 
\end{eqnarray*}
Therefore, ignoring the factor $(2 \pi \sqrt{-1})^{k},$ 
\begin{center}
{\bf the evaluation on the Tate curve = 
the classical Fourier expansion,}
\end{center}  
and hence  
\begin{center}
{\bf a modular form is integral {\boldmath $\Longleftrightarrow$} 
its Fourier coefficients are integral.}
\end{center}  
\vspace{1ex}

\noindent
\underline{\bf Exercise 12.} 
\begin{itemize}

\item 
Prove that 
${\displaystyle \frac{E_{4}(\tau)}{2 \zeta(4)}},$ 
${\displaystyle \frac{E_{6}(\tau)}{2 \zeta(6)}}$ and 
${\displaystyle \Delta(\tau) \stackrel{\rm def}{=} \frac{1}{1728} 
\left( \left( \frac{E_{4}(\tau)}{2 \zeta(4)} \right)^{3} - 
\left( \frac{E_{6}(\tau)}{2 \zeta(6)} \right)^{2} \right)}$ 
are integral (elliptic) modular forms for $SL_{2}({\bf Z}).$ 

\item 
Using that $\Delta(\tau) \neq 0$ $(\tau \in H)$ and 
that modular forms for $SL_{2}({\bf Z})$ of weight $0$ are constant, 
prove that all integral modular forms for $SL_{2}({\bf Z})$ 
are generated over ${\bf Z}$ by these $3$ modular forms.   
\end{itemize}
\vspace{1ex}

\noindent
\underline{\bf 6.2. Siegel modular forms (SMFs)} 
\vspace{2ex}

\noindent
\underline{\bf Moduli of abelian varieties.}
Let $g$ be a positive integer $>1.$ 
Then in a similar way to constructing moduli of curves given in 5.1, 
it is shown in [FKM] that there exists the fine moduli space ${\cal A}_{g}$ 
as an algebraic stack over ${\bf Z}$ 
classifying principally polarized abelian varieties of dimension $g.$ 
By the correspondence: 
\begin{eqnarray*}
& & Z = \left( \begin{array}{c} 
\mbox{\boldmath $z$}_{1} \\ \vdots \\ \mbox{\boldmath $z$}_{g} \end{array} 
\right) \in H_{g} : \ \mbox{the Siegel upper half space of degree $g$}
\\
& \leftrightarrow & 
\left( {\bf C}^{g} / ({\bf Z}^{g} + {\bf Z}^{g} \cdot Z) ; 
\iota(\mbox{\boldmath $e$}_{i}) = \left\{ \begin{array}{ll}
\mbox{\boldmath $e$}_{i} & (1 \leq i \leq g), 
\\
\mbox{\boldmath $z$}_{i-g} & (g+1 \leq i \leq 2g) \end{array} \right. 
\right), 
\end{eqnarray*}
$H_{g}$ becomes the fine moduli space of principally polarized 
abelian varieties $X$ of dimension $g$ over ${\bf C}$ 
with symplectic isomorphism 
${\bf Z}^{2g} \stackrel{\sim}{\rightarrow} H_{1}(X, {\bf Z}).$ 
Hence the orbifold ${\cal A}_{g}({\bf C})$ is given by 
the quotient stack of $H_{g}$ by the integral symplectic group 
$Sp_{2g}({\bf Z})$ of degree $g:$ 
$$
{\cal A}_{g}({\bf C}) \ = \ \left[ H_{g} / Sp_{2g}({\bf Z}) \right]. 
$$
%\vspace{1ex}

\noindent
\underline{\bf Definition of SMFs.} 
Let $\lambda$ be the {\bf Hodge line bundle} on ${\cal A}_{g}$ 
which is defined by 
\begin{eqnarray*}
& & 
\lambda \stackrel{\rm def}{=} \bigwedge^{g} 
\rho_{\ast} \left( \Omega_{{\cal X}/{\cal A}_{g}} \right) \ \ 
\mbox{$(\rho : {\cal X} \rightarrow {\cal A}_{g}$ 
denotes the universal abelian scheme}) 
\\
& \Rightarrow & 
\lambda(S) = \bigwedge^{g} 
H^{0} \left( X, \Omega_{X/S} \right) \ \ 
\mbox{for abelian schemes $X/S$ of relative dimension $g.$} 
\end{eqnarray*}
Then for $h \in {\bf Z}$ and a ${\bf Z}$-module $M,$ we call elements of 
$$
S_{g,h}(M) \stackrel{\rm def}{=} 
H^{0} \left( {\cal A}_{g}, \lambda^{\otimes h} \otimes_{\bf Z} M \right)
$$
{\bf Siegel modular forms} of degree $g$ and weight $h$ 
with coefficients in $M.$ 

For the natural coordinate $z_{1},..., z_{g}$ on the complex abelian varieties 
$$
X_{Z} = {\bf C}^{g}/({\bf Z}^{g} + {\bf Z}^{g} \cdot Z),
$$
$dz_{1},..., dz_{g}$ give a base of $H^{0}(X_{Z}, \Omega_{X_{Z}}),$ 
and hence as in the elliptic case, 
\begin{eqnarray*}
& & 
\varphi = (2 \pi \sqrt{-1})^{gh} \cdot f \cdot 
\left( dz_{1} \wedge \cdots \wedge dz_{g} \right)^{\otimes h} 
\in S_{g,h}({\bf C}) 
= H^{0} \left( \left[ H_{g}/Sp_{2g}({\bf Z}) \right], 
\lambda^{\otimes h} \right) 
\\
& \stackrel{(\ast)}{\Rightarrow} & 
\left\{ \begin{array}{l} 
\mbox{$f = f(Z)$ is a holomorphic function of $Z \in H_{g}$ such that}
\\
\mbox{${\displaystyle f(G(Z)) = \det(CZ + D)^{h} \cdot f(Z)}$ 
for any ${\displaystyle G = 
\left( \begin{array}{cc} A & B \\ C & D \end{array} \right) 
\in Sp_{2g}({\bf Z})}$}
\end{array} \right.
\end{eqnarray*}
which is known as the usual definition of {\bf analytic Siegel modular forms}. 
In particular, $f(Z)$ is invariant under the transformation
$$
Z \longmapsto Z + B 
$$
by integral symmetric matrices $B$ of degree $g,$ 
and hence it can be expanded as a power series of 
$\exp \left( 2 \pi \sqrt{-1} z_{ij} \right)$ 
$\left( Z = (z_{ij})_{i,j} \in H_{g} \right)$ 
which is called the {\bf (classical) Fourier expansion} of $f.$ 
\vspace{2ex}

\noindent
\underline{\bf Exercise 13.} 
Prove the above $\stackrel{(\ast)}{\Rightarrow}.$
\vspace{2ex}

It is shown by Satake that 
${\cal A}_{g/{\bf C}} = {\cal A}_{g} \otimes_{\bf Z} {\bf C}$ has 
the {\bf Satake compactification}: 
$$
{\cal A}_{g/{\bf C}}^{\ast} = \coprod_{i=0}^{g} {\cal A}_{i/{\bf C}}, 
$$
obtained as the Zariski closure of a projective embedding 
using Siegel modular forms of sufficiently large weight. 
Then the codimension of 
${\cal A}_{g/{\bf C}}^{\ast} - {\cal A}_{g/{\bf C}}$ in 
${\cal A}_{g/{\bf C}}^{\ast}$ is 
$$
\frac{g(g + 1)}{2} - \frac{(g - 1)g}{2} = g > 1, 
$$
and hence ignoring $(2 \pi \sqrt{-1})^{gh} 
\left( dz_{1} \wedge \cdots \wedge dz_{g} \right)^{\otimes h},$ 
\begin{eqnarray*}
& & 
\mbox{$\varphi$ is an analytic Siegel modular form}
\\
& \Rightarrow & 
\mbox{$\varphi$ is an analytic section on ${\cal A}_{g/{\bf C}}^{\ast}$ 
\ \ (by Hartogs' theorem)} 
\\
& \Rightarrow & 
\mbox{$\varphi$ is an algebraic section on 
${\cal A}_{g/{\bf C}}^{\ast}$ 
\ \ (by GAGA's principle of Serre)} 
\\
& \Rightarrow & 
\mbox{$\varphi$ is an algebraic section on ${\cal A}_{g/{\bf C}}$} 
\\
& \Rightarrow & 
\varphi \in S_{g,h}({\bf C}). 
\end{eqnarray*}
Therefore, the above $\stackrel{(\ast)}{\Rightarrow}$ is in fact 
an equivalence $\stackrel{(\ast)}{\Longleftrightarrow},$ 
and $S_{g, h}({\bf C})$ is finite dimensional over ${\bf C}$ 
by the compactness of ${\cal A}_{g/{\bf C}}^{\ast}.$ 
\vspace{2ex}

\noindent
\underline{\bf Fourier expansion of SMFs.} 
By Mumford's theory [Mu3] on degenerating abelian varieties, 
there exists a semiabelian scheme expressed as 
$$
\left. {\bf G}_{m}^{g} \right/ 
\langle (q_{ij})_{1 \leq i \leq g} \ | \ 1 \leq j \leq g \rangle 
$$
over the ring 
$$
{\bf Z} \left[ q_{ij}^{\pm 1} \ (i \neq j) \right] 
[[ q_{11},..., q_{gg} ]], 
$$
where $q_{ij}$ $(1 \leq i, j \leq g)$ are variables 
with symmetry $q_{ij} = q_{ji}.$ 
This semiabelian scheme gives a family of complex abelian varieties 
$$
{\bf C}^{g} / \left( {\bf Z} + {\bf Z} \cdot Z \right) \ \cong \ 
\left( {\bf C}^{\times} \right)^{g} \left/ 
\left\langle \left( \exp(2 \pi \sqrt{-1} z_{ij}) \right)_{1 \leq i \leq g} 
\ | \ 1 \leq j \leq g \right\rangle \right.
$$
when $q_{ij} = \exp(2 \pi \sqrt{-1} z_{ij})$ 
for $Z = (z_{ij})_{i,j} \in H_{g}.$ 
Then the natural coordinates $u_{1},..., u_{g}$ on ${\bf G}_{m}^{g}$ 
give a base $du_{1}/u_{1},..., du_{g}/u_{g}$ of $1$-forms 
on this semiabelian scheme, 
and hence the evaluation of any $\varphi \in S_{g,h}(M)$ gives 
\begin{eqnarray*}
\varphi 
& = & 
F(\varphi) \cdot \left( 
\left( du_{1}/u_{1} \right) \wedge \cdots \wedge \left( du_{g}/u_{g} \right) 
\right)^{\otimes h} 
\\
& = & 
(2 \pi \sqrt{-1})^{gh} \cdot F(\varphi) \cdot 
(dz_{1} \wedge \cdots \wedge dz_{g})^{\otimes h} \ \ 
\mbox{(if $M = {\bf C}$ and $u_{i} = \exp(2 \pi \sqrt{-1} z_{i})$).} 
\end{eqnarray*}
Therefore, we have a linear map: 
$$
F : S_{g,h}(M) \longrightarrow 
{\bf Z} \left[ q_{ij}^{\pm 1} \ (i \neq j) \right] 
[[ q_{11},...,q_{gg} ]] \otimes_{\bf Z} M,
$$
which we call the {\bf arithmetic Fourier expansion.}
\vspace{2ex}

\noindent
\underline{\bf Theorem 6.1.} (Chai and Faltings [FaC]) 
\begin{it}

{\rm (1)} (Arithmetic Fourier expansion) 
$F$ is functorial for $M,$ and if $M = {\bf C},$ then $F(\varphi)$ is 
the classical Fourier expansion by $q_{ij} = 
\exp \left( 2 \pi \sqrt{-1} z_{ij} \right)$ for $(z_{ij})_{i,j} \in H_{g}.$ 
Furthermore, $F$ is injective, 
and for a submodule $N$ of $M$ and $\varphi \in S_{g,h}(M),$ 
$$
\varphi \in S_{g,h}(N) \ \Longleftrightarrow \ 
F(\varphi) \in {\bf Z} 
\left[ q_{ij}^{\pm 1} \right] [[ q_{ii} ]] \otimes_{\bf Z} M. 
$$

{\rm (2)} (Finiteness) 
$S_{g,h}({\bf Z})$ is a free ${\bf Z}$-module of finite rank 
such that $S_{g,h}({\bf Z}) \otimes_{\bf Z} {\bf C} = S_{g,h}({\bf C})$ 
and that $S_{g,0}({\bf Z}) = {\bf Z},$ $S_{g,h}({\bf Z}) = \{ 0 \}$ 
if $n < 0.$ 
Furthermore, the ring of integral Siegel modular forms of 
degree $g$ over ${\bf Z}:$ 
$$
S_{g}^{\ast}({\bf Z}) \stackrel{\rm def}{=} 
\bigoplus_{h \geq 0} S_{g,h}({\bf Z}) 
$$
is a normal ring finitely generated over ${\bf Z}.$ 
\end{it} 
\vspace{2ex}

\underline{\it Sketch of Proof.} 
(1) The functoriality for $M$ and the compatibility 
with the classical Fourier expansion is clear from the above. 
Since ${\cal A}_{g}$ is smooth over ${\bf Z},$ 
we have the following left exactness of $S_{g,h}(M)$ for $M:$ 
\begin{eqnarray*}
& & 
0 \rightarrow N \rightarrow M \rightarrow (M/N) \rightarrow 0 
\\ 
& \Rightarrow & 
0 \rightarrow \lambda^{\otimes h} \otimes_{\bf Z} N \rightarrow 
\lambda^{\otimes h} \otimes_{\bf Z} M \rightarrow 
\lambda^{\otimes h} \otimes_{\bf Z} (M/N) \rightarrow 0 
\\ 
& \Rightarrow & 
0 \rightarrow S_{g,h}(N) \rightarrow S_{g,h}(M) \rightarrow S_{g,h}(M/N). 
\end{eqnarray*}
We prove the injectivity of $F.$ 
Since any ${\bf Z}$-module $M$ is the direct limit of 
finitely generated ${\bf Z}$-modules, and cohomology and tensor product 
commute with direct limit, 
we may assume that $M$ is a finitely generated ${\bf Z}$-module, 
hence by the left exactness for $M,$ 
we may put $M = {\bf Z}$ or $= {\bf Z} / p {\bf Z}$ $(p:$ a prime number). 
Therefore, the injectivity follows from that ${\cal A}_{g} \otimes M$ 
is smooth over the ring $M$ with geometrically irreducible fibers 
which is proved in [FaC]. 
Hence the remains of (1) follows from this injectivity and 
the left exactness of $S_{g,h}.$ 

(2) is derived by the following result in [FaC]: 
there exists an algebraic stack $\overline{\cal A}_{g}$ 
which is proper smooth over ${\bf Z}$ 
and contains ${\cal A}_{g}$ as its open dense substack, 
and any integral Siegel modular form of weight $k$ can be extended to 
a section on $\overline{\cal A}_{g}$ of an extension 
$\overline{\lambda}^{\otimes k}$ of $\lambda^{\otimes k}$ 
(called {\bf Koecher's principle}). 

The finiteness of ${\rm rank}_{\bf Z} S_{g,h}({\bf Z})$ follows from 
these results immediately. 
Further, there is $m \in {\bf N}$ such that $\overline{\lambda}^{\otimes m}$ 
defines a projective morphism 
$\overline{\cal A}_{g} \rightarrow {\mathbb P}^{n}_{\bf Z}$ 
which can be, by the theory of Stein factorization, decomposed as 
$\overline{\cal A}_{g} \rightarrow {\cal A}_{g}^{\ast} \rightarrow 
{\mathbb P}^{n}_{\bf Z}$ such that 
$\overline{\cal A}_{g} \rightarrow {\cal A}_{g}^{\ast}$ has 
connected geometric fibers and 
${\cal A}_{g}^{\ast} \rightarrow {\mathbb P}^{n}_{\bf Z}$ is finite. 
Therefore, replacing $m$ by a multiple $\overline{\lambda}^{\otimes m}$ 
defines a immersion of ${\cal A}_{g}^{\ast},$ and hence 
$\bigoplus_{k \geq 0} 
H^{0}(\overline{\cal A}_{g}, \overline{\lambda}^{\otimes mk})$ and 
$S_{g}^{\ast}({\bf Z})$ are normal rings finitely generated over ${\bf Z}.$ 
QED. 
\vspace{2ex}

\noindent
\underline{\bf Ring of SMFs of degree 2 and 3.} 
(Igusa [Ig1,3], Tsuyumine [Ty1])
For $g > 1$ and $h > g + 1,$ 
the {\bf Eisenstein series} of degree $g > 1$ and weight $h$ 
is a function of $Z \in H_{g}$ defined by 
$$
E_{g,h}(Z) \stackrel{\rm def}{=} 
\sum_{G \in \Gamma_{\infty} \backslash Sp_{2g}({\bf Z})} 
\det \left( C Z + D \right)^{-h} ; \ \ 
G = \left( \begin{array}{cc} A & B \\ C & D \end{array} \right), 
$$
where 
$$
\Gamma_{\infty} \stackrel{\rm def}{=} 
\left\{ \left( \begin{array}{cc} U & B \\ 0 & \mbox{}^{t}U^{-1} \end{array} 
\right) \in Sp_{2g}({\bf Z}) \right\}. 
$$
Then $E_{g,h}$ becomes a Siegel modular form 
with Fourier coefficients in ${\bf Q},$ 
and hence an element of $S_{g,h}({\bf Q}).$ 
Igusa [Ig1] proved that 
$$
S_{2}^{\ast}({\bf C}) \ = \ 
{\bf C} \left[ E_{4}, E_{6}, \Delta_{10}, \Delta_{12} \right] \bigoplus 
\Delta_{35} \cdot 
{\bf C} \left[ E_{4}, E_{6}, \Delta_{10}, \Delta_{12} \right], 
$$
where $E_{h} = E_{2, h},$ 
$\Delta_{10} = E_{4} E_{6} - E_{10},$ 
$\Delta_{12} = 441 E_{4}^{3} + 250 E_{6}^{2} - 691 E_{12}$ and 
$\Delta_{35} \in S_{2, 35}({\bf C})$ is given by Ibukiyama as 
$$
\Delta_{35} \left( \begin{array}{cc} z_{11} & z_{12} \\ z_{12} & z_{22} 
\end{array} \right) 
\ = \ \left| \begin{array}{cccc} 
4 E_{4} & 6 E_{6} & 10 \Delta_{10} & 12 \Delta_{12} 
\\ 
& & & 
\\
{\displaystyle \frac{\partial E_{4}}{\partial z_{11}}} & 
{\displaystyle \frac{\partial E_{6}}{\partial z_{11}}} & 
{\displaystyle \frac{\partial \Delta_{10}}{\partial z_{11}}} & 
{\displaystyle \frac{\partial \Delta_{12}}{\partial z_{11}}} 
\\ 
& & & 
\\
{\displaystyle \frac{\partial E_{4}}{\partial z_{12}}} & 
{\displaystyle \frac{\partial E_{6}}{\partial z_{12}}} & 
{\displaystyle \frac{\partial \Delta_{10}}{\partial z_{12}}} & 
{\displaystyle \frac{\partial \Delta_{12}}{\partial z_{12}}} 
\\ 
& & & 
\\
{\displaystyle \frac{\partial E_{4}}{\partial z_{22}}} & 
{\displaystyle \frac{\partial E_{6}}{\partial z_{22}}} & 
{\displaystyle \frac{\partial \Delta_{10}}{\partial z_{22}}} & 
{\displaystyle \frac{\partial \Delta_{12}}{\partial z_{22}}} 
\end{array} \right|. 
$$
In [I5], this result was extended to $S_{2}^{*}(R)$, 
where $R$ is a ${\bf Z}$-algebra in which $6$ is invertible. 

Igusa [Ig3] determined generators of $S_{2}^{*}({\bf Z})$, 
and Tsuyumine [Ty1] gave explicit generators of $S_{3}^{*}({\bf C})$. 
\vspace{2ex}

\noindent
\underline{\bf 6.3. Teichm\"{u}ller modular forms (TMFs)} 
\vspace{2ex}

\begin{tabular}{ccl}
Analytic 
& : & automorphic functions on the Teichm\"{u}ller space 
\\
& = & automorphic forms on the moduli space of Riemann surfaces, 
\\ 
Algebraic 
& : & global sections of line bundles on the moduli of curves. 
\end{tabular}
\vspace{2ex}
\\
This naming is an analogy of 
\vspace{1ex}

\begin{tabular}{ccl}
&   & Siegel modular forms (SMFs) \\ 
& = & automorphic functions on the Siegel upper half space 
\\
& = & global sections of line bundles \\ 
&   & on the moduli of principally polarized abelian varieties. 
\end{tabular}
\vspace{2ex}

\noindent
\underline{\bf Definition of TMFs.}
Let $\pi : {\cal C} \rightarrow {\cal M}_{g}$ be the universal curve 
over the moduli stack of proper smooth curves of genus $g > 1,$ 
and let $\lambda \stackrel{\rm def}{=} 
\bigwedge^{g} \pi_{*} \left( \Omega_{{\cal C}/{\cal M}_{g}} \right)$ 
be the {\bf Hodge line bundle}. 
Then for a ${\bf Z}$-module $M$, we call elements of 
$$
T_{g,h}(M) \stackrel{\rm def}{=} 
H^{0}({\cal M}_{g}, \lambda^{\otimes h} \otimes_{\bf Z} M) 
$$
{\bf Teichm\"{u}ller modular forms} of degree $g$ and weight $h$ 
with coefficients in $M$. 
By the pullback of the Torelli map 
$\tau : {\cal M}_{g} \rightarrow {\cal A}_{g}$ 
sending curves to their Jacobian varieties with canonical polarization, 
we have a linear map 
$$
\tau^{\ast} : S_{g,h}(M) \longrightarrow T_{g,h}(M)
$$
for ${\bf Z}$-modules $M.$ 
If $g = 2, 3,$ then the image of the Torelli map is Zariski dense, 
and hence $\tau^{\ast}$ is injective. 

If $n \geq 3,$ then 
\begin{eqnarray*}
{\cal M}_{g,n/{\bf C}} & \stackrel{\rm def}{=} & 
\mbox{the moduli space of proper smooth curves over ${\bf C}$} 
\\
& & \mbox{of genus $g$ with symplectic level $n$ structure,}
\\
{\cal A}_{g,n/{\bf C}} & \stackrel{\rm def}{=} & 
\mbox{the moduli space of principally polarized abelian varieties 
over ${\bf C}$} 
\\
& & \mbox{of dimension $g$ with symplectic level $n$ structure}
\end{eqnarray*}
are given as fine moduli schemes over ${\bf C}.$ 
Let ${\cal M}_{g,n/{\bf C}}^{\ast}$ be the {\bf Satake-type} compactification, 
i.e., normalization of the Zariski closure of 
$$
(\iota \circ \tau) ({\cal M}_{g,n/{\bf C}}) \subset 
{\cal A}_{g,n/{\bf C}}^{\ast}, 
$$ 
where $\tau : {\cal M}_{g,n/{\bf C}} \rightarrow {\cal A}_{g,n/{\bf C}}$ 
denote the Torelli map, 
and $\iota : {\cal A}_{g,n/{\bf C}} \rightarrow {\cal A}_{g,n/{\bf C}}^{\ast}$ 
denote the natural inclusion to the Satake compactification. 
Then each point of ${\cal M}_{g,n/{\bf C}}^{\ast} - {\cal M}_{g,n/{\bf C}}$ 
corresponds to the product $J_{1} \times \cdots \times J_{m}$ 
of Jacobian varieties over ${\bf C}$ with canonical polarization and 
symplectic level $n$ structure such that $\sum_{i=1}^{m} \dim(J_{i}) \leq g$ 
and that $(m, g) \neq \left( 1, \dim(J_{1}) \right).$ 
Therefore, if $g \geq 3,$ 
then ${\cal M}_{g,n/{\bf C}}^{\ast} - {\cal M}_{g,n/{\bf C}}$ 
has codimension $2$ in ${\cal M}_{g,n/{\bf C}}^{\ast},$ 
and hence by applying Hartogs' theorem to 
${\cal M}_{g,n/{\bf C}} \subset {\cal M}_{g,n/{\bf C}}^{\ast}$ 
and GAGA's principle to ${\cal M}_{g,n/{\bf C}}^{\ast},$ 
one can see that analytic TMFs become algebraic TMFs, i.e., 
$$
T_{g,h}({\bf C}) \cong 
\left\{ \begin{array}{l} 
\mbox{holomorphic functions on the Teichm\"{u}ller space $T_{g}$} \\ 
\mbox{with automorphic condition of weight $h$} \\ 
\mbox{for the action of the Teichm\"{u}ller modular group $\Pi_{g}$}
\end{array} \right\}, 
$$
and this space is finite dimensional over ${\bf C}.$ 
\vspace{2ex}

\noindent
\underline{\bf Exercise 14.} 
Give a precise definition of analytic Teichm\"{u}ller modular forms. 
\vspace{2ex}

\noindent
\underline{\bf Expansion of TMFs.} 
Let $C_{\Delta}$ be the generalized Tate curve given in Theorem 4.2 
which is smooth over the ring $B_{\Delta}.$ 
Then as in the elliptic and Siegel modular case, 
the evaluation on $C_{\Delta}$ 
($=$ the expansion by the corresponding local coordinates on ${\cal M}_{g}$) 
gives rise to a homomorphism 
$$
\kappa_{\Delta} : T_{g,h}(M) \longrightarrow B_{\Delta} \otimes_{\bf Z} M. 
$$
\noindent
\underline{\bf Theorem 6.2.} ([I3]). 
\begin{it} Fix $g >1$ and $h \in {\bf Z}$. 

{\rm (1)} $\kappa_{\Delta}$ is injective, 
and for a Teichm\"{u}ller modular form $f \in T_{g,h}(M)$ 
and a submodule $N$ of $M$, 
$$
f \in T_{g,h}(N) \ \Longleftrightarrow \ 
\kappa_{\Delta}(f) \in B_{\Delta} \otimes_{\bf Z} N. 
$$

{\rm (2)} For a Siegel modular form $\varphi \in S_{g,h}(M),$ \end{it}
$$
\kappa_{\Delta}(\tau^{*}(\varphi)) \ = \ F(\varphi)|_{q_{ij} = p_{ij}}, 
$$ 
where $p_{ij}$ are the multiplicative periods of $C_{\Delta}$ 
given in Theorem 4.2 (4). 
\vspace{2ex}

\underline{\it Proof.} 
(1) follows from the fact that $C_{\Delta}$ corresponds to the generic point 
on ${\cal M}_{g},$ and the argument in the proof of Theorem 6.1 (1) replacing 
${\cal A}_{g}$ by ${\cal M}_{g}$ which is proper and smooth over ${\bf Z}$ 
with geometrically irreducible fibers (see 5.2). 
(2) follows from Theorem 4.2 (4). QED. 
\vspace{2ex}

\noindent
\underline{\bf Schottky problem.} 
As an application of Theorem 6.2, 
we can give a solution to the Schottky problem, i.e. 
characterizing Siegel modular forms vanishing on the Jacobian locus, 
is given as follows: 
$$
\tau^{*}(\varphi) = 0 
\ \Longleftrightarrow \ 
F(\varphi)|_{q_{ij} = p_{ij}} = 0. 
$$
\begin{center}
\fbox{\Large{\sf {\boldmath $p_{ij}$} are computable, 
hence {\boldmath $\kappa_{\Delta}$} are computable}}
\end{center}
\vspace{1ex}

\noindent
Using the universal periods $p_{ij}$ given in Example 4.1, 
the above implies the following result of Brinkmann and Gerritzen [BG, G]: 
For the Fourier expansion 
$$
F(\varphi) = \sum_{T = (t_{ij})} a_{T} 
\prod_{1 \leq i < j \leq g} {q_{ij}}^{2 t_{ij}} 
\prod_{1 \leq i \leq g} {q_{ii}}^{t_{ii}} 
$$ 
of a Siegel modular form $\varphi$ vanishing on the Jacobian locus, 
\begin{eqnarray*}
& & \mbox{integers $s_{1},...,s_{g} \geq 0$ satisfy 
${\displaystyle \sum_{i = 1}^{g} s_{i} 
= \min \{ {\rm T}(T) \ | \ a_{T} \neq 0 \}}$} 
\\ 
& \Rightarrow &
\sum_{t_{ii} = s_{i}} a_{T} 
\prod_{i < j} \left( 
\frac{(x_{i} - x_{j})(x_{-i} - x_{-j})}{(x_{i} - x_{-j})(x_{-i} - x_{j})} 
\right)^{2 t_{ij}} = 0 \ \ \mbox{in $A_{0}$ (: given in Example 4.1).} 
\end{eqnarray*}

\noindent
\underline{\bf Schottky's {\boldmath$J$}.}
For $n \equiv 0 \ {\rm mod}(4),$ put 
\begin{eqnarray*}
L_{2n} & \stackrel{\rm def}{=} & 
\left\{ (x_{1},...,x_{2n}) \in {\bf R}^{2n} \left| \ 
2x_{i}, \ x_{i}-x_{j}, \ \frac{1}{2} \sum_{i} x_{i} \in {\bf Z} \right. 
\right\} 
\\ 
& : & \mbox{a lattice in ${\bf R}^{2n}$ with standard inner product,} 
\\
\varphi_{n}(Z) & \stackrel{\rm def}{=} & 
\sum_{(\lambda_{1},...,\lambda_{4}) \in L_{2n}^{4}} 
\exp \left( \pi \sqrt{-1} \sum_{i,j = 1}^{4} 
\langle \lambda_{i}, \lambda_{j} \rangle z_{ij} \right) \ 
\left( Z = (z_{ij})_{i,j} \in H_{4} \right) 
\\
& : & \mbox{a Siegel modular form of degree $4$ and weight $n,$} 
\\
J(Z) & \stackrel{\rm def}{=} & 
\frac{2^{2}}{3^{2} \cdot 5 \cdot 7} (\varphi_{4}(Z)^{2} - \varphi_{8}(Z)) 
: \ \mbox{\bf Schottky's {\boldmath $J$}} 
\\
& : & \mbox{an integral Siegel modular form of degree $4$ and weight $8.$} 
\end{eqnarray*}
Then Schottky and Igusa proved that the Zariski closure of 
the Jacobian locus in ${\cal A}_{4} \otimes_{\bf Z} {\bf C}$ 
is defined by $J = 0.$ 

Brinkmann and Gerritzen [BG, G] checked the above Brinkmann and Gerritzen's 
criterion for Schottky's $J,$ i.e., 
computed the lowest term of $J$ and showed that this is given by 
up to a constant 
$$
F \frac{q_{11} q_{22} q_{33} q_{44}}{\prod_{1 \leq i < j \leq 4} q_{ij}}, 
$$
where $F$ is a generator of the ideal of 
${\bf C} \left[ q_{ij} \ (1 \leq i < j \leq 4) \right]$ 
which is the kernel of the ring homomorphism given by 
$$
q_{ij} \mapsto 
\frac{(x_{i} - x_{j})(x_{-i} - x_{-j})}{(x_{i} - x_{-j})(x_{-i} - x_{j})} 
\in A_{0}. 
$$ 
%\vspace{1ex}

\noindent
\underline{\bf Problem.} 
Let $J'$ be a primitive modular form obtained from $J$ by dividing 
the GCM (greatest common divisor) of its Fourier coefficients. 
Then for each prime $p,$ 
\begin{eqnarray*}
& & \mbox{the closed subset of ${\cal A}_{4} \otimes_{\bf Z} {\bf F}_{p}$ 
defined by $J' \ {\rm mod}(p) = 0$} 
\\ 
& \mbox{\boldmath $\stackrel{\mbox{\Large {\bf ?}}}{=}$} & 
\mbox{the Zariski closure of $\tau({\cal M}_{4} \otimes_{\bf Z} {\bf F}_{p})$ 
in ${\cal A}_{4} \otimes_{\bf Z} {\bf F}_{p}.$}
\end{eqnarray*}
%\vspace{1ex}

\noindent
\underline{\bf Hyperelliptic Schottky problem.} ([I4]) 
Let $p_{ij}$ be the universal periods given in Example 4.1. 
Then 
$$
p'_{ij} \stackrel{\rm def}{=} p_{ij}|_{x_{-k} = - x_{k}} \ \ (1 \leq k \leq g) 
$$
become the multiplicative periods of the hyperelliptic curve $C_{\rm hyp}$ 
over 
$$ 
{\bf Z} 
\left[ \frac{1}{2 x_{i}}, \ \frac{1}{x_{i} \pm x_{j}} (i \neq j) \right] 
[[ y_{1},..., y_{g} ]] 
$$
uniformized by the Schottky group: 
$$
\left\langle \left. 
\left( \begin{array}{cc} x_{k} & -x_{k} \\ 1 & 1 \end{array} \right) 
\left( \begin{array}{cc} 1   & 0  \\  0  & y_{k} \end{array} \right) 
\left( \begin{array}{cc} x_{k} & -x_{k} \\ 1 & 1 \end{array} \right)^{-1} 
\ \right| \ k = 1,..., g \right\rangle. 
$$
Since $C_{\rm hyp}$ is generic in the moduli space of hyperelliptic curves, 
for any Siegel modular form $\varphi$ over a field of characteristic $\neq 2,$ 
$$
\mbox{$\varphi$ vanishes on the locus of hyperelliptic Jacobians} 
\ \Longleftrightarrow \ 
F(\varphi)|_{q_{ij} = p'_{ij}} = 0.
$$
%\vspace{1ex}

\noindent
\underline{\bf Problem.}
Give an explicit lower bound of $n(g) \in {\bf N}$ satisfying that 
$$
\mbox{$\varphi$ vanishes on the locus of hyperelliptic Jacobians} 
\ \Longleftrightarrow \ 
F(\varphi)|_{q_{ij} = p'_{ij}} \in I^{n(g)}, 
$$
where $I$ is the ideal generated by $y_{1},..., y_{g}.$ 
\vspace{2ex}

\noindent
\underline{\bf Theta constants and ring structure.}
\vspace{2ex}

For $g \geq 2,$ let 
$$
\theta_{g}(Z) \ \stackrel{\rm def}{=} \ 
\prod_{\footnotesize{\begin{array}{l}  
\mbox{\boldmath $a$}, \mbox{\boldmath $b$}\in \{ 0, 1/2 \}^{g} \\ 
4 \mbox{\boldmath $a$}^{t} \mbox{\boldmath $b$}: {\rm even} 
\end{array}}} 
\sum_{\footnotesize{\mbox{\boldmath $n$} \in {\bf Z}^{\it g}}} 
\exp \left( 2 \pi \sqrt{-1} \left[ \frac{1}{2} 
(\mbox{\boldmath $n$} + \mbox{\boldmath $a$}) 
Z^{t}(\mbox{\boldmath $n$} + \mbox{\boldmath $a$}) 
+ (\mbox{\boldmath $n$} + \mbox{\boldmath $a$})^{t}\mbox{\boldmath $b$} 
\right] \right)
$$
be the product of even {\bf theta constants} of degree $g.$ 
If $g \geq 3,$ then $\theta_{g}$ is an integral Siegel modular form 
of degree $g$ and weight $2^{g-2}(2^{g}+1).$ 
\vspace{2ex}

\noindent
\underline{\bf Theorem 6.3.} ([I2, 3]). 
\begin{it} For $g \geq 3,$ 

{\rm (1)} $T_{g,h}({\bf Z})$ is a free ${\bf Z}$-module of finite rank 
satisfying that $T_{g,h}({\bf Z}) \otimes_{\bf Z} {\bf C} = T_{g,h}({\bf C})$, 
and that $T_{g,0}({\bf Z}) = {\bf Z},$ 
$T_{g,h}({\bf Z}) = \{ 0 \}$ if $h < 0.$ 
Furthermore, the ring of integral Teichm\"{u}ller modular forms of degree $g:$ 
$$
T_{g}^{*}({\bf Z}) \stackrel{\rm def}{=} 
\bigoplus_{h \geq 0} T_{g,h}({\bf Z}) 
$$ 
becomes a normal ring which is finitely generated over ${\bf Z}$. 

{\rm (2)} For the product $\theta_{g}$ of even theta constants of degree $g$, 
$$
N_{g} \ \stackrel{\rm def}{=} \ 
\left\{ \begin{array}{ll}
-2^{28} & (g = 3), \\ 2^{2^{g-1}(2^{g} - 1)} & (g \geq 4). 
\end{array} \right.
$$
Then $\sqrt{\tau^{*}(\theta_{g}) / N_{g}}$ is a {\it \bfseries primitive} 
element of 
$T_{g, 2^{g-3}(2^{g} + 1)}({\bf Z})$, i.e., 
not congruent to $0$ modulo any prime. 

{\rm (3)} $T_{3}^{*}({\bf Z})$ is generated by Siegel modular forms 
over ${\bf Z}$ and by $\sqrt{\tau^{*}(\theta_{3}) / N_{3}}$ 
which is of weight $9$, hence is not a Siegel modular form. 
\end{it}
\vspace{2ex}

\underline{\it Proof.}
(1) follows from the argument in the proof of Theorem 6.1 (2) replacing 
$$
\left( {\cal A}_{g}, \ \overline{\cal A}_{g}, \ \overline{\lambda} \right) 
\ \ \mbox{by} \ \ 
\left( {\cal M}_{g}, \ \overline{\cal M}_{g}, \ 
\bigwedge^{g} {\pi}_{*} (\omega_{{\cal C} / \overline{\cal M}_{g}}) \right), 
$$ 
where $\pi : {\cal C} \rightarrow \overline{\cal M}_{g}$ 
denotes the universal stable curve over Deligne-Mumford's compactification. 
$\kappa_{\Delta}$ is used to show that 
any integral Teichm\"{u}ller modular form can be extended to 
a global section on $\overline{\cal M}_{g}.$ 

(2) Let $D$ be the divisor of ${\cal M}_{g} \otimes_{\bf Z} \overline{\bf Q}$ 
consisting of curves $C$ which have a line bundle $L$ such that 
$L^{\otimes 2} \cong \Omega_{C}$ and that 
$\dim H^{0}(C, L)$ is positive and even. 
Then as is shown in [Ty2], $2 D$ gives the divisor of 
$\tau^{\ast}(\theta_{g}),$ and hence a Teichm\"{u}ller modular form of weight 
$(\mbox{the weight of $\theta_{g}$})/2$ with divisor $D,$ 
which exists and is uniquely determined up to constant, 
is a root of $\tau^{\ast}(\theta_{g})$ up to constant (see below). 
Since $D$ is defined over $\overline{\bf Q},$ 
a root of $\tau^{\ast}(\theta_{g})$ times a certain number 
is defined and primitive over ${\bf Z}.$ 
To determine this number, $\kappa_{\Delta}$ is used as follows: 
Let $A_{0}, A_{\Delta}, p_{ij}$ be as in Example 4.1. 
Then 
$$
\theta_{g}(Z) \ = \ 
2^{2^{g-1}(2^{g}-1)} 
\left( \prod_{\scriptsize{\begin{array}{c} 
(b_{1},...,b_{g}) \in \{ 0, 1/2 \}^{g} \\ 
\sum_{i} b_{i} \in {\bf Z} 
\end{array}}} 
(-1)^{\sum_{i} b_{i}} \right) P \cdot \alpha^{2}, 
$$ 
where 
\begin{eqnarray*}
\alpha & : & 
\mbox{a primitive element of 
${\bf Z} \left[ q_{ij}^{\pm 1} \ (i \neq j) \right] 
\left[ \left[ q_{11},...,q_{gg} \right] \right],$} 
\\ 
P & = & 
\prod_{\scriptsize{\begin{array}{c} 
(b_{1},...,b_{g}) \in \{ 0, 1/2 \}^{g} \\ 
\sum_{i} b_{i} \in {\bf Z} 
\end{array}}} 
\frac{1}{2} 
\sum_{S \subset \{ 1,...,g \}} (-1)^{\sharp \{k \in S | b_{k} \neq 0 \}} 
\prod_{i \in S, j \not\in S} q_{ij}^{-1/2} 
\\
& \Rightarrow & 
\left. \left( \mbox{the constant term of 
$P|_{q_{ij} = p_{ij}} \in A_{\Delta}$} \right) 
\right|_{x_{1} = x_{-2},...,x_{g} = x_{-1}} = 1. 
\end{eqnarray*}
Hence we have (see Exercise 15 below): 
\begin{eqnarray*}
& & \sqrt{ \mbox{the constant term of $P|_{q_{ij} = p_{ij}}$}} \in A_{0} 
\\
& \Rightarrow & \sqrt{\theta_{g}|_{q_{ij} = p_{ij}}} \in 
\left\{ \begin{array}{ll}
\sqrt{-1} \cdot 2^{27} \cdot A_{\Delta} & (g = 3), \\ 
2^{2^{g-1}(2^{g}-1) - 1} \cdot A_{\Delta} & (g \geq 4). 
\end{array} \right.
\end{eqnarray*}

(3) Recall the result of Igusa [Ig2] that the ideal of 
$S_{3}^{\ast}({\bf C})$ vanishing on the hyperelliptic locus is 
generated by $\theta_{3}.$ 
Since the Torelli map ${\cal M}_{3} \rightarrow {\cal A}_{3}$ is dominant 
and of degree $2,$ if we denote $\iota$ by the multiplication by $-1$ 
on abelian varieties, then 
\begin{eqnarray*}
\bigoplus_{h: \ {\rm even}} T_{3,h}({\bf C}) 
& = & 
\left\{ \left. f \in T_{3}^{\ast}({\bf C}) \ \right| \ \iota(f) = f \right\} 
= S_{3}^{\ast}({\bf C}),  
\\
\bigoplus_{h: \ {\rm odd}} T_{3,h}({\bf C}) 
& = & 
\left\{ \left. f \in T_{3}^{\ast}({\bf C}) \ \right| \ \iota(f) = -f \right\}. 
\end{eqnarray*}
Let $f$ have odd weight. Then by $\iota(f) = -f,$ 
$f = 0$ on the hyperelliptic locus, and hence by Igusa's result, 
$f^{2} / \theta_{3}$ becomes a Siegel modular form. 
Therefore, $T_{3}^{\ast}({\bf C})$ is generated by $S_{3}^{\ast}({\bf C})$ 
and $\sqrt{\tau(\theta_{3})}$ which implies (3) 
because $\sqrt{\tau(\theta_{3})/N_{3}}$ is integral and primitive. QED. 
\vspace{2ex}

\noindent
\underline{\bf Exercise 15.} 
Prove that 
$$
\left( \prod_{\scriptsize{\begin{array}{c} 
(b_{1},...,b_{g}) \in \{ 0, 1/2 \}^{g} \\ 
\sum_{i} b_{i} \in {\bf Z} 
\end{array}}} 
(-1)^{\sum_{i} b_{i}} \right) \ = \ 
\left\{ \begin{array}{ll} 1 & (g = 3), \\ -1 & (g \geq 4). \end{array} 
\right. 
$$
\vspace{1ex}

\noindent
\underline{\bf TMFs of degree 2.} 
Let $k$ be an algebraically closed field $k$ of characteristic $\neq 2.$ 
Then any proper smooth curve $C$ of genus $2$ over $k$ is hyperelliptic, 
more precisely a base of $H^{0}(C, \Omega_{C})$ gives rise to 
a morphism $C \rightarrow {\mathbb P}^{1}_{k}$ of degree $2$ 
ramified at $6$ points, and hence 
\begin{eqnarray*}
{\cal M}_{2} \otimes_{\bf Z} k & \cong & 
\left. \left\{ \left. (x_{1}, x_{2}, x_{3} \in {\mathbb P}^{1}_{k} - 
\{ 0, 1, \infty \} \ \right| \ x_{i} \neq x_{j} \ (i \neq j) \right\} 
\right/ S_{6}, 
\end{eqnarray*}
where each element $\sigma$ of the symmetric group $S_{6}$ degree $6$ 
acts on $(x_{1}, x_{2}, x_{3})$'s such as 
$$
(\sigma(x_{1}), \sigma(x_{2}), \sigma(x_{3}), 0, 1, \infty)
$$ 
is obtained from $\sigma(x_{1}, x_{2}, x_{3}, 0, 1, \infty)$ 
by some M\"{o}bius transformation of $GL_{2}(k).$ 
Therefore, ${\cal M}_{2} \otimes_{\bf Z} k$ becomes an affine variety, 
and $T_{2,h}(k) = H^{0}({\cal M}_{2}, \lambda^{\otimes h} \otimes_{\bf Z} k)$ 
is infinite dimensional. 
In fact, it is proved in [I3] that the ring 
$$
T_{2}^{\ast}({\bf Z}) \stackrel{\rm def}{=} 
\bigoplus_{h \in {\bf Z}} T_{2, h}({\bf Z})
$$
of integral Teichm\"{u}ller modular forms is generated 
by $\tau^{\ast}(S_{2}^{\ast}({\bf Z}))$ and by 
$2^{12}/\left( \tau^{\ast}(\theta_{2}) \right)^{2}$ 
which is of weight $-10.$   
\vspace{2ex}

\noindent
\underline{\bf Construction of TMFs.} 
Assume that $g \geq 3.$ 
Then by results of Mumford [Mu1] and Harer [H1], 
the {\bf Picard group} of ${\cal M}_{g}:$ 
\begin{eqnarray*}
{\rm Pic}({\cal M}_{g}) 
& \stackrel{\rm def}{=} & 
\mbox{the group of linear equivalence classes of 
line bundles on ${\cal M}_{g}.$} 
\end{eqnarray*}
is isomorphic to 
$H^{2}({\cal M}_{g}({\bf C}), {\bf Z}) \cong H^{2}(\Pi_{g}, {\bf Z})$ 
$(\Pi_{g}$ denotes the Teichm\"{u}ller modular group of degree $g),$ 
and this is free of rank $1$ generated by the Hodge line bundle $\lambda$. 
Therefore, 
\begin{eqnarray*}
& & \mbox{$D \neq 0$ is an effective divisor on ${\cal M}_{g}$ 
over a subfield $K$ of ${\bf C}$} 
\\
& \Rightarrow & 
\mbox{there are $h \in {\bf N}$ such that 
${\cal O}_{{\cal M}_{g}}(D) \cong \lambda^{\otimes h}$} 
\\ 
& \Rightarrow & 
\mbox{there is $f \in T_{g,h}(K)$ such that $(f) = D$} 
\\
& & \mbox{(for the application, see the proof of Theorem 6.3 (2)),}
\end{eqnarray*}
and $f$ is uniquely determined (up to a nonzero constant) 
by the existence of the Satake-type compactification of ${\cal M}_{g}.$ 
From this method, one can construct Teichm\"{u}ller modular forms 
and study their rationality using $\kappa_{\Delta}.$
\vspace{2ex}

\noindent
\underline{\bf Mumford's isomorphism.} 
We recall Mumford's isomorphism (for $g > 1$) given in 5.3:  
$$
\det \left( \pi_{*} \left( \Omega_{{\cal C}/{\cal M}_{g}}^{n} \right) \right)  
\cong 
\det \left( \pi_{*} \left( \Omega_{{\cal C}/{\cal M}_{g}} 
\right) \right)^{\otimes (6n^{2} - 6n +1)} 
= \lambda^{\otimes (6n^{2} - 6n +1)}. 
$$
Therefore, by putting $n = 2$, 
$$   
\bigwedge^{3g-3} \pi_{*} \left( \Omega_{{\cal C}/\overline{\cal M}_{g}}^{2} \right) 
\cong \lambda^{\otimes 13}. 
$$ 
In order to express this isomorphism explicitly, 
we consider the morphism 
$$ 
\rho_{g} : {\rm Sym}^{2} \left( \pi_{*} 
\left( \Omega_{{\cal C}/{\cal M}_{g}} \right) \right) 
\ni (s, s') \mapsto s \cdot s' \in 
\pi_{*} \left( \Omega_{{\cal C}/{\cal M}_{g}}^{\otimes 2} \right) 
$$ 
between vector bundles on ${\cal M}_{g}.$ 
\begin{itemize}

\item 
If $g = 1,$ then $\rho_{1}$ is an isomorphism and gives 
$$
\lambda^{\otimes 2} \stackrel{\rho_{1}}{\longrightarrow} 
\pi_{*} \left( \Omega_{{\cal C}/{\cal M}_{2}}^{\otimes 2} \right)  
\cong \lambda^{\otimes 14} 
\ \Rightarrow \ 
{\cal O}_{{\cal M}_{1}} 
\ni 1 \mapsto \pm \Delta(\tau) = q \prod_{n=1}^{\infty} (1 - q^{n})^{24} 
\in T_{1,12}({\bf Z}). 
$$

\item 
If $g = 2,$ then $\rho_{2}$ is an isomorphism and gives 
$$
\lambda^{\otimes 3} \stackrel{\det(\rho_{2})}{\longrightarrow} 
\bigwedge^{3} \pi_{*} 
\left( \Omega_{{\cal C}/{\cal M}_{2}}^{\otimes 2} \right)  
\cong \lambda^{\otimes 13} 
\ \Rightarrow \ 
{\cal O}_{{\cal M}_{2}} 
\ni 1 \mapsto \pm \left( \tau^{*}(\theta_{2}) / 2^{6} \right)^{2} 
\in T_{2,10}({\bf Z}). 
$$

\item 
If $g = 3,$ then $\rho_{3}$ is an isomorphism generically 
and vanishes on the hyperelliptic locus, hence this gives 
$$
\lambda^{\otimes 4} \stackrel{\det(\rho_{3})}{\longrightarrow} 
\bigwedge^{6} 
\pi_{*} \left( \Omega_{{\cal C}/{\cal M}_{3}}^{\otimes 2} \right)  
\cong \lambda^{\otimes 13} 
\ \Rightarrow \ 
{\cal O}_{{\cal M}_{3}} \ni 1 \mapsto \pm \sqrt{\tau^{*}(\theta_{3}) / N_{3}} 
\in T_{3,9}({\bf Z}). 
$$

\end{itemize}
Furthermore, 
[I6] showed the explicit formula of the Mumford isomorphism for any $g$ 
as an infinite product extending $\Delta(\tau)$.

\newpage

\begin{center}
\underline{\large {\bf References}} 
\end{center}
\begin{itemize}

\item[{[BG]}] 
Brinkmann and Gerritzen, 
The lowest term of the Schottky modular form, 
Math. Ann. {\bf 292} (1992) 329--335. 

\item[{[D]}] 
P. Deligne, 
Formes modulaires et repr\'{e}sentations $l$-adiques, 
Lecture Notes in Math. {\bf 179} (1971) p. 139--172. 

\item[{[DM]}] 
P. Deligne and D. Mumford, 
The irreducibility of the space of curves of given genus, 
Inst. Hautes \'{E}tudes Sci. Publ. Math. {\bf 36} (1969) 75--109. 

\item[{[DS]}] 
P. Deligne and J. P. Serre, 
Formes modulaires de poids 1, 
Ann. Sci. Ecole Norm. Sup. {\bf 7} (1974) 507--530. 

\item[{[FaC]}] 
G. Faltings and C. L. Chai, 
Degeneration of Abelian Varieties, 
Ergeb. Math. Grenzgeb., 3. Folge, Bd. 22, 
Springer-Verlag, 1990. 

\item[{[Fay]}] 
J. D. Fay, 
Theta Functions on Riemann Surfaces, 
Lecture Notes in Math., Vol. 352, 
Springer-Verlag, 1973. 

\item[{[FKM]}] 
J. Fogarty, F. Kirwan and D. Mumford, 
Geometric Invariant Theory, 
Ergeb. Math. Grenzgeb., 3. Folge, Bd. 34, 
Springer-Verlag, 1994. 

\item[{[G]}] 
L. Gerritzen, 
Equations defining the periods of totally degenerate curves, 
Israel J. Math. {\bf 77} (1992) 187--210.

\item[{[GP]}] 
L. Gerritzen and M. van der Put, 
Schottky Groups and Mumford curves, 
Lecture Notes in Math. Vol. 817, Springer-Verlag, 1980. 

\item[{[H1]}] 
J. Harer, 
The second homology group of the mapping class group of an orientable surface, 
Invent. Math. {\bf 72} (1983) 221--239.

\item[{[H2]}] 
J. Harer, 
Stability of the homology of mapping class groups of orientable surfaces, 
Ann. of Math. {\bf 121} (1985) 215--251. 

\item[{[HM]}] 
J. Harris and I. Morrison, 
Moduli of Curves, 
Graduate Texts in Math., Vol. 187, Springer-Verlag, 1998. 

\item[{[I1]}] 
T. Ichikawa, 
$P$-adic theta functions and solutions of the KP hierarchy, 
Commun. Math. Phys. {\bf 176} (1996) 383--399. 

\item[{[I2]}] 
T. Ichikawa, 
Theta constants and Teichm\"{u}ller modular forms, 
J. Number Theory {\bf 61} (1996) 409--419. 

\item[{[I3]}] 
T. Ichikawa, 
Generalized Tate curve and integral Teichm\"{u}ller modular forms, 
Amer. J. Math. {\bf 122} (2000) 1139--1174. 

\item[{[I4]}] 
T. Ichikawa, 
Universal periods of hyperelliptic curves and their applications, 
J. Pure Appl. Algebra {\bf 163} (2001) 277--288. 

\item[{[I5]}] 
T. Ichikawa, 
Siegel modular forms of degree 2 over rings, 
J. Number Theory {\bf 129} (2009) 818--823. 

\item[{[I6]}] 
T. Ichikawa, 
A product formula of Mumford forms, 
and the rationality of Ruelle zeta values for Schottky groups, 
arXiv:1058478v1. 

\item[{[Ig1]}] 
J. Igusa, 
On Siegel modular forms of genus two, 
Amer. J. Math. {\bf 84} (1962) 175--200; 
II, Amer. J. Math. {\bf 86} (1964) 392--412.  

\item[{[Ig2]}] 
J. Igusa, 
Modular forms and projective invariants, 
Amer. J. Math. {\bf 89} (1967) 817--855. 

\item[{[Ig3]}] 
J. Igusa, 
On the ring of modular forms of genus two over ${\bf Z},$ 
Amer. J. Math. {\bf 101} (1979) 149--193. 

\item[{[IhN]}] 
Y. Ihara and H. Nakamura, 
On deformation of maximally degenerate stable marked curves and Oda's problem, 
J. Reine Angew. Math. {\bf 487} (1997) 125--151. 

\item[{[IT]}] 
Y. Imayoshi and M. Taniguchi, 
An Introduction to Teichm\"{u}ller Spaces, 
Springer-Verlag, 1992. 

\item[{[MadW]}] 
I. Madsen and M. Weiss, 
The stable moduli space of Riemann surfaces: Mumford's conjecture, 
Ann. of Math. {\bf 165} (2007) 843--941. 

\item[{[ManD]}] 
Yu. I. Manin and V. G. Drinfeld, 
Periods of $p$-adic Schottky groups, 
J. Reine Angew. Math. {\bf 262/263} (1973) 239--247. 

\item[{[MaT]}] 
G. Mason and M. P. Tuite, 
On genus two Riemann surfaces formed from sewn tori, 
Commun. Math. Phys. {\bf 270} (2007) 587--634. 

\item[{[Mi]}] 
E. Y. Miller, 
The homology of the mapping class group, 
J. Differ. Geom. {\bf 24} (1986) 1--14. 

\item[{[Mo]}] 
S. Morita, 
Characteristic classes of surface bundles, 
Invent. Math. {\bf 90} (1987) 551--577. 

\item[{[Mu1]}] 
D. Mumford, 
Abelian quotients of the Teichm\"{u}ller modular group, 
J. Analyse Math. {\bf 18} (1967) 227--244. 

\item[{[Mu2]}] 
D. Mumford, 
An analytic construction of degenerating curves over complete local rings, 
Compositio Math. {\bf 24} (1972) 129--174. 

\item[{[Mu3]}] 
D. Mumford, 
An analytic construction of degenerating abelian varieties 
over complete rings, 
Compositio Math. {\bf 24} (1972) 239--272. 

\item[{[Mu4]}] 
D. Mumford, 
Stability of projective varieties, 
L'Ens. Math. {\bf 23} (1977) 39--110. 

\item[{[Mu5]}] 
D. Mumford, 
Towards an enumerative geometry of the moduli space of curves, 
In Arithmetic and Geometry, II, 
Progress in Math., Vol. 36, 1983, pp. 271--326. 

\item[{[MumSW]}] 
D. Mumford, C. Series and D. Wright, 
Indra's pearls, Cambridge University Press (2002). 

\item[{[Mur]}] 
V. K. Murty, 
Intorduction to Abelian Varieties, 
CRM Monograph Series, Vol. 3, Amer. Math. Soc., 1993. 

\item[{[N]}] 
I. Nakamura, 
Stability of degenerate abelian varieties, 
Invent. Math. {\bf 136} (1999) 659--715. 

\item[{[S]}] 
F. Schottky, 
\"{U}ber eine specielle Function, welche bei einer bestimmten 
linearen Transformation ihres Arguments univer\"{a}ndert bleibt, 
J. Reine Angew. Math. {\bf 101} (1887) 227--272. 

\item[{[Se]}] 
J. P. Serre, 
Modular forms of weight one and Galois representations, 
in Algebraic Number Fields (ed. A. Fr\"{o}hlich), 
Acad. Press, 1977, pp. 193--268. 

\item[{[Si]}] 
J. H. Silverman, 
Advanced Topics in the Arithmetic of Elliptic Curves, 
Graduate Texts in Math., Vol. 151, Springer-Verlag, 1994. 

\item[{[T]}] 
J. Tate, 
A Review of Non-Archimedean Elliptic Functions, 
in : J. Coates and S. T. Yau, (eds.), 
Elliptic Curves, Modular Forms, \& Fermat's Last Theorem, 
Series in Number Theory, Vol. 1, 1995, pp. 162--184. 

\item[{[Ty1]}] 
S. Tsuyumine, 
On Siegel modular forms of degree three, 
Amer. J. Math. {\bf 108} (1986) 755--862. 

\item[{[Ty2]}] 
S. Tsuyumine, 
Thetanullwerte on a moduli space of curves and hyperelliptic loci, 
Math. Z. {\bf 207} (1991) 539--568. 

\end{itemize}

\newpage

\begin{center}
\underline{\large {\bf List of Exercises}} 
\end{center}
\begin{tabular}{llll}
Exercise 1: p.5,    & Exercise 2: p.11,   & Exercise 3: p.13,  & Exercise 4: p.14, 
\\ 
Exercise 5: p.15,  & Exercise 6: p.20,   & Exercise 7: p.20,  & Exercise 8: p.23, 
\\ 
Exercise 9: p.29,  & Exercise 10: p.37, & Exercise 11: p.37, & Exercise 12: p.38, 
\\ 
Exercise 13: p.40, & Exercise 14: p.43, & Exercise 15: p.47. &  
\end{tabular}
\vspace{5ex}

\begin{center}
Takashi Ichikawa \\
Department of Mathematics, 
Graduate School of Science and Engineering \\
Saga University, Saga 840-8502, Japan \\
e-mail: ichikawa@ms.saga-u.ac.jp, 
ichikawn@cc.saga-u.ac.jp 
\end{center}

\end{document}